\PassOptionsToPackage{backend=bibtex}{biblatex} 
\documentclass{zml-article}    
\usepackage[hyphens]{xurl} 
\usepackage{xurl} 
\usepackage{hyperref}
\usepackage[noabbrev,capitalise]{cleveref}  
\usepackage{comment}
\usepackage{booktabs}
\MFUnocap{25th}        

\title{Finite sets, mappings, cardinals, and \protect\\ arithmetic in intuitionistic New Foundations}  
\shorttitle{Finite sets, mappings, cardinals, and  arithmetic in intuitionistic New Foundations}  
\author[1]{Michael Beeson} 
\affil[1]{San Jos\'e State University (Professor Emeritus) and UCSC (Research Associate)}
\corrauthor{Michael Beeson}{profbeeson@gmail.com}
\date{\today}          

\newtheorem{theorem}{Theorem}
\newtheorem{lemma}[theorem]{Lemma}
\newtheorem{corollary}[theorem]{Corollary}
\newtheorem{definition}[theorem]{Definition}
\numberwithin{theorem}{section}  
\numberwithin{lemma}{section} 
\numberwithin{corollary}{section}

\let\sect\S
\def\S{{\mathbb S}}
\def\H{{\mathbb H}}

\def\T{{\mathbb T}}
\def\J{{\mathbb J}}

\def\imp{\ {\rightarrow}\ }
\def\imp{\ {\rightarrow}\ }
\def\FUNC{{\sf FUNC}}
\def\V{{\mathbb V}}

\def\USC{\mathcal P_1}
\def\SSC{\mathcal P_s}
\def\SC{\mathcal P}
\def\FregeN{{\mathbb F}}
\def\F{{\mathbb F}}   

\def\INF{\ensuremath{\mathit{i}\mkern0.4mu\mathsf{NF}}}
\newcommand{\NF}{\ensuremath{\mathsf{NF}}}

\def\NC{{\sf NC}}
\def\Nc#1{\lvert #1\rvert}
\def\FINITE{{\sf FINITE}}
\def\DECIDABLE{{\sf DECIDABLE}}

\def\m{{\bf m}}
\def\zero{{\sf zero}}
\def\one{{\sf one}} 
\def\two{{\sf two}}
\def\three{{\sf three}}
\def\four{{\sf four}}

\def\iff{{\ \leftrightarrow \ }}
\def\twocases#1#2{\left\{
\begin{array}{ll} 
#1\\ 
{}\\
#2
\end{array} \right.}  

\abstract{
\NF\ set theory using intuitionistic logic is called \INF. 
We develop the theories of finite sets and their power sets and mappings, 
finite cardinals and their ordering, cardinal exponentiation, addition, and 
multiplication.  We follow Rosser and Specker with appropriate 
constructive modifications, especially replacing ``arbitrary subset'' 
by ``separable subset'' in the definitions of exponentiation and order. 
It is not known whether \INF\ proves that the set of finite cardinals is 
infinite, so the whole development must allow for the possibility that 
there is a maximum integer;  arithmetical computations might ``overflow''
as in a computer or odometer, and theorems about them must be carefully
stated to allow for this possibility.  The work presented here is 
intended as a basis for further investigations of \INF, including
the development of Bishop-style constructive mathematics in \INF. 
}
\begin{document}
\maketitle

\section{Introduction}
Quine's \NF\ set theory is a first-order theory whose language contains only the binary predicate symbol $\in$,
and whose axioms are extensionality and stratified comprehension.  
The definition of these axioms will be reviewed below; full details 
can be found in~\cite{rosser1978}.
Intuitionistic \NF, or \INF,  is the theory
with the same language and axioms as \NF, but with intuitionistic logic instead of classical.%
\footnote{ 
$\INF$ is briefly mentioned in Forster's thesis~\cite{Forster1976-NF};
next mentioned in~\cite{Dzierzgowski1992}, 
~\cite{Dzierzgowski1995}, and~\cite{Dzierzgowski1996}, 
where the focus is on intuitionistic type theories.
The initial development of $\INF$ may be in~\cite{forster2009},
which first called attention to the problem of interpreting HA in $\INF$.
}
Here we intend to provide a coherent infrastructure of definitions,
theorems, and Lean-checked proofs on which further investigations
can be based.%
\footnote{
The development of \NF\ and its variants has been surveyed by Forster~\cite{forster}, 
and a comprehensive online bibliography of research on set theories with a 
universal set is maintained by Holmes~\cite{HolmesNFBiblio}.
}

The ``axiom'' of infinity is a theorem of \NF, proved by Rosser 
~\cite{rosser1952, rosser1978} and Specker~\cite{specker1953}.
 These proofs use classical logic
in an apparently essential way.  It is still an open question 
whether \INF\ proves the existence of an infinite set.   
 The Stanford
Encyclopedia of Philosophy article on \NF\ says~\cite{sep-quine-nf}

\begin{quote}
 {\em The only known proof (Specker's) of the axiom of infinity in \NF\ has too little constructive content to allow a demonstration that \INF\ 
  admits an implementation of Heyting arithmetic.}
\end{quote}

\noindent 
In attempting to determine whether the quoted statement is true, I found that 
I first needed to develop enough basic mathematics in $\INF$ to tackle 
Specker's proof.%
\footnote{Eventually I came to the conclusion that the 
quoted statement is true; that is beyond the scope of this paper,
but see a short discussion near the end.}  
That mathematical
infrastructure  is presented in this paper.
The purpose
of this development is to provide a basis on which one can:
\begin{itemize}
\item investigate \INF\ further;
\item develop Bishop-style constructive mathematics in \INF\\
(after proving or assuming infinity).
\end{itemize} 
The following questions about \INF\ remain open:%
\footnote{The terms used in these questions can be looked up in the 
index of~\cite{beeson1985a}.}
\begin{itemize}
\item Is the set $\F$ of finite cardinals finite?  Is it infinite?
\item Can one point to any specific instance of the law of the excluded middle
that is 
not provable in 
\INF? 
\item Is there any   
double-negation interpretation from \NF\ to \INF?
\item Is Church's thesis consistent with \INF?  Markov's principle? 
\item Is \INF\ closed under Church's rule?  
\end{itemize}

Regarding whether $\F$ is finite:  For all we know,
  there might be a largest finite cardinal $\m$, which would 
contain a finite set $U$ that is ``unenlargeable'', in the sense 
that we cannot find any $x$ that is not a member of $U$.  Classically,
that would imply $U = \V$, which is a contradiction, since $\V$ is   
not finite.  But intuitionistically, it is an open question. 

Each 
of the lemmas and theorems in this paper is provable in \INF.  
 An important reference for \NF\ is Rosser's book~\cite{rosser1953,rosser1978}.%
\footnote{
The two editions are identical except for the Appendices added to the second
edition, one of which contains Rosser's proof of infinity.}
But the
logical apparatus of Rosser's system includes a Hilbert-style 
epsilon-operator, which is not compatible with an intuitionistic version,
and also, we do not wish to assume the axiom of infinity.
 Since all of Rosser's results are obtained using
classical logic, we cannot rely on Rosser.  

It should be noted that the consistency of classical \NF\ has been
proved~\cite{holmes2025}, and the proof has been checked in Lean.  This
result implies, of course, that the subtheory \INF\ considered here is 
also consistent; but it is otherwise not directly relevant. 

{\em Notational issues}. There is no traditional, universally accepted notation for some of the notions
central to NF.  Rosser~\cite{rosser1978} and Specker~\cite{specker1953} are two 
of the original sources.  Both of these were written prior to the advent of 
\TeX\ and \LaTeX, and for the most part were limited to characters found on a 
typewriter keyboard.  Forster~\cite{forster} used different notation, making
use of \TeX.  Now, however, keyboard characters are back in style,
because they are easier to use in computer proof-checking.
Lean, for example, goes to great lengths to support complicated
typography---but one cannot search for those symbols,  which is quite
annoying.

I therefore used Specker's notation when using Lean and 
in pre-publication versions of this paper; but the referee asked me to 
change it, so I did.  The following table compares the notational styles of Specker
and this paper.  It may prove useful if anyone wants to compare this paper 
to Specker or Rosser, or to the Lean proofs, or to readers who know one style
or the other already.

\begin{table}[htbp]
\caption{Notation}
\centering
\begin{tabular}{@{}ll@{}}   
Specker   &  This paper \\
\midrule
$USC(x)$  & $\USC(x)$ \\
$SSC(x)$  & $\SSC(x)$ \\
$SC(x)$   & $\SC(x)$ \\
$Nc(x)$   & $\Nc{x}$ \\
$\Lambda$ & $\emptyset$ \\
\end{tabular}
\label{default}
\end{table}

{\em Use of computer proof-checking}. 
All the proofs in this paper have been checked in the proof assistant Lean.
Could there still be errors? The possible sources of error are
\begin{itemize}
\item Use of an unstratified definition
\item Lean proof and paper proof might not exactly correspond
\item Lean might have smuggled in classical logic, i.e., used it without telling me.
\item Perhaps the order of theorems in the paper is not strictly the logical order.
\end{itemize}
Regarding the smuggling: Lean's underlying theory is intuitionistic, but 
the library is classical, and even though I didn't use the library, and 
even though Lean experts helped me, the possibility theoretically exists.
Regarding the correspondence: if there are such problems, they are just typos.
Regarding stratification: I allowed the full comprehension axiom, but 
used only stratified instances.  I used a 
computer script {\em ex post facto} to check stratification.%
\footnote{Originally I intended to use a finite axiomatization.  But it 
is often quite complicated to derive simple definitions from a finite 
axiomatization; and then one still has to worry if the finite axiomatization
is really correct.}
Regarding the order of theorems: At least no lemma or theorem is cited before it is proved.
Of course in Lean, the logical order is enforced, but that is often not the 
best order for presentation.
\smallskip

The reader who is worried about errors in Lean has the option to forget 
it was ever mentioned, 
and just read the proofs, which are here presented in complete human-readable
detail. 
\smallskip 

{\em Acknowledgements}. Thanks to Thomas Forster for asking me (once a year for twenty years)
about the strength of \INF.  
Thanks to Randall Holmes,  Albert Visser, and Thomas Forster for 
many emails on this subject.   Thanks to the creators of the proof assistant Lean~\cite{LeanSystem},  which 
has enabled me to state with high confidence that there are no errors in this 
paper.
Thanks to the users of Lean who helped me acquire sufficient 
expertise in using Lean by answering my questions, especially Mario Carneiro.

\section{Axioms of NF, ordered pairs, and functions} \label{section:axioms} 
The axioms of \NF\ are extensionality and stratified comprehension.
The axiom of extensionality says that two sets with the same
elements are equal.   The axiom schema of stratified comprehension says 
that $\{ x : \phi(x) \}$ exists, if $\phi$ is a stratified formula.
A formula is {\bf stratified}, or {\bf stratifiable}, if each of 
its variables (both bound and free) can be assigned a non-negative 
integer (``index'' or ``type'') such that (i) in every subformula $x \in y$,
$y$ gets an index one greater than $x$ gets, and (ii) every occurrence
of each variable gets the same index.  

Thus, the ``universe'' $\V$ can be defined as 
$$ \V = \{x : x=x\}$$
but the Russell set $\{x: x \not \in x\}$ cannot be defined.
The empty set can be defined as     
$$ \emptyset = \{x:x \neq x\}.$$   

We write $\langle x,y \rangle$ for the 
(Wiener-Kuratowski) ordered pair $\{\{x\},\{x,y\}\}$. 
The ordered pair and the corresponding projection functions
are defined by stratified formulas.  To wit, the formula that 
expresses $z = \langle x,y \rangle$ is 
$$ u \in z \iff \forall w \in u\,( w = x) \lor \forall w \in u\, (w = x \lor w = y),$$
which is stratifiable.  Note that the ordered pair gets an index 2 more than
the indices of the paired elements.%
\footnote{The axiom of infinity is needed to construct an
ordered pair that does not raise the type level.  See 
~\cite{rosser1978}, p.~280.} 
Then we have the basic property
\begin{lemma} \label{lemma:ordered_pair_equality}
$\langle x,y \rangle = \langle a,b \rangle \iff x = a \ \land \ y = b$.
\end{lemma}

\noindent\begin{proof} Straightforward application of the definition
and extensionality.  We omit the approximately 70-step proof. 
\end{proof}
\smallskip

As usual, a function is a univalent set of ordered pairs.  We note
that being a function in \NF\ is a strong condition.  For example, 
$\{x\}$ exists for every $x$, but the map $x \mapsto \{x\}$ is not 
a function in NF, since to stratify an expression involving 
ordered pairs, the elements $x$ and $y$ of $\langle x, y\rangle$
must be given the same index, while in the example, $\{x\}$ must 
get one higher index than $x$. 

Because the ordered pair raises types by two levels,  we define
ordered triples by 

\begin{definition}[Ordered triples]\label{definition:triples} 
$$\langle x,y,z\rangle := \langle \langle x,y \rangle, \{\{z\}\}\rangle.$$
\end{definition}

Then a function of two variables is definable in \INF\ if its graph forms
a set of ordered triples $\langle x,y,f(x,y)\rangle$.

We can conservatively add function symbols for binary union $x \cup y$ and intersection $x \cap y$, 
union, intersection,
set difference $x-y$, 
and generally we can add a function symbol $c_\phi$ for any stratified formula
$\phi$, so that $x \in c_\phi(y) \iff \phi(x,y)$.  For a detailed discussion 
of the logical underpinnings of this step, see~\cite{visser2011}.  Function 
symbols for $\{x\}$, $\{x,y\}$, and $\langle x,y\rangle$ are also special cases
of the $c_\phi$; we can add these function symbols even though the ``functions'' they
denote are not functions in the sense that their graphs are definable in \INF.  Thus 
for example we have 

\begin{lemma} \label{lemma:singleton1}  For all $x,u$:  $u \in \{x\} \iff u = x$
\end{lemma}

\noindent\begin{proof}  This is the defining axiom for the function symbol $\{x\}$,
which is really just $\{u : u = x\}$; that is, the function symbol is $c_\phi$
where $\phi(u,x)$ is $u=x$.\end{proof}

\begin{lemma}\label{lemma:minus_members} For all $a,b,x\,(x \in a - b \iff x \in a \land x \not\in b)$.
\end{lemma}

\noindent\begin{proof}  This can be taken as the defining axiom for $a-b$; or it may be derived in a finite axiomatization
from other axioms. 
\end{proof}

\begin{lemma} \label{lemma:single_oneone} For all $x,y$,  $\{x\} = \{y\} \iff x=y$.
\end{lemma}

\noindent\begin{proof}  Right to left is just equality substitution.  Ad left to right:
Suppose $\{x\} = \{y\}$.  Then
\begin{eqnarray*}
u \in \{x\} \iff u \in \{y\} &&\mbox{\qquad by extensionality} \\
u = x \iff u = y    && \mbox{\qquad by Lemma~\ref{lemma:singleton1}}\\
x = y                && \mbox{\qquad by equality axioms}
\end{eqnarray*}
\end{proof}

\noindent
 \subsection*{  Technical details about stratification}
In practice we need to use stratified comprehension in the presence of 
function symbols and parameters; the notion of stratification has to 
be extended to cover these situations.  We define the notion of a 
formula $\phi$ being ``stratified with respect to $x$''.  The 
variables of $\phi(x)$ are of three kinds:  $x$ (the ``eigenvariable''),
variables other than $x$ that occur only on the right of $\in$ (``parameters''),
and all other variables.  An assignment of natural numbers (indices)
to the variables that are not parameters
 is said to stratify $\phi$ with respect to $x$ if for each atomic
 formula $z \in y$, $y$ is assigned an index one larger than the index
 assigned to $z$.   Note that the assignment is to variables, rather than
 occurrences of variable, so every occurrence of $z$ gets the same index.
 Note also that parameters need not be assigned an index.  
 \smallskip
 
 Now when terms are allowed, built up from constants and function symbols
 that are introduced by definitions, an assignment of indices must be 
 extended from variables to terms.  When we introduce a function symbol,
 we must tell how to do this.  
 For example, the ordered pair $\langle x, y\rangle$ must have $x$ and 
 $y$ assigned the same index, and then the pair gets an index two greater.
 The singleton $\{x\}$ must get an index one more than $x$, and so on.
 Stratified comprehension in the extended language says that 
 $\{ x : \Phi(x)\} $ exists, when $\Phi$ is stratified with respect to 
 $x$.  The set so defined will depend on any free variables of $\Phi$
 besides $x$,  some of which may be parameters and some not.  
 
 It is ``well-known'' that stratified comprehension, so defined, is 
 conservative over NF, but it does not seem to proved in the standard
 references on NF; and besides,
 we need that result for \INF\ as well.  The algorithm in 
~\cite{visser2011} meets the need: it will unwind the function symbols
 in favor of their definitions,  preserving stratification.%
 \footnote{The referee pointed out \cite{adlesic},
 which is recent enough that it cites a preprint version of this paper.
  Theorem~3.11 of that paper is the key step in proving that stratified 
  comprehension in the extended language is conservative over NF.}
    The confused reader is advised to work this out on paper for the 
 example of the binary function symbol $\langle x,y\rangle$.  
 
In our work, we repeatedly assert that certain formulas are stratifiable,
and then we apply comprehension, either directly or indirectly by using
mathematical induction or induction on finite sets.  The question then 
arises of ensuring that only correctly 
stratified instances of comprehension are used.  One approach is to 
use a finite axiomatization of \INF.  (It is easy to write one down
following well-known examples for classical NF.)  But that just pushes
the problem back to verifying the correctness of that axiomatization;
moreover it is technically difficult to reduce given particular instances
of comprehension to a finite axiomatization.  Instead, we just made a 
list of each instance of comprehension that we needed.  There were
at some point 154
instances of comprehension in that list (which includes more than 
just the instances used in this paper).   The Lean proof assistant 
does not check that those instances are stratified.  If one is not 
satisfied with a manual check of those 154 formulas,
then one has to write a computer program to check that they are stratified.
We did write one and those 154 formulas passed;  since this paper is 
being presented as human-readable, we rely here on the human reader 
to check each stratification as it is presented; we shall not go
into the technicalities of computer-checking stratification.

\subsection*{Functions and functional notation}

\begin{definition} \label{definition:maps}
$f:X \to Y$  (``$f$ maps $X$ to $Y$'') means for every $x \in X$ there
exists a unique $y \in Y$ such that $\langle x,y \rangle \in f$.
``$f$ is a function'' means 
$$ \langle x,y\rangle \in f \ \land \ \langle x,z \rangle \in f \imp y = z.$$
\end{definition}
The domain and range of $f$ are defined as usual, so $f$ is a function if and only 
if  it maps its domain to its range.

When $f$ is a function, one writes $f(x)$ for that unique $y$.  It is time 
to justify that practice in the context of \INF.%
\footnote{Rosser's version of classical \NF\ has
Hilbert-style choice operator, which gives us ``some $y$ such that $\langle x, y \rangle \in f$.''
But \INF\ does not and cannot have such an operator, so a different formal treatment is needed.}  
Here is how to do that.
We introduce a function symbol $Ap$  (with the idea that we will abbreviate $Ap(f,x)$
to $f(x)$ informally).

\begin{definition} \label{definition:Ap}
$$ Ap(f,x) = \{ u : \exists y\,(\langle x,y\rangle \in f \ \land \ u \in y \}.$$
\end{definition}
It is legal to introduce $Ap$ because it is a special case of a stratified comprehension
term.  One can actually introduce the symbol $Ap$ formally, or one can regard $Ap$ as an 
informal abbreviation for the comprehension term in the definition.  Informally we 
are going to abbreviate $Ap(f,x)$ by $f(x)$ anyway, so $Ap$ will be invisible in the informal
development.  This procedure is justified by the following lemma:

\begin{lemma} \label{lemma:Ap}
If $f$ is a function and $\langle x,y \rangle \in f$, then $y = Ap(f,x)$.
\end{lemma}

\noindent\begin{proof}  Suppose $f$ is a function and $\langle x,y \rangle \in f$.
We must prove $y = Ap(f,x)$.  By extensionality it suffices to show that for all $t$,
\begin{eqnarray}
t \in y \iff t \in Ap(f,x) \label{eq:343}
\end{eqnarray}
{\em Left to right}.Suppose $t \in y$.  Then by the definition of $Ap$, we have $t\in Ap(f,x)$.
\smallskip

{\em Right to left}.  Suppose $t \in Ap(f,x)$.  Then by the definition of $Ap$, for some $z$ we 
have $\langle x,z \rangle \in f$ and $t \in z$.  Since $f$ is a function, $y = z$.  Then $t \in y$.
That completes the right-to-left direction.  
\end{proof}

\subsection*{One-to-one, onto, and similarities} 

The function 
$f:X \to Y$ is {\bf one-to-one} if 
$$y \in Y \ \land \ \langle x,y \rangle \in f 
\imp x \in X$$
and for $x,z \in X$ we have 
$$ f(x)= f(z) \imp x = z.$$
If $f: X \to Y$ is one-to-one then we define
$$f^{-1} = \{ \langle y,x  \rangle : \langle x,y \rangle \in f \}.$$
The definition of $f^{-1}$  can be given by a stratified formula, 
so it is legal
in \INF. 

\noindent
{\em Remark}.  We could also consider the 
notion of ``weakly one-to-one'': 
$$x,y \in X \land x \neq y \imp f(x) \neq f(y).$$
The two notions are not equivalent unless equality 
on $X$ and $Y$ is {\bf stable},
meaning $\neg \neg\, x = y \imp x =y$.  Since equality on 
finite sets is decidable, the two notions do coincide on finite sets,
but we need the stronger notion in general, in particular, to 
make the notion of ``similarity'' in the next definition be an
equivalence relation.  The point is that the stronger notion is 
needed for the following lemma.

\begin{lemma} \label{lemma:finverse}
The inverse of a one-to-one function from $X$ onto $Y$
 is a one-to-one function from $Y$ onto $X$.  That is,
 if $f:X \to Y$ is one-to-one, then $f^{-1}: Y \to X$
 and $f^{-1}$ is one-to-one and onto.   
\end{lemma}

\noindent\begin{proof}  Let $f:X \to Y$ be one-to-one and onto.
 Since $f$ is one-to-one, for each $y \in Y$ there is a unique
$x$ such that $\langle x,y \rangle \in f$.  Then by definition 
of function, $f^{-1}: y \to x$.  Since $f: X \to Y$,
for each $x \in X$ there is a unique $y \in Y$ such that 
$\langle x, y\rangle \in f$.

I say $f^{-1}: Y \to X$.  Let $y \in Y$. Since $f$ is one-to-one,
there exists a unique $x \in X$ such that $\langle x,y\rangle \in f$.
That is, $\langle y,x \rangle\in f^{-1}$. Therefore
$f^{-1}: Y \to X$, as claimed. 

I say $f^{-1}$ is one-to-one from $Y$ to $X$. Let $x \in X$;
since $f:X \to Y$ there is $y\in Y$ such that $\langle x,y \rangle \in f$.
Then $\langle y,x \rangle \in f^{-1}$.   Suppose also
$\langle z,x \rangle \in f^{-1}$ with $z \in Y$.  Then 
$\langle x,z \rangle \in f$.  Since $f:X \to Y$, we have
$x = z$.  Therefore $f^{-1}$ is one-to-one, as claimed.
 
I say $f^{-1}$ maps $Y$ onto $X$. Let $x \in X$.
Let $y = f(x)$.  Then $\langle x,y\rangle \in f$.
Then $\langle y,x \rangle \in f^{-1}$.  Therefore
$f^{-1}$ is onto, as claimed. 
 \end{proof}

\begin{definition}\label{definition:similar} The relation ``$x$ is similar to $y$''
is defined by 
$$ x \sim y \iff \exists f\,(f:x \to y\ \land \mbox{\, $f$ is one-to-one and onto}).$$
In that case, $f$ is a {\bf similarity} from $x$ to $y$.
\end{definition}
\noindent
The defining formula is stratified giving $x$ and $y$ the same type,
so the relation is definable in \INF.

\begin{lemma}\label{lemma:sim} The relation $x \sim y$ is an 
equivalence relation.  
\end{lemma}

\noindent\begin{proof} Ad reflexivity:  $x \sim x$ because the 
identity map from $x$ to $x$ is one-to-one and onto. 
\smallskip

Ad symmetry: Let $x \sim y$.  Then there exists a one-to-one 
function $f:x \to y$. By Lemma~\ref{lemma:finverse},
there exists a function $f^{-1}:y \to x$ that is one-to-one 
and onto.  
Hence $y \sim x$.  That completes the proof of symmetry.
\smallskip

Ad transitivity: Let $x \sim y$ and $y \sim z$.  Then there
exist $f$ and $g$ such that $f:x \to y$ is one-to-one and onto,
and $g: y \to z$ is one-to-one and onto.  Then $f \circ g : x \to z$
is one-to-one and onto.  Therefore $x \sim z$.  That
completes the proof of transitivity.  \end{proof}

\begin{lemma} \label{lemma:similar_to_empty2} For all $x$,
$$ x \sim \emptyset \iff x = \emptyset.$$
\end{lemma}

\noindent\begin{proof} {\em Left to right}. suppose $x \sim \emptyset$.  Let $f:x \to \emptyset$
be a similarity.  Suppose $u \in x$.  Then for some $v$, $\langle u,v \rangle \in f$
and $v \in \emptyset$.  But $v \not\in \emptyset$.  Hence $u \not\in x$.  Since $u$
was arbitrary, $x = \emptyset$, as desired.
\smallskip

{\em Right to left}.  Suppose $x = \emptyset$. We have to show $\emptyset \sim \emptyset$.  But 
$\emptyset: \emptyset \to \emptyset$ is a similarity.  
\end{proof} 

\begin{lemma} \label{lemma:subsets_to_equal} $a \subseteq b \ \land \ b \subseteq a \iff a = b$.
\end{lemma}

\noindent\begin{proof} By the definition of $\subseteq$ and the axiom of extensionality.
\end{proof}

 \section{Finite sets}

\begin{definition} \label{definition:finite}
The set \FINITE\ of finite sets is defined as the  
intersection of all $X$ such that $X$ contains the empty set $\emptyset$
and
$$ u \in X \land z \not \in u \imp u \cup \{ z \} \in X.$$
\end{definition}
The  formula in the definition 
 can be stratified by giving $u$ index 1, $z$ index 0, 
and $X$ index 2, so the definition can be given in \INF.  
\smallskip

This definition was introduced in~\cite{Dzierzgowski1995} as ``$N$-finite.''%
\footnote{
He also defined other notions of ``finite''; for example $K$-finite drops
the requirement $z \not \in u$ from the definition.  That notion,
and the other notion considered {\em op.\,cit.}, do not satisfy the property that 
the cardinality of a finite set is a finite cardinal, i.e., an integer.
For example, $\{c\}$ will be $K$-finite, even if we do not know whether 
or not $c$ is inhabited, so we cannot assign $\{c\}$ a finite cardinal.
}

\begin{definition} \label{definition:decidable}
The set $X$   {\bf has decidable equality} if 
$$\forall x,y \in X(x = y \lor x \neq y).$$
The class  $\DECIDABLE$ is the class of all sets 
having decidable equality.
\end{definition}
The formula defining decidable equality is stratified,
so the class $\DECIDABLE$ can be proved to exist.

\begin{lemma}\label{lemma:finitedecidable} Every finite set
has decidable equality.  That is,  $\FINITE \subseteq \DECIDABLE$.
\end{lemma}

\noindent\begin{proof}  Let $Z$ be the set of finite sets with 
decidable equality.  I say that $Z$ satisfies the closure conditions
in the definition of \FINITE, Definition~\ref{definition:finite}.
The empty set has decidable equality, so the first condition
holds.  Now suppose $Y = X \cup \{a\}$, where $X \in Z$ and $a \not \in X$.
 We must 
show $Y \in Z$.  Let $x,y \in Y$.  Then $x \in X \lor x = a$ and 
$y \in X \lor y = a$.  There are thus four cases to consider:
If both $x$ and $y$ are in $X$, then by the induction hypothesis, we
have the desired $x = y \lor x \neq y$. If one of $x,y$ is in $X$
and the other is $a$, then $x \neq y$, since $a \not \in X$; hence
$x = y \lor x \neq y$.  Finally if both are equal to $a$, then $x=y$
and hence $x=y\lor x \neq y$.  Therefore, as claimed, $Z$ satisfies
the closure conditions.  Hence every finite set belongs to $Z$.
\end{proof}

\begin{lemma}\label{lemma:empty_or_inhabited}
A  finite set is empty or it is inhabited (has a member).
\end{lemma}

\noindent\begin{proof} Define
$$ Z = \{ X \in \FINITE\ : X = \emptyset \ \lor \ \exists u\,(u \in X)\}.$$
We will show $Z$ satisfies the closure conditions in the definition
of \FINITE.  Evidently $\emptyset \in Z$.  Now suppose 
$X \in Z$ and $Y = X \cup \{a\}$ with $a \not\in X$.
  We must show $Y \in Z$.
Since $X \in Z$, $X$ is finite.  Therefore $Y$ is finite.
Since $a \in Y$ we have $Y \in Z$.
 \end{proof}

\begin{corollary} [Finite Markov's principle] \label{lemma:markov}
For every finite set $X$
$$ \neg\neg\,\exists u\,(u \in X) \imp \exists u\,(u \in X).$$
\end{corollary}

\noindent\begin{proof}  Let $X$ be a finite set. 
Suppose $ \neg\neg\,\exists u\,(u \in X)$.  
That is, $X$ is nonempty.  By Lemma~\ref{lemma:empty_or_inhabited},
$X$ has a member.  \end{proof}
 
\begin{lemma} \label{lemma:lambda_finite}  $\emptyset \in \FINITE$.
\end{lemma}

\noindent\begin{proof} $\emptyset$ belongs to every 
set $W$ containing $\emptyset$ and containing $u \cup \{e\}$ whenever
$u \in W$ and $e \not\in W$.  Since $\FINITE$ is the intersection of 
such sets $W$,  $\emptyset \in \FINITE$.  \end{proof} 
 
\begin{lemma} \label{lemma:finite_adjoin}
If $x \in \FINITE$ and $c \not\in x$, then $x \cup \{c\} \in \FINITE$.
\end{lemma}

\noindent\begin{proof} Let $x \in \FINITE$.  Then $x$ belongs to every 
set $W$ containing $\emptyset$ and containing $u \cup \{e\}$ whenever
$u \in W$ and $e \not\in W$.  Let $W$ be any such set. Then $x \cup \{c\} \in W$.
Since $W$ was arbitrary, $x \cup \{c\} \in \FINITE$.  
\end{proof}

\begin{lemma} \label{lemma:finite_structure}
If $z \in \FINITE$, then $z = \emptyset$ or 
there exist $x \in \FINITE$ and $c \not \in x$
such that $z = x \cup \{c\}$.
\end{lemma}

\noindent\begin{proof} The formula is stratified,
giving $c$ index 0, and $x$ and $z$ index 1.  $\FINITE$ is 
a parameter.   We prove the formula by induction on finite sets.  Both 
the base case (when $z = \emptyset$) and the induction step are 
immediate.   \end{proof} 

\begin{lemma} \label{lemma:singletons_finite}
Every unit class $\{ x\}$ is finite.
\end{lemma}

\noindent\begin{proof}  We have
\begin{eqnarray*}
\emptyset \in \FINITE && \mbox{\qquad by Lemma~\ref{lemma:lambda_finite}}\\
x \not \in \emptyset  && \mbox{\qquad by the definition of $\emptyset$}\\
\emptyset \cup \{x\} \in \FINITE && \mbox{\qquad by Lemma~\ref{lemma:finite_adjoin}}\\
\{x\} = \emptyset \cup \{x\}  && \mbox{\qquad by the definitions of $\cup$ and $\emptyset$}\\
\{x\} \in \FINITE && \mbox{\qquad by the preceding two lines}
\end{eqnarray*}
\end{proof}

\begin{definition} $\USC(x) := \{ \{y\} : y \in x\}$.  
\end{definition}                                       

\begin{lemma} \label{lemma:uscfinite}
 $\USC(x)$ is finite if and only if $x$ is finite.
\end{lemma}

\noindent\begin{proof}  Left-to-right: we have to prove
\begin{eqnarray}
\forall y \in \FINITE\, \forall x\,( y = \USC(x)  \imp x \in \FINITE) \label{eq:449}
\end{eqnarray}
The formula is weakly stratified with respect to $y$, as 
we are allowed to give the two occurrences of $\FINITE$ different
types.  So we may prove the formula by induction on finite sets $y$.
\smallskip

{\em Base case}.  When $y = \emptyset = \USC(x)$ we have $x = \emptyset$,
so $x \in \FINITE$. 
\smallskip

{\em Induction step}.  Suppose $y \in \FINITE$ has 
the form $y = z \cup \{w\} = \USC(x)$ and $w \not\in z$,
and $z \in \FINITE$. 
Then $w = \{c\}$ for some $c\in x$.
 The induction hypothesis is 
 \begin{eqnarray}
 \forall w\,( z = \USC(w)  \imp w \in \FINITE) \label{eq:320}
\end{eqnarray}
  Then 
\begin{eqnarray*}
z &=& y - \{w\} \mbox{\qquad\qquad\qquad since $y = z \cup \{w\}$} \\
&=& \USC(x) - \{\{c\}\} \mbox{\qquad\ since $y = \USC(x)$ and $w=\{c\}$} \\
&=& \USC(x-\{c\}).
\end{eqnarray*}
Since $y \in \FINITE$ and $\{c\} \in y$, we have 
\begin{eqnarray*}
q \in y \imp q = \{c\} \ \lor \ q \neq \{c\} && 
            \mbox{\qquad by Lemma~\ref{lemma:finitedecidable}}\\
u \in x \imp \{u \} = \{c\} \ \lor \ \{u\} \neq \{c\} && 
            \mbox{\qquad since $y = \USC(x)$} \\
u \in x \imp  u = c \ \lor \ u \neq c 
\end{eqnarray*}
It follows that 
\begin{eqnarray}
(x-\{c\}) \cup \{c\} &=& x \label{eq:338}
\end{eqnarray} 
 By 
the induction hypothesis (\ref{eq:320}), with $x-\{c\}$ substituted
for $w$,  we have
\begin{eqnarray*}
x-\{c\} \in \FINITE  &&\\ 
(x-\{c\}) \cup \{c\} \in \FINITE && \mbox{\qquad by definition of $\FINITE$}\\
x \in \FINITE &&\mbox{\qquad by (\ref{eq:338})}
\end{eqnarray*}
That completes the induction step.
That completes the proof of the left-to-right implication.
\smallskip

Right-to-left: We have to prove
\begin{eqnarray} 
x \in \FINITE \imp \USC(x) \in \FINITE \label{eq:474}
\end{eqnarray}
Again the formula is weakly stratified since $\FINITE$ is a parameter.
We proceed by induction on finite sets $x$.
\smallskip

{\em Base case}\,: $\USC(\emptyset) = \emptyset \in \FINITE$.
\smallskip

{\em Induction step}\,: 
We have for any $x$ and $c \not\in x$,  
$$\USC(x \cup \{c\}) = \USC(x) \cup \{\{c\}\}.$$
Let $c \not\in x$ and $x \in \FINITE$.
By the induction hypothesis (\ref{eq:474}), $\USC(x)$ is finite,
and since $c \not\in x$, we have $\{c\} \not\in \USC(x)$.
Then $\USC(x) \cup \{\{c\}\}$ is finite. Then 
$\USC(x \cup \{c\})$ is finite.  That completes the induction step.
\end{proof}

\begin{lemma}\label{lemma:union}
The union of two disjoint finite sets is finite.
\end{lemma}

\noindent\begin{proof}  We prove by induction on finite sets $X$ that
$$ \forall\, Y \in \FINITE\,(X \cap Y = \emptyset \imp X \cup Y \in \FINITE).$$
{\em Base case}\,:  $\emptyset \cup Y = Y$ is finite.
\smallskip

{\em Induction step}\,:  Suppose $X = Z \cup \{b\}$ with $b \not\in Z$ 
and $Y \cap (Z \cup \{b\}) = \emptyset$ 
and $Z$ finite. 
Then 
\begin{eqnarray}
X \cup Y &=& (Z \cup Y) \cup \{b\}\nonumber \\
X \cup Y &=& Z \cup (Y \cup \{b\}) \label{eq:677} 
\end{eqnarray}
Since $Y \cap (Z \cup \{b\}) = \emptyset$, $b \not \in Y$. 
Then by the definition of $\FINITE$, $Y \cup \{b\}$ is finite.
We have
$$ Z \cap (Y \cup \{ b\}) = Y \cap ( Z \cup \{b\}) = \emptyset.$$
Then by the induction hypothesis, $Z \cup (Y \cup \{b\})$
is finite.  Then by (\ref{eq:677}), $X \cup Y$ is finite.
That completes the induction step.
\end{proof} 

\begin{lemma} \label{lemma:similar_decidable}
If $x$ has decidable equality, and $x \sim y$,
then $y$ has decidable equality.
\end{lemma}

\noindent\begin{proof}  Suppose $x \sim y$. Then
there exists $f:x \to y$ with $f$ one-to-one and onto.
 By Lemma~\ref{lemma:finverse},
$f^{-1}: y \to x$ is a one-to-one function.  Then we have for 
$u, v \in y$, 
\begin{eqnarray}
 u = v \iff f^{1}(u) = f^{-1}(v). \label{eq:427}
 \end{eqnarray}
Since $x$ has decidable equality, we have
$$  f^{1}(u) = f^{-1}(v) \ \lor \  f^{1}(u)\neq f^{-1}(v).$$
By (\ref{eq:427}), 
$$ u = v \ \lor \ u \neq v.$$
Therefore $y$ has decidable equality. \end{proof}
\smallskip

\begin{lemma} \label{lemma:similarityrestricted}
Let $f:z \cup \{c\} \to y$ be one-to-one and onto.
Suppose $c \not\in z$, and let $g$ be $f$ 
restricted to $z$.  Then $g:z \to y - \{f(c)\}$ is 
one-to-one and onto.
\end{lemma}

\noindent{\em Remark.} Somewhat surprisingly,
it is not necessary to assume that $z \cup \{c\}$
has decidable equality.  That is not important as
decidable equality is available when we use this lemma.
\medskip

\noindent\begin{proof} Let $q = f(c)$. Then
$g:z \to y - \{q\}$. Suppose $g(u) = g(v)$.
Then $f(u) = f(v)$. Since $f$ is one-to-one, 
$u = v$.  Hence $g$ is one-to-one.  Suppose 
$v \in  y-\{q\}$.  Since $f$ is onto, $v = f(u)$
for some $u \in z \cup \{c\}$; but $u \neq c$ since
if $u=c$ then $v = f(u)= q$, but $v \neq q$ since
$v \in y-\{q\}$. Then $u \in z$. Hence $g$ is onto.
\end{proof}

\begin{lemma}\label{lemma:finitesimilar}
A set that is similar to a finite set is finite.
\end{lemma}

\noindent\begin{proof}  We prove by induction on finite sets
$x$ that 
$$ \forall y\, (y \sim x \imp y \in \FINITE).$$
The formula is stratified, so induction is legal.
\smallskip

{\em Base case}\,: When $x = \emptyset$.  Suppose $y \sim \emptyset$.
Then $y = \emptyset$, so $y \in \FINITE$.  That completes the base case.
\smallskip

{\em Induction step}\,: Suppose the finite set $x$ has the form 
$x = z \cup \{c\}$ with $c \not\in z$,
and $x \sim y$.  By Lemma~\ref{lemma:finitedecidable},
$x$ has decidable equality.
Then by Lemma~\ref{lemma:similar_decidable},
$y$ has decidable equality. 
Let $f:z \cup \{c\} \to y$ be $f$ one-to-one and onto.
Let $q = f(c)$.  Then $\langle c,z\rangle \in f$. 
Let $g$ be $f$ restricted to $z$.  By 
Lemma~\ref{lemma:similarityrestricted},
 $g:z \to y - \{q\}$
is one-to-one and onto.   Then by the induction hypothesis, 
$y - \{q\}$ is finite.  Then $(y-\{q\}) \cup \{q\} \in \FINITE$,
by the definition of $\FINITE$.  But since $y$ has 
decidable equality, we have 
$$ y =  (y-\{q\}) \cup \{q)\}.$$
Therefore $y \in \FINITE$.  That completes the induction step.
\end{proof}

\begin{definition}\label{definition:sc} The {\bf  power set}
of a set $X$ is defined as the set of   subclasses of $X$:
$$ \SC(X) = \{ Y: Y \subset X   \}.$$
\end{definition}

We shall not make use of $\SC(X)$, because there are ``too many''
subclasses of $X$.  Consider, by contrast, the
 separable subclasses of $X$:

\begin{definition} \label{definition:SSC}
We define the set of {\bf separable subclasses} of $X$ by 
$$ \SSC(X):= \{ u : u \subseteq X \ \land \ X = u \cup (X-u)\}
$$
\end{definition}
That is, $u$ is a separable subclass (or subset, which is synonymous)
 of $X$ if and only if $\forall y \in X\,(y \in u \ \lor \ y \not\in u)$.  Classically, of course,
every subset is separable, so we have $\SSC(X) = \SC(X)$,
but that is not something we can assert constructively. 
The formula in the definition is stratified, so the definition
can be given in \INF.   When working with finite sets, $\SSC(X)$
is a good constructive substitute for $\SC(X)$.  We illustrate
this by proving some facts about $\SC(X)$, before returning to 
the question of the proper constructive substitute for $\SC(X)$ 
when $X$ is not necessarily finite.

\begin{lemma}\label{lemma:finitepowerset}
Let $x$ be a finite set. Then $\SSC(x)$ is also a finite set.
\end{lemma}

\noindent{\em Remark}.  We cannot prove this with $\SC(x)$ 
in place of $\SSC(x)$.
\medskip 

\noindent\begin{proof} The formula to be proved is 
$$ x \in \FINITE \imp \SSC(x) \in \FINITE.$$
The formula is weakly stratified because the two occurrences of 
the parameter $\FINITE$ may receive different indices.
Therefore we can proceed by induction on finite sets $x$.
\smallskip

{\em Base case}\,:  $\SSC(\emptyset) = \{\emptyset\}$ is finite.
\smallskip

{\em Induction step}\,: Suppose $x$ is finite and consider $x \cup \{c\}$
with $c \not\in x$.  Then $x \cup \{c\}$ is finite and hence,
by Lemma~\ref{lemma:finitedecidable}, it has decidable equality.

By the induction hypothesis, $\SSC(x) \in \FINITE$. 
I say that the map $u \mapsto u \cup \{c\}$ is definable in \INF:
$$ f := \{ \langle u,y \rangle :  u \in \SSC(x) \ \land \ 
                       y = u \cup \{c\}.\}$$
The formula can be stratified by giving $c$ index 0,
$u$ and $y$ index 1, $\SSC(x)$ index 2; then $\langle u,y\rangle$
has index 3 and we can give $f$ index 4.  Hence
$f$ is definable in \INF\ as claimed. $f$ is a function
since $y$ is uniquely determined as $u \cup \{c\}$ when 
$u$ is given. Also $f$ is one-to-one, 
 since if  $u\subseteq x$ and 
 $v \subseteq x$ and $c \not \in x$, and
$u \cup \{c\} = v \cup \{c\}$, then $u=v$. 
Define 
$$ A:= Range(f). $$
Then
\begin{eqnarray}
A = \{ u \cup \{c\}: u \in \SSC(x)\}. \label{eq:593}
\end{eqnarray}
Then $\SSC(x) \sim A$, because $f:\SSC(x) \to A$ is 
one-to-one and onto. 
Since  $\SSC(x)$ is finite (by the induction 
hypothesis),  by Lemma~\ref{lemma:finitedecidable},
$\SSC(x)$ has decidable equality.  Then 
$A$ has decidable equality, by Lemma~\ref{lemma:similar_decidable}.
Since $A$ has decidable equality, and is similar
to the finite set $\SSC(x)$,  $A$ is finite,
 by Lemma~\ref{lemma:finitesimilar}. 
 
I say that 
\begin{eqnarray}
\SSC(x \cup\{c\}) &=& A \cup \SSC(x). \label{eq:461}
\end{eqnarray}
By extensionality, it suffices to show that the two 
sides of (\ref{eq:461}) have the same members.
\smallskip

Left-to-right:  Let $v \in \SSC(x \cup \{c\})$.
Then $v$ is a separable subset of $x \cup \{c\}$.
Then $c \in v \ \lor \ c \not\in v$.  If $c \not \in v$ then 
$v \in \SSC(x)$.  If $c \in v$ 
\begin{eqnarray*}
x \cup \{c\} \in \FINITE & \mbox{\qquad since $x \in \FINITE$ and $c \not\in x$} \\
x \cup \{c\} \mbox{\ \ has decidable equality} & \mbox{\qquad by Lemma~\ref{lemma:finitedecidable}}\\
v \mbox{ \ \ has decidable equality} & \mbox{\qquad since $v \subseteq x \cup \{c\}$} \\
v  = (v-\{c\}) \cup \{c\} & \mbox{\qquad since $x \in v \imp x = c \ \lor \ x \neq c$} 
\end{eqnarray*}
We have $v -\{c\} \in \SSC(x)$, since $v \subset x \cup \{c\}$ and $v$ 
has decidable equality. Then $f(v -\{c\}) \in Range(f) = A$.
But $f(v-\{c\}) = (v-\{c\}) \cup \{c\} = v$.  Therefore $v \in A$.
Therefore $v \in A \cup \SSC(x)$, 
as desired.
That completes the proof of the left-to-right direction 
of (\ref{eq:461}).
\smallskip

Right-to-left. 
Let $v \in A \cup \SSC(x)$.  Then $v \in A \ \lor \ v \in \SSC(x)$.
\medskip

Case~1, $v \in A$. Then by (\ref{eq:593}), $v$
has the form $v= u \cup \{c\}$ for some $u \in \SSC(x)$.
Then $u \cup \{c\} \in \SSC(x \cup \{c\})$ as required.
\smallskip

Case~2, $v \in \SSC(x)$. 
First we note that if $c \not \in x$ then
\begin{eqnarray*}
\SSC(x) \subseteq \SSC(x \cup \{c\})  
\end{eqnarray*}
Therefore, since $v \in \SSC(x)$, we have 
$v \in \SSC(x \cup \{c\})$.  
 That completes
the proof of (\ref{eq:461}).  
\smallskip

Note that $A$ and $\SSC(x)$ are disjoint, since every member of $A$
contains $c$, and no member of $\SSC(x)$ contains $c$, since $c \not\in x$.
Then by Lemma~\ref{lemma:union}
and (\ref{eq:461}), $\SSC(x \cup \{c\}) \in \FINITE$,
as desired.  \end{proof}

\begin{lemma} \label{lemma:finiteseparable}
A finite subset of 
a finite set is a separable subset. 
\end{lemma}

\noindent\begin{proof}  Let $a \in \FINITE$.
 By induction on finite sets $b$
we prove
\begin{eqnarray}
b \in \FINITE \imp b \subseteq a \imp a = (a-b) \cup b. \label{eq:711}
\end{eqnarray}
The formula is stratified, so induction is legal.
\smallskip

{\em Base case}\,: Suppose   $b = \emptyset$.  Then 
$b \subseteq a$, so we have to prove $a = (a - \emptyset) \cup \emptyset$,
which is immediate. That completes
the base case.
\smallskip

{\em Induction step}\,: Suppose $b \in \FINITE$ and $c \not \in b$ 
and $b \cup \{c\} \subseteq a $.  We must show
\begin{eqnarray*}
a &=& (a - (b \cup \{c\})) \cup (b \cup \{c\})  
\end{eqnarray*}
By extensionality, it suffices to show that 
\begin{eqnarray}
x \in a &\iff& x \in (a - (b \cup \{c\}) \cup (b \cup \{c\}) \label{eq:723}
\end{eqnarray}
Since $a$ is finite, $a$ has decidable equality, by 
Lemma~\ref{lemma:finitedecidable}. 

Ad left-to-right of (\ref{eq:723}):  Let $x \in a$.
Then by decidable equality on $a$, we have
\begin{eqnarray}
x = c \ \lor \ x \neq c   \label{eq:729}
\end{eqnarray}
By the induction hypothesis (\ref{eq:711}), we have 
\begin{eqnarray}
x \in b \ \lor \ x \not \in b  \label{eq:736}
\end{eqnarray}
By (\ref{eq:729}) and (\ref{eq:736}) we have 
\begin{eqnarray}
x \in b \cup \{c\} \ \lor \ x \not \in b \cup \{c\} \label{eq:739}
\end{eqnarray}
Therefore 
$$ x \in (a - (b \cup \{c\}) \cup (b \cup \{c\}).$$
That completes the left-to-right implication in (\ref{eq:723}).
\smallskip

Ad right-to-left: Suppose 
$$ x \in (a - (b \cup \{c\}) \cup (b \cup \{c\}).$$
We must show $x \in a$.  If $x \in (a - (b \cup \{c\})$
then $x \in a$.  If $x \in   (b \cup \{c\})$ then 
$x \in a$, since by hypothesis $b \cup \{c\} \subseteq a$.
That completes the right-to-left direction.
That completes the induction step. 
\end{proof}

\begin{lemma} \label{lemma:separablefinite}
 Every separable subset of a finite set is finite.
\end{lemma} 

\noindent\begin{proof}
By  induction on finite sets $X$.
When $X$ is the empty set,  every subset of $X$ is the empty 
set, so every subset of $X$ is empty, and hence finite.
Now let $X = Y \cup \{a\}$ with $a \not \in Y$ and $Y$ finite, 
and let $U$ be a separable subset of $X$; that is, 
\begin{eqnarray}
\forall z \in X\,(z\in U \ \lor \ z \not\in U).\label{eq:1479}
\end{eqnarray}
We have to show $U$ is finite.
Since $U$ is 
separable, $a \in U \ \lor a \not \in U$; we argue by cases 
accordingly.
\smallskip

Case~1: $a \not\in U$. Then $U \subseteq Y$, so by the induction 
hypothesis, $U$ is finite.
\smallskip

Case~2: $a \in U$.    
Let $V = U - \{a\}$.  Then $V \subseteq Y$.
I say that $V$ is a separable subset of $Y$; that is,
\begin{eqnarray}
\forall z \in Y\,(z\in V \ \lor \ z \not \in V) \label{eq:1484}
\end{eqnarray}
Let $z \in Y$.  Since $U$ is a separable subset of $X$,
$z \in U \lor z \not\in U$.  By Lemma~\ref{lemma:finitedecidable},
$X$ has decidable equality, so 
$z = a \lor z \neq a$. Therefore $z \in V \lor z \not \in V$,
as claimed in (\ref{eq:1484}). Then, by the induction hypothesis, $V$ is finite. 
Since $a \not \in V$,  also $V \cup \{a\}$ is finite.  I say 
that $V \cup\{a\} = U$.  If $x \in V \cup \{a\}$ then $x \in U$,
since $V\subseteq U$ and $a \in U$.  Conversely if $x \in U$ then
$x = a \ \lor x \neq a$, since $a$ and $x$ both are members of $X$
and $X$ has decidable equality by Lemma~\ref{lemma:finitedecidable}.
If $x= a$ then $x \in \{a\}$ and if $x \neq a$ then $x \in V$,
so in either case $x \in V \cup \{a\}$. Therefore 
$V \cup\{a\} = U$ as claimed.  Since $V$ is finite and $a \not \in V$,
$V \cup \{a\}$ is finite.  Since $U = V \cup \{a\}$, 
$U$ is finite.
 That completes the induction step.  
\end{proof}

\begin{lemma} \label{lemma:finitedif}
Let $a$ and $b$ be finite sets with $b \subseteq a$.
  Then $a-b$ is also a finite set.
\end{lemma}

\noindent\begin{proof} 
We first prove the special case when $b$ is a singleton, $b = \{c\}$.
That is, 
\begin{eqnarray}
a \in \FINITE \ \land \ c \in a \imp
 a - \{c\} \in \FINITE \label{eq:776}
\end{eqnarray}
By Lemma~\ref{lemma:finitedecidable}, $a$ has decidable equality.
Hence $a-\{c\}$ is a separable subset of $a$.  Then by 
Lemma~\ref{lemma:separablefinite}, it is finite.  That completes the 
proof of (\ref{eq:776}).
\smallskip

We now turn to the proof of the theorem proper.
 By induction on finite sets $a$ we prove
$$  \forall b \in \FINITE\, (b \subseteq a \imp (a-b) \in \FINITE).$$
{\em Base case}\,: $\emptyset - b = \emptyset$ is finite.
\smallskip

{\em Induction step}.  Let $a = p \cup \{c\}$, with $c \not\in p$.
Let $b$ be a finite subset of $a$.  
 We have $c \in b \ \lor c \not\in b$
by Lemma~\ref{lemma:finiteseparable}. We argue 
by cases accordingly.
\smallskip

Case~1: $c \in b$. Then 
\begin{eqnarray*}
a-b &=& p \cup \{c\} - b \\
&=&  p-b \\
&=&  p- (b -\{c\}) \mbox{\qquad since $c \not\in p$ and $c \in b$}
\end{eqnarray*}
Since $b$ is finite, also $b-\{c\}$ is finite,
by (\ref{eq:776}).
Since $b - \{c\} \subseteq p$, by the induction 
hypothesis we have 
$$ p - (b - \{c\} ) \in \FINITE.$$
Therefore $p-b \in \FINITE$.  Therefore $a-b \in \FINITE$. 
 That completes Case~1.
\smallskip   

Case~2: $c \not\in b$. Then  $b \subseteq p$, so by
the induction hypothesis $p-b$ is finite.

\begin{eqnarray*}
a-b &=&( p \cup \{c\}) -b \\
   &=& (p-b) \cup \{c\}  \mbox{\qquad since $c \not \in b$}
\end{eqnarray*}
 Therefore
$a-b$ is finite.  That completes Case~2. 
That completes the induction step.  
\end{proof}

\begin{lemma} [Bounded quantification] \label{lemma:boundedquantification2}  
Let $X$ be any
set with decidable equality, and $B$ a finite subset of $X$.
Let $Y$ be any set, with $R$  a separable subset of  $X \times Y$.
  Let $P$ be 
defined by 
$$ z \in P \iff z \in X \ \land \ \exists u \in B\, \langle u,z \rangle \in R$$ 
Then  $P$ is a separable subset of $X$.
With complete precision: 
\begin{eqnarray*}
&& \forall u,v \in X\,(u=v \ \lor \ u \neq v) \ \land \ \\
&& B \in \FINITE   \ \land \ B \subseteq X 
\ \land \ \forall u \in X\, \forall z \in Y\,(\langle u,z \rangle \in R \ \lor \ \neg\, \langle u,z \rangle \in R) \\
&& \imp 
\forall z \in X\,  ( \exists u \in B\, \langle u,z \rangle \in R \ \lor \ \neg\, \exists u \in B\, \langle u,z \rangle \in R)
\end{eqnarray*}
\end{lemma}

\noindent{\em Remark}.  We may express the lemma informally as 
``The decidable sets are closed under bounded quantification.''
\smallskip

\noindent\begin{proof} The formula to be proved is stratified, with 
 $\FINITE$ as a parameter, giving $u$ and $z$ index 0, $B$ index 1,
 and $R$ index 3. 
 Therefore it is legal to prove it by induction on finite sets $B$. 
\smallskip

{\em Base case}\,: $B = \emptyset$.  Then $X \times Y = \emptyset$, so $R = \emptyset$.  Then 
$\forall z\, \neg\,\exists u \in B\, \langle u,z \rangle \in R$, and therefore
$$\forall z\in X\, (\exists u \in B\,\langle u,z \rangle \in R \ \lor \ 
                  \neg\, \exists u \in B\, \langle u,z \rangle \in R).$$ 
  That completes the base case.
\smallskip

{\em Induction step}.  Suppose $B = A \cup \{c\}$ with $A$ finite and 
$c \not\in A$.  Then 
\begin{eqnarray}
 \forall z \in X\,( \exists u \in B\,\langle u,z \rangle \in R  
 &\iff& (\exists u \in A\, \langle u,z \rangle \in R) \ \lor \ \langle c,z \rangle \in R) \label{eq:865}
\end{eqnarray}
We have to prove 
\begin{eqnarray}
&&(\exists u \in B\, \langle u,z \rangle \in R) \ \lor \ \neg\,(\exists u \in B\, \langle u,z \rangle \in R) \label{eq:869} 
\end{eqnarray}
By (\ref{eq:865}), that is equivalent to 
\begin{eqnarray*}
&& (\exists u \in A\, \langle u,z \rangle \in R \ \lor \ \langle c,z \rangle \in R)
\ \lor \ \neg\,(\exists u \in A\, \langle u,z \rangle \in R \ \lor \ \langle c,z \rangle \in R) \\
&\iff& (\exists u \in A\, \langle u,z \rangle \in R \ \lor \ \langle c,z \rangle \in R)
\ \lor \ (\neg\,(\exists u \in A\, \langle u,z \rangle \in R) \ \land \ \langle c,z \rangle \not\in R)) \\
&\iff& (\exists u \in A\, \langle u,z \rangle \in R \ \lor  \neg\,\exists u \in A\, \langle u,z \rangle \in R \ \lor \langle c,z \rangle \in R) \\
&&\ \land \ (\exists u \in A\, \langle u,z \rangle \in R \ \lor  \neg\,\exists u \in A\, \langle u,z \rangle \in R \ \lor \neg\,\langle c,z \rangle \in R)
\end{eqnarray*}
Since $\langle c,z \rangle \in R \ \lor \ \neg\,\langle c,z \rangle \in R$, the last formula is 
equivalent to 
$$ \exists u \in A\, \langle u,z \rangle \in R \ \lor  \neg\,\exists u \in A\, \langle u,z \rangle \in R.$$
But by the induction hypothesis, that holds.  
That completes the induction step.  
\end{proof}

\begin{lemma}[swap similarity] \label{lemma:swap_similarity}
Let $X$ have decidable equality and let $U \subseteq X$ and $b,c \in X$
with $b \in U$ and $c \not\in U$.  Let $Y = U - \{b\} \cup \{c\}$.
Then $U \sim Y$.
\end{lemma}

\noindent\begin{proof} 
 Since $b \in U$ and $c \not\in U$,
we have $b \neq c$. 
Define $f: U \to Y$ by 
$$ f(x) = \twocases { c  \mbox{\qquad if $x = b$} }
                    { x  \mbox{\qquad otherwise}}
$$
Since $X$ has decidable equality, $f$ is well-defined on $X$,
and from the definitions of $f$ and $Y$ we see that $f:U \to Y$
and $f$ is onto.  Ad one-to-one:  suppose $f(u)= f(v)$.  Since
$X$ has decidable equality, $u$ and $v$ are either equal or not.
If $u=v$, we are done. If $u \neq v$ then exactly one of $u,v$
is equal to $b$, say $u = b$ and $v \neq b$.  Then $f(u)= c$
and $f(v) = v$.  Since $f(u) = f(v)$ we have $v =c$.
But $v \in U$ and $c \not\in U$, contradiction.
\end{proof}

\begin{definition} [Dedekind]\label{definition:infinite} \ \\
The class $X$ is {\bf  infinite}  if  $X\sim Y$ 
for some   $Y \subseteq X$ with $Y \neq X$.%
\footnote{ Alternate definitions one might consider:
$X$ is infinite if there is a similarity from $X$ to a 
subset of $X$ that omits a value;  $X$ is infinite if $X-A$
is inhabited, for every finite set $A$.  Whether Dedekind
infinite implies these properties is not known. }
\end{definition} 
 
 \begin{theorem} \label{theorem:infiniteimpliesnotfinite}
Let $X$ be infinite, in the sense that it is similar to 
some $Y \subset X$ with $Y \neq X$.  Then $X$ is not finite.
\end{theorem}

\noindent{\em Remark}.  We expressed the theorem as
``infinite implies not finite'',  but of course it is 
logically equivalent to ``finite implies not infinite'', 
since both forms amount to ``not both finite and infinite.''
\medskip

\noindent\begin{proof}  It suffices to show that every finite set
is not infinite.  The formula to be proved is
$$ X \in \FINITE  \imp  \forall Y\, (Y \subseteq X \imp X \sim Y \imp X = Y)\}.$$
That formula  is stratified, giving $X$ and $Y$ index 1, 
since the similarity relation can be defined in \INF.  Therefore
induction is legal.
\smallskip

{\em Base case}, $X = \emptyset$. The only subset of $\emptyset$ is $\emptyset$, 
so any subset of $X$ is equal to $X$.  That completes the base case.
\smallskip

{\em Induction step}.
Suppose $X = A \cup \{b\}$, with $A \in \FINITE$ and $b \not\in A$.   
 Then
\begin{eqnarray*}
Y \subseteq X && \mbox{\qquad by hypothesis} \\
X \in \FINITE &&\mbox{\qquad by Lemma~\ref{lemma:finite_adjoin}}\\
X \in \DECIDABLE && \mbox{\qquad by Lemma~\ref{lemma:finitedecidable}}\\
X \sim Y \ \land \ Y \subseteq X &&\mbox{\qquad assumption}\\
f:X\to Y && \mbox{\qquad with $f$ one-to-one and onto, by definition of $X \sim Y$}\\ 
Y \in \FINITE  &&\mbox{\qquad by Lemma~\ref{lemma:finitesimilar}}\\
Y \in \DECIDABLE && \mbox{\qquad by Lemma~\ref{lemma:finitedecidable}}
\end{eqnarray*}
Let $c =f (b)$ and $U = Y - \{c\}$. 
Let $g$ be $f$ restricted to $A$. 
Then $g:A \to U$ is one-to-one and onto (140 steps omitted). 
  Thus $A \sim U$.
\smallskip

 Since $X$ has decidable 
equality, $b = c \ \lor \ b \neq c$.
By Lemma~\ref{lemma:finiteseparable}, $Y$ is a separable subset of $X$.
Therefore $b \in Y \ \lor  b \not\in Y$.   We can therefore
argue by three cases:  $b = c$, or $b \neq c$ and $b \in Y$, or 
$b \neq c$ and $b \not\in Y$.  
\smallskip

Case~1, $b=c$.  Then $U \subseteq A$.   By the induction hypothesis, we have 
$A = U$. Then $X = A \cup \{b\}  = U \cup \{b\} = U \cup \{c\} = Y$.
That completes Case~1.
\smallskip

Case~2, $b \neq c$ and $b \in Y$.  Then
\begin{eqnarray*}
f(p) = b  &&\mbox{\qquad for some $p \in A \cup\{b\}$, since $f$ is onto $Y$}\\
p \neq b  && \mbox{\qquad since $f(b) = c \neq b$, and $f$ is one-to-one}
\end{eqnarray*}
Define 
$$ g = (f - \{ \langle b,c \rangle \} - \{\langle p, b\rangle\}) \cup \{\langle p, c\rangle\}.$$
Then one can check that $g: A \to Y - \{b\}$ is one-to-one and onto.  (It requires
more than six hundred inference steps, here omitted.) We note that $A \cup \{b\}$
is finite, and therefore has decidable equality, which allows us to argue by cases
whether $x = b$ or not, and whether $x = p$ or not.) 
Then
\begin{eqnarray*}
A \sim Y - \{b\}  &&\mbox{\qquad since $g$ is a  similarity}\\
Y - \{b\} \subseteq A && \mbox{\qquad since $X = A \cup \{b\}$ and $Y \subseteq $X }
\end{eqnarray*}
Thus $A$ is similar to its  subset $Y - \{b\}$.  Then by the induction 
hypothesis, 
\begin{eqnarray}
A =Y-\{b\} \label{eq:1112}
\end{eqnarray}
Therefore $Y = A \cup \{b\} = X$. 
That completes Case~2. 
\smallskip

Case~3: $b \neq c$ and $b \not\in Y$.
Since $Y \subseteq X = A \cup \{b\}$,  and $b \not\in Y$,
we have $Y \subseteq A$.
Then $f:A \to Y - \{c\} \subseteq A$.  Then by the induction hypothesis, 
\begin{eqnarray}
Y - \{c\} = A.  \label{eq:1128}
\end{eqnarray}
Then $c \not\in A$.  But $X = A \cup \{b\}$, and $c = f(b) \in Y \subseteq A$,
so $c \in A$.  That contradiction completes Case 3, and that completes
the proof of the induction step.  \end{proof}

 \begin{lemma} \label{lemma:finiteunion}
 A finite union of finite disjoint sets is finite.  That is,
 \begin{eqnarray*}
 && x \in \FINITE \ \land \ \forall u\, (u \in x \imp u \in \FINITE) \\
&& \land\ \forall u,v \in x\,( u \neq v \imp u \cap v = \emptyset)  \\
&& \imp \bigcup x \in \FINITE.
 \end{eqnarray*}
 \end{lemma}
 
 \noindent\begin{proof}  By induction on the finite set $x$.  
 \smallskip
 {\em Base case}, $x = \emptyset$. Then $\bigcup x = \emptyset$, which is finite.
 \smallskip
 
 Induction step, $x = y \cup \{c\}$ with
 $c \not\in y$. The induction hypothesis is that if all members of 
 $y$ are finite, and any two distinct members of $y$ are disjoint,
 then $\bigcup y$ is finite.
 We have to prove that if all members of $x$ are finite
 and any two distinct members of $x$ are disjoint,
  then $\bigcup x \in \FINITE$.  Assume all members of $x$ are finite
  and any two distinct members of $x$ are disjoint.
    Since the members of $y$
 are members of $x$, all the members of $y$ are finite, and any
 two distinct members of $y$ are disjoint.  Then by the 
 induction hypothesis, $\bigcup y$ is finite.  A short argument 
 from the definitions of union and binary union proves 
 $$ \bigcup (y \cup \{ c\}) = \left( \bigcup y\right) \cup c.$$
 Since $x = y \cup \{c\}$, we have
 \begin{eqnarray}
 \bigcup x   &=& \left( \bigcup y\right) \cup c \label{eq:1057}
 \end{eqnarray}
 Now $c$ is finite, since every member of $x$ is finite and $c \in x$.
 We have   $\bigcup y \cap c = \emptyset$, since if $p$ belongs to both
 $\bigcup y$ and $c$,  then for some $w\in y$ we have $p \in w \cap c$,
 contradicting the hypothesis that any two distinct members of $x = y \cup \{c\}$
 are disjoint.  
 Then $\bigcup y \cup c$ is finite, by Lemma~\ref{lemma:union}.
 Then $\bigcup x$ is finite, by (\ref{eq:1057}).  That completes
 the induction step. \end{proof}

 \begin{lemma}\label{lemma:ssc_adjoin} Suppose $c \not \in x$. Then 
 $$ \SSC(x) \subseteq \SSC(x \cup \{c\}).$$
 \end{lemma}
 
 \noindent\begin{proof} About 30 straightforward steps, which we choose to omit here.
 \end{proof}

\begin{lemma} \label{lemma:intersectionseparable} 
Let $A$ be any set. 
Then the intersection and union 
of two separable subsets of $A$ are also separable subsets of $A$.
\end{lemma}

\noindent\begin{proof}  Let $X$ and $Y$ be two separable subsets of $A$.
Let $u \in A$.  By definition of separability, we have 
$$(u \in X \lor u \not \in X) \ \land \ (u \in Y \lor u \not\in Y).$$
I say that $X \cap Y$ is a separable subset of $A$.  To prove that,
we must prove 
\begin{eqnarray}
  u \in X \cap Y  \ \lor \ u \not \in X \cap Y. \label{eq:1660}
\end{eqnarray}
This can be proved by cases; there are four cases according to 
whether $u$ is in $X$ or not, and whether $u$ is in $Y$ or not.
In each case, (\ref{eq:1660}) is immediate.  Hence $X \cap Y$ is 
a separable subset of $A$, as claimed.  Similarly, $X \cup Y$
is a separable subset of $A$.  \end{proof}

\begin{lemma}[Finite DNS] \label{lemma:finiteDNS}
For every finite set $B$ we have
$$ \forall P\, (\forall x \in B\, ( \neg\neg\, x\in P )) \imp \neg\neg\, \forall x \in B\,( x \in P).$$
\end{lemma}

\noindent{\em Remark}. DNS stands for ``double negation shift.''  Generally 
it is not correct to move a double negation leftward through $\forall x$; but 
this lemma shows that it is OK to do so when the quantifier is bounded
by a finite set.  
\medskip

\noindent\begin{proof}  The formula of the lemma is stratified, giving
$x$ index 0, $B$ index 1, and $P$ index 1.  Therefore we may proceed
by induction on finite sets $B$.  (Notice that the statement being 
proved by induction is universally quantified over $P$--that is important
because in the induction step we need to substitute a different set for $P$;
the proof does not work with $P$ a parameter.) 
\smallskip

{\em Base case}, $B = \emptyset$.  The conclusion $\forall x \in \emptyset\, x \in P$
holds since $x \in \emptyset$ is false.  
\smallskip

{\em Induction step}.  Suppose $c \not \in B$ and $B \in \FINITE$ and 
\begin{eqnarray*}
\forall x \in B \cup \{c\}\, (\neg\neg\, x \in P). 
\end{eqnarray*}
By Lemma~\ref{lemma:empty_or_inhabited}, $B$ is empty or inhabited.
We argue by cases.  
\smallskip

Case~1, $B$ is empty.  Then $B \cup \{c\} = \{c\}$,  so we must prove
\begin{eqnarray*}
\forall x\, (x \in \{c\} \imp \neg\neg\, x \in P) \imp 
\neg\neg \forall x\, (x \in \{c\} \imp x \in P)
\end{eqnarray*}
That is equivalent to 
\begin{eqnarray*}
\forall x\, (x=c \imp \neg\neg\, x \in P) \imp \neg\neg\, \forall x\, (x = c \imp x \in P)
\end{eqnarray*}
That is, $ \neg\neg\, c \in P \imp \neg\neg\, c \in P$, which is 
logically valid.  That completes case 1.
\smallskip

Case~2, $B$ is inhabited.  Fix $u$ with $u \in B$.
Then
\begin{eqnarray*}
\forall x \in B\, (\neg\neg\, x \in P)   \\
\neg\neg\, c \in P   
\end{eqnarray*}
Since $x$ does not occur in $c \in P$ we have
\begin{eqnarray*}
\forall x \in B\, (\neg\neg\, x \in P \ \land \ \neg\neg\, c \in P) \\
\forall x \in B\,  \neg\neg\,( x \in P \ \land \   c \in P)
\end{eqnarray*}
Define $Q = \{ x : x \in P \ \land \ c \in P\}$, which is legal
since the defining formula is stratified.  Then
\begin{eqnarray}
\forall x \in B\,(\neg\neg\, x \in Q) \label{eq:1288}
\end{eqnarray}
Since $P$ is quantified in the formula being proved by induction,
we are allowed to substitute $Q$ for $P$ in the induction hypothesis;
then with (\ref{eq:1288}) we have 
\begin{eqnarray}
\neg\neg\, \forall x \in B\,( x \in Q )&&\mbox{\qquad by the induction hypothesis} \nonumber\\
\neg\neg\,  \forall x \in B\,( x \in P \ \land \ c \in P)
    &&\mbox{\qquad by the definition of $Q$  \label{eq:1296} }  
\end{eqnarray}
Now we would like to infer
\begin{eqnarray}
\neg\neg\, ((\forall x \in B\, (x \in P)) \ \land \ c \in P), \label{eq:1300}
\end{eqnarray}
which seems plausible as  $x$ does not occur in `$c \in P$'.
In fact we have the equivalence
of (\ref{eq:1296}) and~\ref{eq:1300}), since  $B$ is inhabited. 
(That was why we had to break the proof into cases according as $B$
is empty or inhabited.)  Then indeed (\ref{eq:1300}) 
follows.  By the definitions of union and unit class we have
$$
(\forall x \in B\, (x \in P)) \ \land \ c \in P \iff \forall x \in (B \cup \{c\})\,(x \in P).
$$
Applying that equivalence to (\ref{eq:1300}),
we 
have the desired conclusion,
$$ 
\neg\neg\, \forall x \in (B \cup \{c\})\,(x \in P). 
$$
That completes the induction step.  \end{proof} 

\begin{lemma}\label{lemma:notnotseparable}
Every subset of a finite set is not-not separable and not-not finite.
\end{lemma}

\noindent{\em Remark}. We already know that separable subsets of a finite
set are finite, and finite subsets of finite set are separable, but 
one cannot hope to prove every subset of a finite set is finite, because
of sets like $\{ x \in \{\emptyset\} : P\}$.   That set is finite if and 
only if $P \ \lor \ \neg P$, by Lemma~\ref{lemma:empty_or_inhabited}.
\smallskip

\noindent\begin{proof}  Let $X$ be a finite set, and $A \subseteq X$.
By Lemma~\ref{lemma:separablefinite}, if $A$ is a separable subset of $X$
then $A$ is finite.  Double-negating that implication, if $A$ is not-not
separable, then it is not-not finite.  Hence, it suffices to prove that 
not-not $A$ is a separable subset of $X$.   More formally, we 
must prove
\begin{eqnarray}
\neg\neg\, X = A \cup (X-A) \label{eq:1315}
\end{eqnarray}
 We have
\begin{eqnarray*}
  \forall t \in X\, \neg\neg\, (t \in A \ \lor \ t \not\in A) &&\mbox{\qquad by logic}\\
\neg\neg\, \forall t \in X\, t \in A \ \lor \ t \not\in A) &&\mbox{\qquad by Lemma~\ref{lemma:finiteDNS}} \\
\neg\neg\, X = A \cup (X-A) &&\mbox{\qquad by the definitions of union and difference}
\end{eqnarray*}
That is (\ref{eq:1315}).
\end{proof}
\medskip

\begin{lemma} \label{lemma:union2} Let $x\in \FINITE$ and $y \in \FINITE$.  
Then $\neg\neg\, (x \cup y \in \FINITE)$.
\end{lemma}

\noindent{\em Remark.} Lemma~\ref{lemma:union} shows the double negation 
can be dropped if $x$ and $y$ are assumed to be disjoint.  It cannot be
dropped in general, as $\{a\} \cup \{ b\} \in \FINITE$ implies $a = b \ \lor a \neq b$,
so if we could drop the double negation in this lemma, then every set would 
have decidable equality.  
\smallskip

\noindent\begin{proof} The formula is stratified, so we can prove it by 
induction on finite sets $y$, for a fixed finite set $x$. 
\smallskip

{\em Base case}, $y = \emptyset$.  We have $x \cup \emptyset = x$, which is finite by 
hypothesis.  That completes the base case.
\smallskip

{\em Induction step}.  Suppose $y \in \FINITE$,  $x \cup y \in \FINITE$,
and $c \not \in y$.  Then  I say
\begin{eqnarray}
 c \not \in x  \imp x \cup (y \cup \{c\}) \in \FINITE  \label{eq:1429}
\end{eqnarray}
To prove that:
\begin{eqnarray*}
  x \cup (y \cup \{c\}) = (x \cup y) \cup \{c\}  && \mbox{\qquad by definition of union}\\
 c \not \in x \cup y  && \mbox{\qquad since $c \not \in x$}\\
 x \cup (y \cup \{c\}) \in \FINITE 
 \end{eqnarray*}
 That completes the proof of (\ref{eq:1429}).
 \smallskip
 
 We also have
 \begin{eqnarray}
 c \in x \imp  x \cup (y \cup \{c\}) \in \FINITE  \label{eq:1442}
 \end{eqnarray}
 since $ (x \cup y ) \cup \{c\} = x \cup y \in \FINITE$.
\smallskip

We have by intuitionistic logic
$$ \neg\neg\,(c \in x \ \lor \ c \not \in x).$$
and by the induction hypothesis we have $\neg\neg\, x \cup y \in \FINITE$.
Then by (\ref{eq:1429}) and (\ref{eq:1442}), we have
$$ \neg\neg\, (x \cup y) \cup \{c\} \in \FINITE.$$
That completes the induction step.
\end{proof}

\begin{lemma} \label{lemma:oneout}  Let $X$ be a finite set and $c \in X$.  Then 
$X-\{c\}$ is finite.
\end{lemma}

\noindent\begin{proof} 
\begin{eqnarray*}
X \mbox{\ has decidable equality } && \mbox{\qquad by Lemma~\ref{lemma:finitedecidable}}\\
X - \{c\} \mbox{\ is a separable subset of $X$} && \mbox{\qquad by the definition of separable} \\
X -\{c\} \in \FINITE && \mbox{\qquad by Lemma~\ref{lemma:separablefinite}}
\end{eqnarray*}
\end{proof}
\smallskip

Before leaving this section, we shall state a technical lemma about
similarities, arising from the details of the definitions of ``maps'' and 
``similarity''.  The issue is that $f:X \imp Y$  does not require that the domain of $f$
be exactly $X$; it is allowed to be larger.  That is generally a good thing, as once we have
defined $X$ and proved it maps $X$ to $Y$, it automatically maps subsets of $X$ to $Y$.
But to be a similarity from $X$ to $Y$, the domain of $f$ must be exactly $X$ and 
the range exactly $Y$.  The following lemma is the price we must pay for allowing the 
domain of $f$ to be larger in ``maps''.  Stating it here allows us to cite it,
without revisiting this issue in future work.

\begin{lemma} \label{lemma:similarity_helper}
Suppose $f: A \imp B$ and $f$ is one-to-one and onto $B$.  Let
$R$ be the range of $f$.   Suppose $R \subseteq B$ and
the domain of $f$ is $A$.  Then $f$ is 
a similarity from $A$ to $R$.
\end{lemma} 

\noindent\begin{proof}  We omit the proof, which takes 111 lines of Lean,
because it is just a straightforward unwinding of the definitions involved.
\end{proof}

 \section{Frege cardinals}

The formula in the following definition is stratifiable, so the definition 
can be given in \INF.  Specifically, we can give $a$ index 0, $x$ and $z$
index 1, and $\kappa$ index 2.  Then $\kappa^+$ gets index 2, so 
the successor function 
$\kappa \mapsto \kappa^+$ is a function in \INF.
\begin{definition} \label{definition:cardinalsuccessor}
The successor of any set $\kappa$,  
denoted $\kappa^+$, is defined as 
$$ \kappa^+ = \{ x : \exists z,a\,(z \in \kappa \land a \not\in z \land x = z \cup \{a\})\}.$$
\end{definition}

\begin{definition} \label{definition:Fregezero}
\begin{eqnarray*}
\zero &=& \{\emptyset\} 
\end{eqnarray*}
\end{definition}

\begin{definition}\label{definition:FregeN}
The set $\FregeN$ of finite Frege cardinals is the least set 
containing $\zero = \{\emptyset\}$   and containing $\kappa^+$
whenever it contains $\kappa$ and $\kappa^+$ is inhabited.  More precisely,
$$ \kappa \in \FregeN \iff \forall w\, (\zero \in w \ \land \ 
\forall \mu\, (\mu \in w \ \land \ (\exists z\,( z \in \mu^+)) \imp \mu^+ \in w)
\imp \kappa \in w).$$
\end{definition}

\noindent{\em Remarks}.  The formula defining $\FregeN$ is stratified,
so the definition can be given in \INF.  According to that definition, if 
there were a largest finite cardinal $\kappa$, then 
 $\kappa^+$ would be the empty set,
not $\zero$, which is $\{\emptyset\}$.
So in that case, the successor of the largest finite cardinal $\kappa$ 
would not belong to $\FregeN$, which does not contain $\emptyset$.
Instead, in that case the result would be that successor does not 
map $\FregeN \to \FregeN$.  Of course $\emptyset^+ = \emptyset$, so 
once that happened, more applications of successor would do nothing more.
 Note also that in general a finite cardinal is not
a finite set; rather, the members of a finite cardinal are finite sets.
\smallskip 

 \begin{lemma} \label{lemma:finitecardinals1}
Let $\kappa \in \FregeN$ and $x \in \kappa$.  Then $x$ is a finite set.
\end{lemma}

\noindent\begin{proof}  Define
$$ Z = \{ x \in \FregeN: \forall y \in x\,(y \in \FINITE)\}.$$
The formula in the definition is stratifiable, so the definition is legal.
We will show that $Z$ is closed under the conditions defining $\FregeN$.
First, $\zero$ $ = \{ \emptyset\}$ is in $Z$, since $\emptyset$ is finite. 
To verify the second condition, assume $\kappa \in Z$ and $\kappa^+$
is inhabited; we must show 
$\kappa^+ \in Z$.  Let $u \in \kappa^+$. Then there exists 
$x \in \kappa$ and there exists $a$ such that $u = x \cup \{a\}$.
Since $\kappa \in Z$,  $x$ is finite. Then by definition of \FINITE,
$u$ is finite.  That completes the proof that $Z$ satisfies the second
condition.  Hence $\FregeN \subseteq Z$.  
\medskip

\begin{lemma}[Stratified induction] \label{lemma:induction}  Let $\phi$ be a stratified
formula (or weakly stratified with respect to $x$), so $\{x : \phi(x)\}$
exists.  Then 
$$ (\phi(\zero) \ \land \ \forall x\,(\phi(x) \ \land \exists u\,(u \in \ x^+ ) \imp \phi(x^+))) \imp \forall x\, \phi(x).$$
\end{lemma}
\noindent\begin{proof}  $Z:= \{x : \phi(x)\}$ is definable and satisfies
the closure conditions that define $\FregeN$. Therefore $\FregeN \subseteq Z$.
\end{proof}
\smallskip

{\em Remark}. When carrying out a proof by induction,
during the induction step we get to assume that $x^+$ is inhabited.
\medskip

We follow Rosser  (\cite{rosser1953}, p.~372) in defining 
cardinal numbers:  a cardinal number, or just ``cardinal'',
 is an equivalence class of 
the similarity relation $x \sim y$ of one-to-one correspondence:

\begin{definition}\label{definition:cardinal} 
The class $\NC$ of cardinal numbers is 
defined by
$$ \NC = \{ \kappa : \  \forall u\in \kappa\, \forall v\, (v \in \kappa \iff u \sim v)\}.$$
\end{definition}
\noindent
{\em Remark.} It would not do to use $\exists u$ instead of $\forall u$, since
then $\emptyset$ would not be a cardinal, but allow for that possibility.
We note that Rosser's definition requires cardinals to be inhabited.  In the 
work presented here, it makes no difference, as we work only with finite cardinals.
\medskip

The following two lemmas show that the members of $\FregeN$ are 
indeed cardinals in that sense.  

\begin{corollary} \label{lemma:cardinalsinhabited} Every 
finite cardinal is inhabited.
\end{corollary}

\noindent\begin{proof} Lemma~\ref{lemma:induction} justifies us 
in proving $\exists u\,(u\in \kappa)$ by induction on $\kappa$.
\smallskip

{\em Base case}\,: $\zero = \{ \emptyset\}$ is inhabited.  
\smallskip

{\em Induction step}\,:  Suppose $\kappa^+$ is inhabited.  Then
$\kappa^+$ is inhabited. (We do not even need to use the 
induction hypothesis.)
\end{proof}

\begin{lemma}\label{lemma:finitecardinals0}
If $\kappa \in \FregeN$ and $x \in \kappa$ and $x \sim y$,
then $y \in \kappa$.
\end{lemma}

\noindent{\em Remarks.} This lemma shows that finite cardinals
are cardinals, in the sense of equivalence classes under
similarity.   

\medskip

\noindent\begin{proof} Define
$$Z = \{ \kappa \in \FregeN: \forall x \in \kappa\, \forall y\,  
(x \sim y \imp y \in \kappa)\}.$$
That formula can be stratified, since we have already shown 
that $x \sim y$ is definable in \INF.  Therefore the 
definition of $Z$ is legal.

We will show $Z$ contains $\zero$ and is closed under Frege successor.
$Z$ contains $\zero$ since the only member of $\zero$ is the empty
set, and the only set in one-to-one correspondence with the empty set
is $\emptyset$ itself.   
\smallskip

Ad the closure under Frege successor:  Suppose $\kappa \in Z$, 
and $x \in \kappa^+$,
and $f:x \to y$ is one-to-one and onto.
Then $x = u \cup \{a\}$ for some $u \in \kappa$ and $a \not\in u$.
Let $g$ be $f$ restricted to $u$, and let $v$ be the range of $g$.
Then $g: u\to v$ is one-to-one and onto.  Since $\kappa \in Z$ 
and $u \in \kappa$, we have $v \in \kappa$.  Let $b = f(a)$. 
Then $b \not\in v$, since $f$ is one-to-one.  Then $v \cup \{b\} \in 
\kappa^+$.  
\smallskip

I say that $v \cup \{b\} = y$.   Let $p \in y$.  Then
$p = f(q)$ for some $q \in x$, since  $f$ maps $x$ onto $y$.
By Lemma~\ref{lemma:finitecardinals1}, since $x \in \kappa^+$,
 $x$ is finite. 
Since $x$ is finite, it has decidable equality by Lemma~\ref{lemma:finitedecidable}.  Therefore 
$q = a \ \lor\  q \neq a$.  If $q = a$ then $p = f(a) = b \in \{ b\}$.
If $q \neq a$ then since $q\in x = u \cup \{a\}$ and $q \in x$,
we have $q \in u$.  Then by definition of $v$,  $p = f(q) \in v$.
 Therefore
$p \in v \cup \{b\}$.  Since $p$ was an arbitrary member of $y$,
we have proved $y \subseteq v \cup \{b\}$.  But
$v \cup \{b\} \subseteq y$ is immediate, since $v \subseteq y$
and $b \in y$.  Therefore $v \cup \{b\} = y$, as claimed.

Since $v \in \kappa$, it follows that $y \in \kappa^+$ 
as desired.  Thus $Z$ is closed under Frege successor.  By the 
definition of $\FregeN$,  we have $\FregeN \subseteq Z$.  
\end{proof}

\begin{lemma} \label{lemma:finitecardinals2} Let $\kappa \in \FregeN$
and $x,y \in \kappa$.  Then $x \sim y$.
\end{lemma}

\noindent\begin{proof}  By induction on $\kappa$.  Similarity
is defined by a stratified formula, so induction is legal.
The base case is immediate as $\zero$ has only one member. 
For the induction step, let $x$ and $y$ belong to $\kappa^+$. 
Then there exist $u,v,a,b$ such that $u,v \in \kappa$ and 
$a \not \in u$ and $b \not \in v$ and $x = u \cup \{a\}$ and 
$y = v \cup \{b\}$.  By the induction hypothesis, there is a one-to-one
correspondence $g: u\to v$.  We define $f:x \to y$ 
by 
$$ f(x) = \twocases{ g(x) \mbox{\qquad if $x \in u$} }
                   { b  \mbox{\qquad \quad \ \,if $x = a$}}
$$
By Lemma~\ref{lemma:finitecardinals1}, $x$ is finite. 
Since $x$ is finite, it has decidable equality by Lemma~\ref{lemma:finitedecidable}.  Since $a \not\in u$,
$f$ is a function.  Hence the domain of $f$ is 
$x$.  By Lemma~\ref{lemma:finitecardinals1}, $y$ is finite.
Therefore  by Lemma~\ref{lemma:finitedecidable}, $y$ has decidable
equality, so the range of $f$ is $y$.  I say that 
$f$ is one-to-one.  Suppose $f(x) = f(z)$.  We 
must show $x = z$.  Since $y$ has
decidable equality, we may argue by the following cases:
\smallskip

Case~1:  $f(x) = f(z) = b$.  Then since $b\not \in v$, $x$ and $z$
are not in $u$, so $x = a$ and $z = a$.  Then $x=z$ as desired.
\smallskip

Case~2:  $f(x) = g(x)$ and $f(z) = g(z)$.  Then $g(x) = g(z)$.
Since $g$ is one-to-one, we have $x=z$ as desired.
\smallskip

Therefore $f$ is one-to-one, as claimed. Therefore $x \sim y$.
That completes the induction step.
\end{proof} 

\begin{definition}\label{definition:Nc} Following Rosser,
we define the cardinal of $x$ to be 
$$ \Nc{x} = \{ u : u \sim x\}.$$
\end{definition}
\noindent
Then the inhabited cardinals, that is, the 
inhabited members of $\NC$,  are exactly the sets of the form
$\Nc{x}$ for some $x$.

\begin{lemma} \label{lemma:xinNcx}
For all $x$,  $x \in \Nc{x}$
\end{lemma}

\noindent{lemma} By Lemma~\ref{lemma:sim}, we have $x \sim x$.
Then $x \in \Nc{x}$ by Definition~\ref{definition:Nc}.
\end{proof} 

\begin{lemma} \label{lemma:cardinalequality}
$\Nc{x} = \Nc{y}$ if and only if $x \sim y$.
\end{lemma}

\noindent\begin{proof}  By Lemma~\ref{lemma:sim}, which says
that the relation
$\sim$ is an equivalence relation.  \end{proof}

\begin{lemma} \label{lemma:Ncsuccessor}
$c \not \in x \imp \Nc{x \cup \{c\}} = \Nc{x}^+$.
\end{lemma}

\noindent\begin{proof}  By extensionality, it suffices to 
show that the two sides have the 
same members.  That is, we must show, under the assumption $c \not\in x$,
\begin{eqnarray}
u \sim x \cup \{c\} &\iff& \exists  b, v   \,( b \not \in v \ \land \  
v \sim x \ \land \ u = v \cup \{b\}). \label{eq:1236}
\end{eqnarray}

{\em Right to left}.  Suppose $b \not \in v$ and $v \sim x$ and 
$u = v \cup \{b\}$.  Let $f:v \to x$ be a similarity, and extend 
it to $g$ defined by $g = f \cup  \{\langle b, c \rangle \}$.  Then $g$ is 
a similarity from $v \cup \{b\}$ to $x \cup\{c\}$.  That completes
the right-to-left direction.
\smallskip

{\em Left to right}. Suppose $f: u \to x \cup \{c\}$ is a similarity.
Since $f$ is onto, there exists $b \in x$ with  $f(b)= c$.  
Let $v = u - \{b\}$.   Use this $b$ and $v$ on the right.  
Then 
$g = f - \{ \langle b, c \rangle \}$ is a similarity from $v$ to $x$.
It remains to show $u = v \cup \{b\} = (u-\{b\}) \cup \{b\}$.
That is,  
$$z \in u \imp z\in u \imp z \neq b \ \lor \ z = b.$$
Let $z \in u$.  Since $f$ is a similarity from $u$ to $x \cup \{c\}$,
there is a unique $y \in x \cup \{c\}$ such that $\langle z,y \rangle\in f$.
Then $z = b \iff y = c$.  Since $c \not \in x$,  and 
$y \in x \cup \{c\}$, $y = c \ \lor \ y \neq c$.  Therefore $z = b \ \lor z \neq b$, as desired.  Note that it is not necessary that $z$ have
decidable equality. 
That completes the left-to-right direction.  
\end{proof}
 \begin{lemma} \label{lemma:Nc_empty}
 $\Nc{\emptyset} = \zero$.
 \end{lemma}
 
 \noindent\begin{proof}  By Definition~\ref{definition:Fregezero}, $\zero = \{\emptyset\}$.
 By definition, $\Nc{\emptyset}$ contains exactly the sets similar to $\emptyset$.
 By Lemma~\ref{lemma:similar_to_empty2},  $\emptyset$ is the only set similar to 
 $\emptyset$.  Therefore $\Nc{\emptyset} = \{\emptyset\}$.  Then $\Nc{\emptyset} = \zero$ since
 both are equal to $\{\emptyset\}$.  \end{proof}

\begin{lemma}\label{lemma:successorinhabited}
  For  every set $\kappa$, if $\kappa^+$ is inhabited, then  
$\kappa^+$ contains an 
inhabited set, and  every member of $\kappa^+$ is inhabited.  
\end{lemma}

\noindent{\em Remark}. Note that $\kappa$ is not assumed to 
be a finite cardinal, or even a cardinal.  Successor cannot take
the value $\zero = \{ \emptyset\}$ on any set.   
\medskip

\noindent\begin{proof}   By definition the members of
$\kappa^+$ are exactly the sets of the form $x \cup \{a\}$ with 
$x \in \kappa$ and $a \not\in x$.  (That is true whether or not
there are any such members.)  But if $\kappa^+$ is inhabited, then
there is at least one such member, and each such member $ x \cup \{a\}$ is 
 inhabited, since it contains $a$.  
 \end{proof}
 
\begin{lemma} \label{lemma:Fregesuccessoromits0} 
  Frege successor does not take the value $\zero$ on  any set at all:
$\forall x\, (x^+ \neq \zero)$.
\end{lemma}
 
\noindent{\em Remark}.  This does not depend on the finiteness or 
not-finiteness of $\FregeN$.  If $\FregeN$ is finite then eventually 
$\kappa^+$ is $\emptyset$, rather than $\zero$, which is $\{\emptyset\}$,
so even in that case $\zero$ does not occur as a successor.
\medskip
 
\noindent\begin{proof} If $\kappa^+ =  \{\emptyset\}$ then $\kappa^+$
is inhabited, but contains no inhabited set, contradicting Lemma~\ref{lemma:successorinhabited}. \end{proof}

\begin{lemma} \label{lemma:nonzeroissuccessor}   
  Every finite cardinal is either equal to $\zero$ or is the successor of an element of $\F$.
\end{lemma} 

\noindent\begin{proof} 
The set $Z = \{ \kappa \in \FregeN : 
\kappa = \zero \ \lor \ \exists \mu\, ( \mu \in \F \ \land \ \kappa = \mu^+)\}$ is definable
in \INF, since its defining formula is stratified.  $Z$ contains
$\zero$ and is closed under successor.  Therefore, by definition of 
$\FregeN$, $\FregeN \subseteq Z$.  \end{proof} 

\begin{lemma}\label{lemma:zeroF} $\zero \in \F$.
\end{lemma}

\noindent\begin{proof} Let $W$ be one of the sets whose intersection defines $\FregeN$, i.e.,
$W$ contains $\zero$ and is closed under inhabited successor.  Then 
$W$ contains $\zero$.  Since $W$ was arbitrary, $\zero \in \F$. 
\end{proof}

\begin{lemma} \label{lemma:successorF}
  $\FregeN$ is closed under inhabited successor.
\end{lemma}

\noindent\begin{proof}
  Suppose $\kappa \in \FregeN$ and $\kappa^+$ is inhabited.
Let $W$ be one of the sets whose intersection defines $\FregeN$, i.e.,
$W$ contains $\zero$ and is closed under inhabited successor.  By 
induction on $\kappa$, we can prove $\kappa \in W$.
Since $W$ is closed under inhabited successor,  and $\kappa \in W$,
and $\kappa^+$ is inhabited, we have $\kappa^+ \in W$.  
Since $\FregeN$ is the intersection of all such sets $W$, and
$\kappa^+$ belongs to every such $W$, we have $\kappa^+ \in \FregeN$
as desired.  \end{proof}

\begin{lemma} \label{lemma:finitecardinals3}
The cardinal of a finite set is a finite cardinal.
That is, 
$$ \forall x \in \FINITE\,( \Nc{x} \in \FregeN).$$
\end{lemma}

\noindent\begin{proof} The formula to be proved is 
  stratified, so we can prove it by induction
on finite sets.  
\smallskip

{\em Base case}\,: By Lemma~\ref{lemma:Nc_empty}, $\Nc{\emptyset} = \zero$.
By Lemma~\ref{lemma:successorF}, $\zero \in \FregeN$.
\smallskip

{\em Induction step}\,:  Let $x \in \FINITE$ and $c \not \in x$. Consider 
$\Nc{{x \cup \{c\}}}$,  which by Lemma~\ref{lemma:Ncsuccessor}  is 
$\Nc{x}^+$.  By the induction hypothesis, $\Nc{x} \in \FregeN$.
By definition of $\FregeN$, $\Nc{x}^+ \in \FregeN$. That completes
the induction step.  \end{proof} 

\begin{lemma} \label{lemma:Finhabited}  Every member of $\F$ is inhabited.
\end{lemma}

\noindent\begin{proof}  By induction we prove 
$$\forall m\, (m \in \F \imp \exists u\, (u \in m)).$$
The formula is stratified, giving $u$ index 0 and $m$ index 1.
For the base case, $\zero = \{\emptyset\}$ by definition, so $\zero$ 
is inhabited.  For the induction step, we always suppose $m^+$ is 
inhabited, so there is nothing more to prove.  
\end{proof}
  
To put the proof directly: the set of inhabited 
members of $\F$ contains $\zero$ and is closed under inhabited successor,
so it contains $\F$.

\begin{lemma} \label{lemma:similar_to_finite} A set similar to a finite set is finite.
\end{lemma}

\noindent\begin{proof}  Let $a$ be finite and $a \sim b$. Let 
$\kappa = \Nc{a}$.  Then
\begin{eqnarray*}
\kappa \in \FregeN &&\mbox{\qquad by Lemma~\ref{lemma:finitecardinals3}}\\
a \sim a  && \mbox{\qquad by Lemma~\ref{lemma:sim}}\\
a \in  \kappa  && \mbox{\qquad by definition of $\Nc{a}$}\\
b \in \kappa   && \mbox{\qquad by Lemma~\ref{lemma:finitecardinals0}}\\
b \in \FINITE  && \mbox{\qquad by Lemma~\ref{lemma:finitecardinals1}}
\end{eqnarray*}
\end{proof}

\begin{lemma} \label{lemma:cardinalsdisjoint} \ \\
(i) If two finite cardinals have a common member, then they are equal.

\noindent
(ii) Two distinct finite cardinals are disjoint.
\end{lemma}

\noindent\begin{proof} Part (ii) is the contrapositive of (i),
so it suffices to prove (i).  Let $\kappa$ and $\mu$ belong to $\FregeN$.
Suppose $x$ belongs to both $\kappa$ and $\mu$.  We must show 
$\kappa = \mu$.  By extensionality, it suffices to show that 
$\kappa$ and $\mu$ have the same members.  Let $y \in \kappa$.
Then by Lemma~\ref{lemma:finitecardinals2}, $y \sim x$. 
By Lemma~\ref{lemma:finitecardinals0}, $y \in \mu$.  Therefore
$\kappa \subseteq \mu$.  Similarly $\mu \subseteq \kappa$.
\end{proof}

\begin{lemma} \label{lemma:finitecardinals4}
Let $x$ and $y$ be finite sets.  Then 
$$ x \sim y \imp \Nc{x} = \Nc{y}.$$
\end{lemma}

\noindent\begin{proof}  Assume $x\in \FINITE$ and $y \in \FINITE$ and $x \sim y$. Then
\begin{eqnarray*}
x \in \Nc{x}  && \mbox{\qquad by Lemma~\ref{lemma:xinNcx}}\\
y \in \Nc{y}  && \mbox{\qquad by Lemma~\ref{lemma:xinNcx}}\\
\Nc{x} \in \F && \mbox{\qquad by Lemma~\ref{lemma:finitecardinals3}}\\
\Nc{y} \in \F && \mbox{\qquad by Lemma~\ref{lemma:finitecardinals3}}\\
y \in \Nc{x}  && \mbox{\qquad by Lemma~\ref{lemma:finitecardinals0}}\\
\Nc{x} = \Nc{y}  && \mbox{\qquad by Lemma~\ref{lemma:cardinalsdisjoint}}
\end{eqnarray*}
\end{proof}

\section{Order on the cardinals} 
In this section, $\kappa, \mu$, and $\lambda$ will always
be cardinals.  We start with Rosser's classical definition
(which is not the one we use). 

\begin{definition} [Rosser] \label{definition:orderRosser} 
\begin{eqnarray*}
\kappa \le \mu &:=&  \exists a,b\,(a \in \kappa \ \land \ b \in \mu \ \land\ a \subseteq b)\\
\kappa < \mu &:=& \kappa \le \mu \ \land \ \kappa \neq \mu.
\end{eqnarray*}
\end{definition}

For constructive use, we need to add the requirement
$b = a \cup (b-a)$, which says that 
$b$ is a separable subset of $a$.
Classically, every subset is 
separable, so the definition is classically equivalent to Rosser's.

\begin{definition}  \label{definition:order} 
For cardinals $\kappa$ and $\mu$:
\begin{eqnarray*}
 \kappa \le \mu &:=&  \exists a,b\,(a \in \kappa \ \land \ b \in \mu 
\ \land \ a \subseteq b\ \land\ b = a \cup (b-a)) \\
  \kappa < \mu &:=& \kappa \le \mu \ \land \ \kappa \neq \mu.
\end{eqnarray*}
\end{definition}

\begin{definition} \label{definition:image}
The {\bf image of $a$ under $f$}, written $f``a$,  is defined by
$$  f`` a := Range(f \cap (a \times \V)).$$
\end{definition}
If $f$ is a function then $f``a$ is the set of values $f(x)$ for $x \in a$.

\begin{lemma} \label{lemma:separable_similarity}
The image of a separable subset under a similarity is a separable
subset.  More precisely, let $f: b \to c$ be a similarity and
suppose $b = a \cup (b-a)$.  Let $e = f``a$ be the image of $a$
under $f$.  Then $c = e \cup (c-e)$.
\end{lemma}

\noindent\begin{proof} 
  We have 
 \begin{eqnarray}
  e \cup (c-e)  &\subseteq& c  \label{eq:1366}
\end{eqnarray}
since $e \subseteq c$ and $c-e \subseteq c$.  We have 
\begin{eqnarray}
 c &\subseteq& e \cup (c-e) \label{eq:1370}
\end{eqnarray}
since if $q \in c$ then $q = f(p)$ for some $p \in b$, 
and $p \in a \ \lor \ p \in b - a$, since $b = a \cup (b-a)$,
and if $p \in a$ then $q \in e$, while if $p \in b-a$ then $q \in c-e$.
Combining (\ref{eq:1366}) and (\ref{eq:1370}),  we have
 $ c = e \cup (c-e)$ as desired.  \end{proof}
 
\begin{lemma} \label{lemma:similarity_image}
Let $f:a \to b$ be a similarity, and let $x \subseteq a$.  
Let $g$ be $f$ restricted to $x$.  Then $g: x \to f``x$ is 
a similarity.
\end{lemma}

\noindent\begin{proof}  Straightforward; requires about 75 
inferences that we choose to omit here. \end{proof}
 
\begin{lemma} \label{lemma:le_transitive} The ordering relation $\le$
is  transitive on $\FregeN$. 
\end{lemma}

\noindent\begin{proof} 
  Suppose $\kappa \le \lambda $ and $\lambda \le \mu$.
We must show $\kappa \le \mu$.  Since $\kappa < \lambda$ and 
$\lambda < \mu$, there exist $a \in \kappa$, $b,c \in \lambda$,
and $d \in \mu$ such that $a \subseteq b$ and $c \subseteq d$,
and $b = a \cup (b-a)$, and $d = c \cup (d-c)$. 
By Lemma~\ref{lemma:finitecardinals2},
 $b \sim c$,  since both belong to $\lambda$.  Let $f:b \to c$
be one-to-one and onto. Let $e = f``a$.  Then $e \subseteq c$
and $a \sim e$.  So $e \in \kappa$, by Lemma~\ref{lemma:finitecardinals0}.
Then $e \subseteq d$. 
By Lemma~\ref{lemma:separable_similarity} we have
\begin{eqnarray}
c &=& e \cup (c-e)    \label{eq:1375}
\end{eqnarray}
Now I say that  $d = e \cup (d-e)$.
\begin{eqnarray*}
e \cup (d-e) &=& e \cup ((c \cup (d-c)) - e)\mbox{\hskip1.2cm since $d = c \cup (d-c)$}\\
             &=& e \cup (c-e) \cup ((d-c) - e)
                   \mbox{\quad\ \ since $(p \cup q) -r = (p-r) \cup (q-r)$} \\
             &=&  c \cup ((d-c)-e) \mbox{\hskip 2cm by (\ref{eq:1375}) }\\
             &=&  c \cup (d-c)  \mbox{\hskip2.9cm  since $e \subseteq c$} \\
             &=& d   \mbox{\hskip4.3cm since $d = c \cup (d-c)$} 
\end{eqnarray*}
as desired. 
Then $\kappa \le \mu$ as desired.  
 \end{proof}

\begin{lemma} \label{lemma:lessthan2} For finite cardinals $\kappa$ and $\mu$,
$$ \kappa < \mu \iff 
 \exists x,y\,(x \in \kappa \ \land \ y \in \mu
  \ \land \ x \subset y \ \land\  y = x \cup (y-x)).$$
\end{lemma}

\noindent\begin{proof}  {\em Left to right}.
Suppose $\kappa < \mu$.  Then by definition of $<$,
$\kappa \le \mu$ and $\kappa \neq \mu$.  By definition of $\le$, 
 there exist $x$ and $y$ with $x \in \kappa$, $y \in \mu$, and 
$x \subseteq y$ and $y = x \cup (y-x)$.  
By Lemma~\ref{lemma:cardinalsdisjoint},
which applies because $\kappa \neq \mu$,  we have $x \neq y$. 
Therefore $x \subset y$ as desired.  That completes the proof of the
left-to-right implication. 
\smallskip

{\em Right to left}.   Suppose $x \in \kappa$ and $y \in \mu$ and $x \subset y$
and $y = x \cup (y-x)$. 
Then $\kappa \le \mu$ by definition.  We must show $\kappa \neq \mu$.
If $\kappa = \mu$ then $y \sim x$, by 
Lemma~\ref{lemma:finitecardinals2}.  Then $y$
 is similar to a proper subset of $y$, namely $x$. 
 Since $y \in \mu$ and $\mu \in \FregeN$, by Lemma~\ref{lemma:finitecardinals1},
 $y$ is finite. 
Since $y$ is similar to a proper subset of itself (namely $x$),
 Theorem~\ref{theorem:infiniteimpliesnotfinite} implies that
  $y$ is not finite, which is a contradiction.
\end{proof}

 \begin{lemma}\label{lemma:le2}  Let $\kappa,\mu \in \FregeN$, 
 with $\mu$  inhabited. Then
$$\kappa \le \mu \iff \forall b \in \mu\, \exists a \in \kappa\,
( a \subseteq b \ \land \ b = a \cup (b-a)).$$
\end{lemma}

\noindent\begin{proof}  Left-to-right.  Suppose $\kappa \le \mu$.
Then by definition of $\le$, 
there exist $x\in \kappa$ and $y \in \mu$ with $x \subseteq y$
and 
\begin{eqnarray}
y &=& x \cup (y-x) \label{eq:1454}
\end{eqnarray} 
Let $b \in \mu$; we must show  there exists $a \in \kappa$ 
with $a \subseteq b$ and $b = a \cup (b-a)$.
 
 We have $b \sim y$ by Lemma~\ref{lemma:finitecardinals2}.
So $y \sim b$.   Let $f:y \to b$ be one-to-one and onto. 
Let $a = f``(x)$.  Then $a \subseteq b$ and $x \sim a$.
 By Lemma~\ref{lemma:finitecardinals0},  $a \in \kappa$.
By Lemma~\ref{lemma:separable_similarity} and (\ref{eq:1454}),
we have  $b = a \cup (b-a)$.  
That completes the proof of the left-to-right implication.
\smallskip

Right-to-left. Suppose $\forall b \in \mu\,  \exists a \in \kappa\,
(a \subseteq b \ \land \  b = a \cup (b-a))$.   Since $\mu$ is a 
inhabited, there exists $b \in \mu$.  Then $\exists a \in \kappa\,
(a \subseteq b\ \land \ b = a \cup (b-a))$.  \end{proof}

\begin{lemma} \label{lemma:cardinalpredecessor}
Suppose $\kappa  \in \FregeN$ and $x \in \kappa^+$ 
and $c \in x$.  Then $x-\{c\} \in \kappa$.
\end{lemma}

\noindent{\em Remark}.  We will use this in the proof
that successor is one-to-one, so we cannot use that 
fact to prove this lemma.
\medskip

\noindent\begin{proof}  Since $x \in \kappa^+$, 
there exists $z \in \kappa$ and $a \not\in z$
such that $x = z \cup \{a\}$. Since $c\in x$, 
we have $c\in z \ \lor \ c = a$.  If $c = a$ 
then $z = x-\{c\} \in \kappa$ and we are done.
Therefore we may assume $c  \in z$ and $c \neq a$.

Since $a \neq c$ we have 
\begin{eqnarray}
(z \cup \{a\}) - \{c\} = (z - \{c\}) \cup \{a\} \label{eq:1487}
\end{eqnarray}
  Since $x \in \kappa^+$,
$x$ is finite, by Lemma~\ref{lemma:finitecardinals1}.
By Lemma~\ref{lemma:finitedecidable}, $x$ has decidable
equality.  Then
\begin{eqnarray*}
z &\sim& (z -\{c\}) \cup \{a\} \mbox{\qquad by Lemma~\ref{lemma:swap_similarity}}\\  
&=& (z \cup \{a\}) - \{c\}  \mbox{\qquad by (\ref{eq:1487})} 
\end{eqnarray*}
Then by Lemma~\ref{lemma:finitecardinals0} and the fact that $z \in \kappa$,  
we have \begin{eqnarray}
(z \cup \{a\}) - \{c\} \in \kappa \label{eq:1494} 
\end{eqnarray}
 Since $z \cup \{a\} = x$, that implies $x-\{c\} \in \kappa$,
 which is the conclusion of the lemma.  
\end{proof}

\begin{lemma}\label{lemma:ordersuccessor}
For finite cardinals $\kappa$ and $\mu$, if   
 $\mu^+$ is inhabited, we 
have 
$$ \kappa \le \mu \iff \kappa^+ \le \mu^+.$$
\end{lemma}

\noindent\begin{proof} 
{\em Left to right}.  Suppose
$\kappa \le \mu$. 
 Since $\mu^+$ is inhabited, 
there is some $y \in \mu$ and some $c \not\in y$, so 
$y \cup \{c\} \in \mu^+$.  By Lemma~\ref{lemma:le2},
there is a separable subset $x \subseteq y$ with $x \in \kappa$.
Then 
$x \cup \{c\} \in \kappa^+$ and $x \cup \{c\} \subset y \cup \{c\}$.
We have to show that 
\begin{eqnarray}
y \cup \{c\} &=& (x \cup \{c\})\cup( y\cup \{c\} - (x \cup \{c\})).
\label{eq:1077}
\end{eqnarray}
 Left-to-right of (\ref{eq:1077}):  Suppose $u \in y \cup \{c\}$.
 Then $u \in y$ or $u = c$. If $u=c$ then $u \in x \cup \{c\}$,
 so $u$ belongs to the right side of (\ref{eq:1077}). Now 
 $y \cup \{c \}$ is
  finite (by Lemma~\ref{lemma:finitecardinals1}),
 and hence has decidable equality by Lemma~\ref{lemma:finitedecidable}.
 Therefore $u = c \ \lor \ u \neq c$; so we can assume $u \neq c$.
 If $u \in y$ then, since $y = x \cup (y-x)$, $u \in x \ \lor \ u \not\in x$.
 If $u \in x$ then $u \in x \cup \{c\}$ and hence $u$ belongs
 to the right side of (\ref{eq:1077}).  If $u \not\in x$ then 
 $u \in  y\cup \{c\} - (x \cup \{c\})$, and hence $u$ belongs to 
 the right side of (\ref{eq:1077}).  That completes the proof 
 of the left-to-right direction of (\ref{eq:1077}).
 \smallskip
 
 Right-to-left of (\ref{eq:1077}).  Since $x \subseteq y$ we have 
 $$x \cup \{c\} \subseteq y \cup \{c\}$$
 and
 $$ y\cup \{c\} - (x \cup \{c\}) \subseteq y \cup \{c\}.$$
 Hence the right side of (\ref{eq:1077}) is a subset of the left side.
 That completes the proof of~(\ref{eq:1077}).
 \smallskip
 
 Therefore $\kappa^+ \le \mu^+$.  That completes the proof of the
left-to-right direction of the lemma.
\smallskip

{\em Right to left}.  Suppose $\kappa^+ \le \mu^+$.  Then 
there exist $x \in \kappa^+$ and $y \in \mu^+$ with $x \subseteq y$
and $y = x \cup (y-x)$. 
By Lemma~\ref{lemma:successorinhabited}, $x$ is 
inhabited, so 
 there exists $c \in x$. Since $x \subseteq y$, also $c \in y$. 
 Then by Lemma~\ref{lemma:cardinalpredecessor},
 $x-\{c\} \in \kappa$
and $y-\{c\} \in \mu$. 
Since $y \in \mu^+$, $y$ is finite, by Lemma~\ref{lemma:finitecardinals1}.
By Lemma~\ref{lemma:finitedecidable}, $y$ has decidable equality. 
Then 
\begin{eqnarray}
u \in y \imp u = c \ \lor \ u \neq c \label{eq:1567}
\end{eqnarray}
Since $y = x \cup (y-x)$, we have 
\begin{eqnarray}
 u \in y \imp u \in x \ \lor \ u \not\in x  \label{eq:1571}
 \end{eqnarray}
Then by (\ref{eq:1567}) and (\ref{eq:1571}), we have 
\begin{eqnarray}
 u \in y \imp u \in (x - \{c\}) \ \lor \ u \not \in (x - \{c\}).\label{eq:1575}
\end{eqnarray} 
It follows from (\ref{eq:1575}) that 
\begin{eqnarray*}
y-\{c\} &=& ((y-\{c\}) -(x-\{c\})) \cup (x-\{c\})  
\end{eqnarray*}
Therefore $\kappa \le \mu$.  \end{proof}

\begin{lemma} \label{lemma:successoroneone}
For $\lambda$ and $\mu$ in $\FregeN$, if $\lambda^+$ and $\mu^+$ are inhabited,
then
$$ \lambda = \mu \iff \lambda^+ = \mu^+.$$
\end{lemma}

\noindent\begin{proof}  Left to right is immediate.  We take up the 
right to left implication.   
Suppose $\kappa^+ = \mu^+$. By Lemma~\ref{lemma:cardinalsdisjoint},
it suffices to show that $\kappa \cap \mu$ is inhabited.  Since
$\kappa^+$ is inhabited, there exists $y \in \kappa^+$.  By definition
of successor, $y$ has the form $y = x \cup \{a\}$ for some $x \in \kappa$
and $a \not\in x$.  
 We will prove $x \in \mu$.  
  Since $\mu^+ = \kappa^+$ we 
have $x \cup \{a\} \in \mu^+$.  
 Then by Lemma~\ref{lemma:cardinalpredecessor},
 $x \cup \{a\}-\{a\} \in \mu$.
Since $x \cup \{a\} \in \mu^+$, $x \cup \{a\}$ is finite, by Lemma~\ref{lemma:finitecardinals1}.
By Lemma~\ref{lemma:finitedecidable}, $x \cup \{a\}$ has decidable equality. 
Then $x \cup \{a\}-\{a\} = x$,  so $x \in \mu$.  Then 
$x \in \kappa \cap \mu$ as claimed. 
\end{proof}

\begin{lemma} \label{lemma:difference_nonempty}
Let $x$ be a separable subset of $y$, that is, $x \subseteq y$ and 
$y = x \cup (y-x)$.  Then $y- x = \emptyset\iff y = x$.
\end{lemma}

\noindent\begin{proof} Suppose $x \subseteq y$ and 
$y = x \cup (y-x)$.  {\em Left to right}. suppose $y-x = \emptyset$; we 
must show $y=x$.  If $u \in x$ then by $y = x \cup (y-x)$ we have 
$u \in y$. Conversely, if $u \in y$ then $u \in x \ \lor \ u \not \in x$.
If $u \in x$ we are done.  If $u \not \in x$ then $u \in y - x$, so $u \in y$.
That completes the left-to-right direction.  
{\em Right to left}.  Suppose $y=x$.  Then $y-x = x-x = \emptyset$.  
\end{proof}

\begin{lemma}\label{lemma:successorstrict}
For finite cardinals $\kappa$ and $\mu$, if $\kappa^+$ 
and $\mu^+$ are inhabited we 
$$ \kappa < \mu \iff \kappa^+ < \mu^+.$$
\end{lemma}

\noindent\begin{proof}
 Left-to-right.  Suppose $\kappa < \mu$.  By definition that 
means $\kappa \le \mu$ and $\kappa \neq \mu$.  By Lemma~\ref{lemma:ordersuccessor},
$\kappa^+ \le \mu^+$. We have to show $\kappa^+ \neq \mu^+$.
Suppose $\kappa^+ = \mu^+$.   Since $\mu^+$ is inhabited,
 there is an element $y \cup \{c\}$ of $\mu^+$ 
with $y \in \mu$ and $c \not \in y$.  Since $\kappa^+ = \mu^+$,
we also have $y \cup \{c\} \in \kappa^+$. 
Since $y \in \mu$, by Lemma~\ref{lemma:finitecardinals1}, $y$ is finite.
Since $\mu^+$ is inhabited, $\mu$ is also inhabited.  
 Since $\kappa < \mu$, by Lemma~\ref{lemma:le2},
there exists a separable subset $x$ of $y$ with 
$x \in \kappa$. By Lemma~\ref{lemma:finitecardinals1}, 
$x$ is finite.  By Lemma~\ref{lemma:finitedif}, $y-x$ is finite. 
 Since $\kappa \neq \mu$,  we have $x \neq y$,
by Lemma~\ref{lemma:cardinalsdisjoint}. 
Then, since $x$ is a separable subset of $y$,
  $y-x$ is not empty, by Lemma~\ref{lemma:difference_nonempty}.
   Since it is finite, 
by Lemma~\ref{lemma:empty_or_inhabited}, $y-x$ 
is inhabited.   Hence there exists some $b\in y$ with $b \not\in x$.
Then $x \cup \{b\} \in \kappa^+$.  Then $x \cup \{b\}$ and 
$y \cup \{c\}$ both belong to $\kappa^+$.  

Note that $x \cup\{b\}$
and $y \cup \{c\}$ are finite (by Lemma~\ref{lemma:finitecardinals1}), 
and hence have decidable equality (by Lemma~\ref{lemma:finitedecidable}).
Hence $y = (y \cup \{c\}) - \{c\}$; then by Lemma~\ref{lemma:cardinalpredecessor}
we have $y \in \kappa$.  But from the start we had $y \in \mu$.  
Then by Lemma~\ref{lemma:cardinalsdisjoint}, we have $\kappa = \mu$,
contradicting the hypothesis $\kappa < \mu$.  Hence the assumption 
$\kappa^+ = \mu^+$ has led to a contradiction.  Hence $\kappa^+ < \mu^+$.
That completes the proof of the left-to-right direction of the lemma. 
  
Right-to-left:  Suppose $\kappa^+ < \mu^+$.  Then $\kappa^+ \le \mu^+$
and $\kappa^+ \neq \mu^+$.  By Lemma~\ref{lemma:ordersuccessor},
 $\kappa \le \mu$, and since
successor is a function, $\kappa \neq \mu$.  
\end{proof}

\begin{definition} \label{definition:Fregeonetwo}
We define names for the first few integers (repeating
the definition of $\zero$, which has already been given).
\begin{eqnarray*}
\zero &=& \{\emptyset\} \\
\one  &=& \zero^+ \\
\two  &=& \one^+ \\
\three &=& \two^+\\
\four &=& \three^+ 
\end{eqnarray*} 
\end{definition} 

\begin{lemma} \label{lemma:oneF} $\one \in \F$.
\end{lemma}

\noindent\begin{proof}  
\begin{eqnarray*}
\zero \in \F &&\mbox{\qquad by Lemma~\ref{lemma:zeroF}}\\
\one = \zero^+ &&\mbox{\qquad by the definition of $\one$}\\
 \emptyset \in \zero && \mbox{\qquad by the definition of $\zero$}\\
\zero \not\in \zero && \mbox{\qquad since $\zero = \{\emptyset\}$ and $\zero \neq \emptyset$}\\
\emptyset \cup \{\zero\} \in \zero^+ &&\mbox{\qquad by definition of successor}\\
\exists u\, (u \in \one)  && \mbox{\qquad since $\one = \zero^+$}\\
\one \in \F &&\mbox{\qquad by Lemma~\ref{lemma:successorF}}
\end{eqnarray*}
\end{proof} 

\begin{lemma} \label{lemma:zero_or_not_zero}  We have
$$\forall \kappa \in \FregeN\,(\kappa = \zero \ \lor \ \kappa \neq \zero).$$
\end{lemma}

\noindent\begin{proof}  By induction on $\kappa$. More explicitly, 
define 
$$W := \FregeN \cap ((\FregeN - \{ \zero \}) \cup \{ \zero \}).$$
  We will show that 
$W$ satisfies the conditions defining $\FregeN$.  Specifically 
$0 \in W$ (which is immediate from the definitions of $W$ and union),
and $W$ is closed under (inhabited) Frege successor.  Suppose $\kappa \in W$
and $\kappa^+$ is inhabited. 
We have to show $\kappa^+ \in W$.  By Lemma~\ref{lemma:Fregesuccessoromits0},
$\kappa^+ \neq \zero$.  By definition of $W$, $\kappa \in \FregeN$.
By definition of $\FregeN$, $\kappa^+ \in \FregeN$;
therefore $\kappa^+ \in \FregeN - \{\zero\}$.  Therefore $\kappa^+ \in W$,
as claimed.
\smallskip

Then by definition of $\FregeN$ (or, if you prefer, ``by induction on $\kappa$''),
$\FregeN \subseteq W$.  Then by the definition of union, 
$\kappa \in \FregeN \imp \kappa = \zero \ \lor\  \kappa \neq \zero$.
\end{proof}

\begin{theorem} \label{theorem:finitetrichotomy}
For finite cardinals $\kappa$ and $\mu$, we have
$$\kappa < \mu \lor \kappa= \mu \lor \mu < \kappa$$
and 
$$ \neg\, (\kappa < \mu \ \land \ \mu < \kappa).$$
 
\end{theorem}

\noindent\begin{proof}  We prove by induction on $\kappa$
that for all $\mu$ we have the assertion in the statement
of the lemma. Lemma~\ref{lemma:induction} justifies this 
method of proof.   The formula is 
stratified since the relation $x < y$ is definable.
\smallskip

{\em Base case}\,:  
 We have to prove $$\zero < \mu \ \lor \ \zero = \mu \ \lor\ \mu < \zero$$
and exactly one of the three holds.   If $\mu \le \zero$,
then we would have $x \in \mu$ and $x$ a separable
subset of $y$ and $y \in \zero$;  but the only member of \zero\ is $\emptyset$,
so $x = y = \emptyset$.  Then $\emptyset \in \mu$ and $\emptyset \in \zero$,
so by Lemma~\ref{lemma:cardinalsdisjoint},  $\mu = \zero$.  Thus $\mu < \zero$ is impossible
and $\mu \le  \zero$ if and only if $\mu = \zero$.  If $\mu \in \FregeN$ then 
by Lemma~\ref{lemma:zero_or_not_zero},
$\mu = \zero \lor \mu \neq \zero$; 
and if $\mu \neq \zero$ then $\zero < \mu$, since $\emptyset$ is a separable
subset of any $x \in \mu$.  
\smallskip

{\em Induction step}\,:  Suppose $\kappa^+$ is inhabited. 
We have to prove 
\begin{eqnarray}
\kappa^+ < \mu \ \lor \ 
\kappa^+ =  \mu \ \lor \ \mu < \kappa^+ \label{eq:1698}
\end{eqnarray}  
By Lemma~\ref{lemma:zero_or_not_zero}, we have $\mu = \zero \ \lor \ \mu \neq \zero$.
If $\mu = \zero$, we
are done by the base case.  If $\mu \neq \zero$, then 
by Lemma~\ref{lemma:nonzeroissuccessor}, $\mu = \lambda^+$
for some $\lambda \in \FregeN$.  By Corollary~\ref{lemma:cardinalsinhabited},
$\lambda^+$ is inhabited.   We have to prove
\begin{eqnarray}
 \kappa^+ < \lambda^+ \ \lor \kappa^+ = \mu^+ \ \lor \ \mu^+ < \lambda^+.
 \label{eq:1154}
\end{eqnarray}
By the induction hypothesis we have 
$$ \kappa < \lambda \ \lor \kappa = \mu \ \lor \ \mu < \lambda.$$
and exactly one of the three holds.
By Lemma~\ref{lemma:successorstrict} and Lemma~\ref{lemma:successoroneone},
 each disjunct is equivalent 
to one of the disjuncts of (\ref{eq:1154}).  That completes
the induction step.  \end{proof}

\begin{corollary} \label{lemma:FregeNdecidable}
$\FregeN$ has decidable equality.  Precisely, 
$$\forall \kappa, \mu \in \FregeN\, (\kappa = \mu \ \lor \ \kappa \neq \mu).$$
\end{corollary}

\noindent\begin{proof} Let $\kappa, \mu \in \FregeN$.  
We must show $\kappa = \mu \ \lor \kappa \neq \mu$.
By Theorem~\ref{theorem:finitetrichotomy}, we have
 $\kappa < \mu$ or $\kappa = \mu$ or $\mu < \kappa$,
 and exactly one of the disjuncts holds. Therefore 
 $\kappa \neq \mu$ is equivalent to $\kappa < \mu \ \lor \ \mu < \kappa$.
 \end{proof}

\begin{lemma} \label{lemma:le_reflexive} For all $\kappa \in \FregeN$,
we have $\kappa \le \kappa$.
\end{lemma}

\noindent\begin{proof}  Suppose $\kappa\in \FregeN$. By Corollary~\ref{lemma:cardinalsinhabited},
$\kappa$ is inhabited.  Let $a \in \kappa$.  Since $a$ is a separable subset
of $a$, we have $\kappa \le \kappa$ by the definition of $\le$.  
\end{proof}

\begin{lemma} \label{lemma:letolessthan} For $\kappa,\mu \in \FregeN$
we have $$ \kappa \le \mu \iff \kappa < \mu \ \lor \ \kappa = \mu.$$
\end{lemma}

\noindent\begin{proof} Suppose $\kappa, \mu \in \FregeN$. 
By Theorem~\ref{theorem:finitetrichotomy} we have 
$\kappa < \mu \ \lor \ \kappa = \mu \ \lor\ \mu < \kappa$,
and exactly one of the three disjuncts holds.
\smallskip

{\em Left to right}. Suppose $\kappa \le \mu$. By Definition~\ref{definition:order},
there exist $a$ and $b$ with $a \in \kappa$, $b \in \mu$, $a \subseteq b$,
and $b = a \cup (b-a)$. By Lemma~\ref{lemma:finitecardinals1}, $a$ and 
$b$ are finite.  By Lemma~\ref{lemma:finitedif}, $b-a$ is finite. 
By Lemma~\ref{lemma:empty_or_inhabited}, $b-a$ is empty or inhabited. 
\smallskip

Case 1, $b-a = \emptyset$. 
I say $b = a$. By extensionality, it suffices to prove 
$t \in b \iff t \in a$.  {\em Left to right}.assume $t \in b$.
  Since $b = a \cup (b-a)$ we have 
$t \in a \ \lor \ t \in b-a$.  But $t \not \in b-a$, since $b-a = \emptyset$.
Therefore $t \in a$.  {\em Right to left}.  assume $t \in a$.  Since 
$a \subseteq b$ we have $t \in b$. Therefore $b=a$ as claimed.
\smallskip

Then $a \in \kappa \cap \mu$.  Then 
by Corollary~\ref{lemma:cardinalsdisjoint}, $\kappa = \mu$.
That completes Case~1.
\smallskip

Case 2, $b-a$ is inhabited.  Then $a$ is a proper subset of $b$.
By Lemma~\ref{lemma:lessthan2}, $\kappa < \mu$.
 That completes Case~2.  That completes the left to right direction.
\smallskip

{\em Right to left}.   Suppose $\kappa < \mu$.  Then by definition of $<$, we have
$\kappa \le \mu$.  On the other hand, if $\kappa = \mu$ then $\kappa \le \mu$
by Lemma~\ref{lemma:le_reflexive}.
\end{proof}  

\begin{lemma}\label{lemma:finitetrichotomy2}
For $\kappa, \mu \in \FregeN$ we have
\begin{eqnarray*}
 \kappa \le \mu  \ \land \ \mu \le \kappa \imp \kappa = \mu.  
\end{eqnarray*}
\end{lemma}

\noindent\begin{proof}  By Lemma~\ref{lemma:letolessthan}, it suffices 
to prove
\begin{eqnarray}
 (\kappa < \mu \ \lor \ \kappa = \mu) \ \land \ (\mu < \kappa \lor \mu = \kappa)
\imp \kappa = \mu.  \label{eq:1953}
\end{eqnarray}
By Theorem~\ref{theorem:finitetrichotomy},
$$ \kappa < \mu \ \lor \ \kappa = \mu \ \lor \ \mu < \kappa$$
and exactly one of the three disjuncts holds. 
Now (\ref{eq:1953}) follows by propositional logic. 
\end{proof}
\smallskip

We next prove two variations on trichotomy that are frequently useful.

\begin{lemma} \label{lemma:le_transitive2}
Suppose $\kappa < \mu \le \lambda$, where $\kappa, \mu, \lambda \in \F$. 
Then $\kappa < \lambda$.
\end{lemma}

\noindent\begin{proof}  By Lemma~\ref{lemma:le_transitive}, we have 
$\kappa \le \lambda$.  We must show $\kappa \neq \lambda$.  Suppose 
$\kappa = \lambda$.  Since $\kappa < \mu$ we have $\lambda < \mu$.
Hence $\lambda \le \mu$.  By hypothesis $\mu \le \lambda$.  By 
Lemma~\ref{lemma:finitetrichotomy2}, $\mu = \lambda$, contradicting
$\mu < \lambda$.  \end{proof}  

\begin{lemma} \label{lemma:finitetrichotomy3}
Let $\kappa,\mu \in \F$.  Then 
$$ \kappa < \mu \lor \mu \le \kappa.$$
\end{lemma}

\noindent\begin{proof} 
\begin{eqnarray}
\kappa < \mu \ \lor \ \kappa = \mu \ \lor \ \mu < \kappa &&\mbox{\qquad by Theorem~\ref{theorem:finitetrichotomy}}
\end{eqnarray}

{\em Case~1}, $\kappa < \mu$.  Then we are done.
\smallskip

{\em Case~2}, $\kappa = \mu$. Then $\kappa \le \mu$ by Lemma~\ref{lemma:le_reflexive}.
\smallskip

{\em Case~3}, $\mu < \kappa$.  Then $\mu \le \kappa$ by the definition of $<$.
\end{proof}

\begin{lemma} \label{lemma:finitedichotomy}
Let $\kappa,\mu \in \F$.  Then 
$$ \kappa \le \mu \lor \mu < \kappa.$$
\end{lemma}

\noindent\begin{proof} 
\begin{eqnarray}
\kappa < \mu \ \lor \ \kappa = \mu \ \lor \ \mu < \kappa &&\mbox{\qquad by Theorem~\ref{theorem:finitetrichotomy}}
\end{eqnarray}

{\em Case~1}, $\kappa < \mu$.  Then $\kappa \le \mu$ by the definition of $<$.
\smallskip

{\em Case~2}, $\kappa = \mu$. Then $\kappa \le \mu$ by Lemma~\ref{lemma:le_reflexive}.
\smallskip

{\em Case~3}, $\mu < \kappa$.  Then we are done.
\end{proof}

\begin{lemma} \label{lemma:le_transitive3}
Let $\kappa, \lambda, \mu \in \FregeN$ and suppose $\kappa \le \lambda < \mu$.
Then  $\kappa < \mu$.
\end{lemma}

\noindent\begin{proof}  By Lemma~\ref{lemma:le_transitive}, we have 
$\kappa \le \mu$.  Since $\kappa < \mu$
is defined as $\kappa \le \mu$ and $\kappa \neq \mu$,  it only remains
to show $\kappa \neq \mu$.  Suppose $\kappa = \mu$.  Then 
$\kappa \le \lambda$ and $\lambda \le \kappa$.  By Theorem~\ref{theorem:finitetrichotomy}, we have $\kappa = \lambda$,
contradiction.  \end{proof}

\begin{lemma} \label{lemma:lessthan_transitive}
Let $\kappa, \lambda, \mu \in \FregeN$ and suppose $\kappa < \lambda < \mu$.
Then  $\kappa < \mu$.
\end{lemma}

\noindent\begin{proof}  Since $\kappa < \lambda$ we have $\kappa \le \lambda$,
by the definition of $<$.  Then by Lemma~\ref{lemma:le_transitive2},
$\kappa < \mu$.  \end{proof}

\begin{lemma}\label{lemma:lessthansuccessor}
Let $\kappa^+ \in \FregeN$. Suppose $\kappa^+$ is inhabited.
  Then $\kappa < \kappa^+$.
\end{lemma}

\noindent\begin{proof}  Since $\kappa^+ \in \FregeN$ and $\kappa^+$
is inhabited,  there exists $x \in \kappa^+$.  Then $x = y \cup \{c\}$
for some $y \in \kappa$ and $c \not \in x$.  
Then $x-y = \{c\}$. 
By Lemma~\ref{lemma:finitecardinals1}, since $x \in \kappa^+$, 
  $x$ is finite.  By Lemma~\ref{lemma:finitedecidable},
$x$ has decidable equality. Therefore 
$x = y \cup \{c\} = y \cup (x-y)$.   Then $y \subseteq x$.
It is a proper subset, since $c \in x$ but $c \not \in y$. 
Now, we will use the right-to-left direction of
 Lemma~\ref{lemma:lessthan2}, 
substituting $\kappa^+$ for
 $\mu$.  That gives us
 $$ \exists x,y\, (x \in \kappa \ \land \ y \in \kappa^+ \ \land \ 
 x \subset y \ \land \ y = x \cup (y-x) \imp \kappa < \kappa^+).$$
 Then take $(y,x)$ for $(x,y)$ in 
 the hypothesis.  That yields
 $$ y \in \kappa \ \land \ x \in \kappa^+ \ \land \ y \subset x \ \land \ 
 x = y \cup (x-y) \imp \kappa < \kappa^+.$$
 Since we have verified all four hypotheses, we may conclude
 $\kappa < \kappa^+$. 
 \end{proof}
 
 \begin{lemma} \label{lemma:successorincreasing}
 For all $m \in \F$,  we do not have $m^+ \le m$.
 \end{lemma}
 
 \noindent\begin{proof}
Suppose $m\in \F$ and $m^+ \le m$.
 By the definition of $\le$, $m^+$ is inhabited.  Then
by Lemma~\ref{lemma:lessthansuccessor},
we have $ m < m ^+ $,  which contradicts Theorem~\ref{theorem:finitetrichotomy}, since $m^+ \le  m$.
\end{proof} 

\begin{lemma} \label{lemma:xnotlessthanx} For $x \in \F$, $x \not<x$.
\end{lemma}

\noindent\begin{proof} Immediate from Theorem~\ref{theorem:finitetrichotomy},
since $x = x$. \end{proof}

\begin{lemma} \label{lemma:xnotlessthanzero}   For $x \in \F$ we have $x \not < \zero$.
\end{lemma}

\noindent\begin{proof}  Suppose $x < \zero$.  We will derive a contradiction.  
\begin{eqnarray*}
x \le \zero &&\mbox{\qquad by definition of $<$}\\
a \in x \land a \subset b \land b \in \zero &&\mbox{\qquad for some $a,b$, 
 by definition of $\le$ } \\
b \in \{\emptyset\}   &&\mbox{\qquad since $\zero = \{\emptyset\}$}\\
b = \emptyset  && \mbox{\qquad by Lemma~\ref{lemma:singleton1}}\\
a = \emptyset  && \mbox{\qquad since $a \subset b$}\\
\emptyset \in x \cap \zero  && \mbox{\qquad by definition of intersection}\\
x = \zero && \mbox{\qquad by Lemma~\ref{lemma:cardinalsdisjoint}} \\
\zero < \zero && \mbox{\qquad since $x < \zero$}\\
\neg\, \zero < \zero && \mbox{\qquad by Lemma~\ref{lemma:xnotlessthanx}}
\end{eqnarray*}
That is the desired contradiction.
\end{proof}

\begin{lemma} \label{lemma:noinsertions}
 For $\kappa, \mu \in \FregeN$, if $\kappa < \mu$,  then $\kappa^+ \le \mu$.
 \end{lemma}
 
 \noindent\begin{proof}  Suppose $\kappa < \mu$. Then there exists $a \in \kappa$
 and $b \in \mu$ such that $b = a \cup (b-a)$. 
 Then 
 \begin{eqnarray*}
 b \in \FINITE \ \land \ a \in \FINITE  &&\mbox{\qquad by Lemma~\ref{lemma:finitecardinals1}}\\
 b-a \in \FINITE &&\mbox{\qquad by Lemma~\ref{lemma:finitedif}}\\
 b-a = \emptyset \ \lor \exists u\,(u \in b-a) &&\mbox{\qquad by
Lemma~\ref{lemma:empty_or_inhabited}}
\end{eqnarray*}
We argue by cases.   
\smallskip

Case 1,  $b-a = \emptyset$. Then 
 $b=a$,  so $a \in \kappa \cap \mu$,
 so by Lemma~\ref{lemma:cardinalsdisjoint}, $\kappa = \mu$,
 contradicting $\kappa < \mu$.
\smallskip

Case 2, $\exists c\, (c \in b-a)$.  Fix $c$.   Then
\begin{eqnarray*}
a \cup \{c\} \in \kappa^+ && \mbox{\qquad by  the definition of successor}\\
a \cup \{c\} \subseteq b &&\mbox{\qquad since $c \in b$}\\
b = (a \cup \{c\}) \cup ( b - (a \cup \{c\})) && \mbox{\qquad by
 Lemma~\ref{lemma:FregeNdecidable} }\\
\kappa^+ \le \mu &&\mbox{\qquad by the definition of $\le$.}
\end{eqnarray*} 
 \end{proof} 
 
\begin{lemma}\label{lemma:successorbounded} 
If $a < b$ and $a,b \in \F$, then $a^+ \in \F$.
\end{lemma} 

\noindent\begin{proof}  Suppose $a < b$ and $a,b \in \F$.  By the definition 
of $<$, we have $a \le b$ and $a \neq b$.  By the definition of $\le$, there exists $v \in b$
and $u \in a$ with $u \in \SSC(v)$.  Then 
\begin{eqnarray*}
v \in \FINITE &&\mbox{\qquad by Lemma~\ref{lemma:finitecardinals1}}\\
u \in \FINITE && \mbox{\qquad by Lemma~\ref{lemma:separablefinite}}\\
v-u \in \FINITE &&\mbox{\qquad by Lemma~\ref{lemma:finitedif}}\\
v-u \neq \emptyset &&\mbox{\qquad by Lemma~\ref{lemma:cardinalsdisjoint}, since $a \neq b$}\\ 
\exists c\, (c \in v-u) &&\mbox{\qquad by Lemma~\ref{lemma:empty_or_inhabited}}\\
 c \in v-u  &&\mbox{\qquad fixing $c$}\\
 u \cup \{c\} \in a^+  &&\mbox{\qquad by definition of successor}\\
 a^+ \in \F && \mbox{\qquad  by Lemma~\ref{lemma:successorF}}
\end{eqnarray*}
\end{proof}

\begin{lemma}\label{lemma:lessthansuccessor2} For $\kappa, \mu \in \F$, we have 
$$ \kappa \le \mu^+ \imp \kappa \le \mu \ \lor \ \kappa = \mu^+.$$
If we also assume $\mu^+ \in \F$ then we have 
$$ \kappa \le \mu^+ \iff \kappa \le \mu \ \lor \ \kappa = \mu^+.$$
\end{lemma}

\noindent{\em Remark}. We cannot replace the $\imp$ with $\iff$ without 
the extra assumption, 
because if $\kappa \le \mu$ there is no guarantee that $\mu^+ \in \F$.
\medskip

\noindent\begin{proof}   Suppose $\kappa \le \mu^+$. 
Then by Lemma~\ref{lemma:letolessthan}, 
\begin{eqnarray*}
\kappa < \mu^+ \lor \kappa = \mu^+.
\end{eqnarray*}
If $\kappa = \mu^+$ we are done; so we may suppose $\kappa < \mu^+$.
Then
\begin{eqnarray*}
\kappa^+ \le \mu^+ &&\mbox{\qquad by Lemma~\ref{lemma:noinsertions}} \\
\exists u\,(u \in \mu^+) &&\mbox{\qquad by the definition of $\le$}\\
\exists u\,(u \in \kappa^+) &&\mbox{\qquad by the definition of $\le$}\\
\kappa \le \mu    && \mbox{\qquad by Lemma~\ref{lemma:ordersuccessor}} 
\end{eqnarray*}
\end{proof}

\begin{lemma}\label{lemma:lessthansuccessor3} For $\kappa, \mu \in \F$, we have 
$$ \kappa < \mu^+ \imp \kappa < \mu \ \lor \ \kappa = \mu.$$
If we also assume $\mu^+ \in \F$ then we have 
$$ \kappa < \mu^+ \iff \kappa < \mu \ \lor \ \kappa = \mu.$$
\end{lemma}

\noindent\begin{proof}  Left to right.
Suppose $\kappa < \mu^+$.  Then by the definition of $<$,
$\kappa \le \mu^+$ and $\kappa \neq \mu^+$.  By Lemma~\ref{lemma:lessthansuccessor2},
$\kappa \le \mu$.  By Lemma~\ref{lemma:letolessthan}, $\kappa < \mu \ \lor \ \kappa = \mu$
as desired.
\smallskip

Right to left.   Assume $\mu^+ \in \F$.  Then $\mu^+$ is inhabited,
by Corollary~\ref{lemma:cardinalsinhabited}.  If $\kappa = \mu $ then $\kappa < \mu^+$ by Lemma~\ref{lemma:lessthansuccessor}.  If $\kappa < \mu$ then $\kappa < \mu^+$ by
Lemma~\ref{lemma:lessthan_transitive}.  
\end{proof}

\begin{lemma} \label{lemma:nothinglessthanzero} 
 $\forall m \in \F\, (\neg\,(m < \zero))$.
\end{lemma}

\noindent\begin{proof}  By definition, $\zero = \{ \emptyset\}$.  Suppose
$m \in \F$ and $m < \zero$.  By definition of $<$, $m \le \zero$ and $m \neq \zero$.
By definition of $\le$, there exist $a$ and $b$ with $a \in m$ and $b \in \zero$ 
and $a \in \SSC(b)$.  Since $\zero = \{\emptyset\}$ we have $b = \emptyset$.  The only 
separable subset of $\emptyset$ is $\emptyset$, so $a = \emptyset$.  Then by 
Lemma~\ref{lemma:cardinalsdisjoint}, $m = \zero$.  But that contradicts
 $m \neq \zero$.  Therefore the assumptions $m \in \F$ and $m < \zero$ are 
 untenable. \end{proof} 
 
  \begin{lemma}\label{lemma:finitemaximal}
 Every nonempty finite subset of $\F$ has a maximal element.
 \end{lemma}
 
 \noindent{\em Remark}.  By Lemma~\ref{lemma:empty_or_inhabited},
 it does not matter whether use ``nonempty'' or ``inhabited''
 to state this lemma.
 \medskip 
 
\noindent\begin{proof}  The formula to be proved 
is 
\begin{eqnarray*}
&&\forall x \in \FINITE\, (x \subseteq \F \imp x \neq \emptyset \imp \exists m \in x \,  \forall t\,(t \in x \imp t \le m)) 
\end{eqnarray*}
 The formula is stratified, giving $m$ and $t$ index 0 and $x$ index 1.  $\F$ and $\FINITE$
 are parameters, and do not require an index.   Therefore we may proceed by 
 induction on finite sets. 
 \smallskip
 
{\em Base case}\,: immediate, since $\emptyset \neq \emptyset$. 
\smallskip

{\em Induction step}.  Let $x$ be a finite subset of $\F$ and $c  \in \F-x$. 
By Lemma~\ref{lemma:empty_or_inhabited}, $x$ is empty or inhabited.  If 
$x = \emptyset$, then $c$ is the maximal element of $x \cup \{c\}$, and we
are done.  So we may assume $x$ is inhabited.  Then by the induction hypothesis,
$x$ has a maximal element $m$.  By Theorem~\ref{theorem:finitetrichotomy},
$c \le m$ or $m < c$.  If $c \le m$, then $m$ is the maximal element of 
$x \cup \{c\}$.  If $m < c$, then $c$ is the maximal element of $x \cup \{c\}$,
by the transitivity of $\le$.  \end{proof} 

\begin{lemma} \label{lemma:xnotequalsuccessorx}
For $x \in \F$ and $x^+ \in \F$,  we have $x \neq x^+$.
\end{lemma}

\noindent\begin{proof}   Suppose $x = x^+$; then 
\begin{eqnarray*}
z \in x^+ &&\mbox{\qquad for some $z$, by Corollary~\ref{lemma:cardinalsinhabited}}\\
z = u \cup \{c\} && \mbox{\qquad for some $u \in c$ and $c \not\in u$, by definition of successor} \\
u \cup \{c\} \in x^+ &&\mbox{\qquad by the previous two lines}\\
u \cup \{c\} \in x  &&\mbox{\qquad since $x = x^+$}\\
u \cup \{c\} \in \FINITE && \mbox{\qquad by Lemma~\ref{lemma:finitecardinals1}}\\
u \sim u \cup \{c\}  &&\mbox{\qquad by Lemma~\ref{lemma:finitecardinals2}} \\
u \cup \{c\} \neq u  && \mbox{\qquad since $c \not \in u$} 
\end{eqnarray*}
Now $u \cup \{c\}$ is a finite set, similar to a proper subset of itself 
(namely $u$).  Then by definition, $u \cup \{c\}$ is infinite.
By
Theorem~\ref{theorem:infiniteimpliesnotfinite}, it is not finite. But it 
is finite.   That contradiction shows $x \neq x^+$. \end{proof}

\begin{lemma} \label{lemma:xlessthansuccessorx}
For $x \in \F$ and $x^+ \in \F$,  we have $x < x^+$.
\end{lemma}
 
\noindent\begin{proof} Let $u \in x$ and $u \cup \{c \} \in x^+$, 
with $c \not \in x$.  Then by definition of $\le$ we have $x \le x^+$.
By Lemma~\ref{lemma:xnotequalsuccessorx}, we have $x \neq x^+$.
Then by definition of $<$, we have $x < x^+$. 
\end{proof}
  
\section{Power sets and similarity} 

We will replace Rosser and Specker's use of the full power set $\SC$
by the separable power set $\SSC$.   In this section we prove some 
lemmas from Specker \sect2, and some other similar lemmas.
 For finite sets $a$, since finite
sets have decidable equality, every unit subclass is separable, which
is helpful.  We begin with 
  Specker's Lemma~2.6, which we take in two steps with the 
next two lemmas, and after that Specker 2.4 and 2.3. 

\begin{lemma} \label{lemma:subset_usc} Let $y \in \SSC( \USC(a))$.
Then there exists $z \in \SSC(a)$ such that $y = \USC(z)$.
\end{lemma}

\noindent\begin{proof}  Suppose $y \in \SSC( \USC(a))$. Define 
 \begin{eqnarray}
 z := \{ u : \{u\} \in y \}. \label{eq:1870}
 \end{eqnarray}
 That definition is legal since the 
formula is stratified giving $u$ index 0 and $y$ index 2. 
Then $y = \USC(z)$ since the members of $y$ are the singletons
of the members of $z$.   I say that $z \subseteq a$:
Suppose 
$u \in z$. Then 
\begin{eqnarray*}
\{u\} \in y && \mbox{\qquad by (\ref{eq:1870})} \\
\{u\} \in \USC(a) && \mbox{\qquad since $y \subseteq \USC(a)$}\\
u \in a          && \mbox{\qquad by definition of $\USC(a)$} 
\end{eqnarray*}
Therefore $z \subseteq a$, as claimed.   It remains to show 
that $z$ is a separable subset of $a$; it suffices to show that 
for $u \in a$, we have $u \in z \ \lor \ u \not\in z$.  
Suppose $u \in a$.  Then by (\ref{eq:1870}), 
\begin{eqnarray*}
 && u \in z \ \lor u \not \in z \\
&\iff& \{u\} \in y \ \lor \{u\} \not\in y 
\end{eqnarray*}
and that is true since $y$ is a separable subset of $\USC(a)$.
\end{proof}

\begin{lemma}[Specker 2.6]\label{lemma:sscusc}
$\Nc{{\SSC(\USC(a))}} = \Nc{{\USC(\SSC(a))}}$.
\end{lemma}

\noindent{\em Remarks}.  Or course Specker has $\SC$ instead of $\SSC$.
We follow the proof from~\cite{rosser1953}, p.~368,
that Specker cites, checking it constructively with $\SSC$ in place of $\SC$.
But fundamentally, this lemma is just about shuffling brackets.  We have
$ \{\{p\},\{q\},\{r\}\} \in \SSC(\USC(a))$ corresponding to $\{ \{ p,q,r\}\} \in \USC(\SSC(a))$.
It is a useful result but not a deep one. 
\medskip

\noindent\begin{proof}  Let 
$$ W:= \{ u: \exists z\,( u = \langle \{z\}, \USC(z)\rangle)\}.$$
The definition is stratified giving  
  $z$   index 1, so $\{z\}$ and $\USC(z)$ both get index 2, and $u$ gets
index 4.  It follows that $W$ is a relation (contains only ordered pairs) and 
\begin{eqnarray}
 \langle x,y \rangle \in  W  \iff
  \exists z\,(x = \{z\} \ \land \ y = \USC(z)). \label{eq:1881}
\end{eqnarray}

I say that $W$ is (the graph of) a one-one-function mapping
$\USC(\SSC(a))$ onto $\SSC(\USC(a))$. (Formally there is no distinction
between a function and its graph.)  For if $x$ is given, then $z$ is 
uniquely determined, so $y$ is uniquely determined;  and if $y$ is 
given with $y = \USC(z)$, then $z = \bigcup y$ is unique, so $x = \{z\}$
is unique.  Hence $W$ is a function and one-to-one.  It remains to show 
that $W$ is onto.  Let $y \in \SSC(\USC(a))$. 
By Lemma~\ref{lemma:subset_usc}, there exists $z \in \SSC(a)$ such that
  $y = \USC(z)$.  Then 
$\langle \{z\},y \rangle \in W$.  Hence $y$ is in the range of $W$. 
Since $y$ was an arbitrary member of $\SSC(\USC(a))$, it follows that 
$W$ is onto. 
\smallskip

We have shown that $W$ is a similarity from 
$\SSC(\USC(a))$ to $\USC(\SSC(a))$.  Therefore those two sets 
have the same cardinal.   \end{proof}

\begin{lemma} \label{lemma:singletons_similar} 
 Any two unit classes
are similar.
\end{lemma}

\noindent\begin{proof}
Let $\{a\}$ and $\{b\}$ be unit classes. 
Define $f = \{ \langle a, b\rangle\}$.  One can verify that $f: \{a\} \to \{b\}$
is a similarity.  We omit the 75 inferences required to do so.  \end{proof}
 
\begin{lemma} \label{lemma:similar_to_singleton}  
Any set similar to a unit class is a unit class. 
\end{lemma}

\noindent\begin{proof} Let $x \sim \{a\}$.
Then let $f:x \to \{a\}$ be a similarity.  Since $f$ is onto, there exists $c \in x$
with $f(c) = a$.  Let $e \in x$.  Then $f(e) \in \{a\}$, so $f(e) = a$.  Since 
$f$ is one-to-one, $e = c$.  Then $x = \{c\}$.  
\end{proof}

\begin{lemma} \label{lemma:one_members} We have 
 $$ u \in \one \iff \exists a\, (u = \{a\}).$$
\end{lemma}

\noindent\begin{proof} By definition, $\one = \zero^+$ and $\zero = \{\emptyset\}$.
For any $a$, we have $a \not \in \emptyset$, so 
$$\emptyset \cup \{a\} = \{a\} \in \zero^+ = \one.$$
Conversely, if $u \in \one$, then $u= \emptyset \cup \{a\}$ for some $a$,
by definition of successor, so $u = \{a\}$.  
\end{proof}

\begin{lemma} \label{lemma:usc_subset3}
Suppose $a$ and $b$ are finite sets.   Then
$$  a \in \SSC(b) \imp \USC(a) \in \SSC (\SSC (b)) .$$
\end{lemma}

\noindent\begin{proof}  Suppose $a \in \SSC(b)$. 
Since $b$ is finite, it has decidable equality,
by Lemma~\ref{lemma:finitedecidable}. Therefore
$\USC(b) \subseteq \SSC(b)$.  Since $\USC(a) \subseteq \USC(b)$,
we have  
\begin{eqnarray}
\USC(a) \subseteq \SSC(b) \label{eq:2570}
\end{eqnarray}
It remains to show that $\USC(a)$ is a separable
subset of $\SSC(b)$; that is, 
$$\SSC(b) = \USC(a) \cup (\SSC(b)-\USC(a)).$$
By extensionality and the definitions of subset and union,
it suffices to show 
\begin{eqnarray}
t \in \SSC(b) & \iff & t \in \USC(a) \ \lor \ 
                     ( t \in \SSC(b) \ \land \ t \not\in \USC(a))
            \label{eq:2580}
\end{eqnarray}
{\em Right to left}.   It suffices to show $t \in \USC(a) \imp t \in \SSC(b)$.
Let $t \in \USC(a)$.  Then $t = \{c\}$ for some $c \in a$.  Since $b$ has
decidable equality, $t$ is a separable subset of $b$.  That completes the 
right-to-left direction.
\smallskip

Left to right:
suppose $t \in \SSC(b)$.  Then $t \in \FINITE$, by Lemma~\ref{lemma:separablefinite}.
Then $\Nc{t} \in \F$, by Lemma~\ref{lemma:finitecardinals3}. 
Then by Lemma~\ref{lemma:FregeNdecidable}, 
$$\Nc{t} = \one \ \lor \ \Nc{t}\neq \one. $$  
Case 1, $\Nc{t} = \one$. By Lemma~\ref{lemma:one_members}, $t$
is a unit class.   Since $a \in \SSC(b)$,  we have 
$$x \in b \imp x \in a \ \lor \ x \not \in a.$$
Since $t \in \USC(a)$ if and only if for some $x$ we have $t=\{x\} \ \land \ x \in a$,
we have 
$$t \in \SSC(b) \imp  t \in \USC(a) \ \lor \ t \not\in \USC(a).$$
That completes Case~1.
\smallskip

Case~2, $\Nc{t} \neq \one$. Then $\Nc{t}$ is not a unit class, 
by Lemma~\ref{lemma:one_members} and Lemma~\ref{lemma:cardinalsdisjoint},
so the second disjunct on the right holds.
\end{proof}

\begin{lemma}[Specker 2.4] \label{lemma:uscsimilar} For any sets $a$ and $b$
$$a \sim b   \iff  \USC(a) \sim \USC(b).$$
\end{lemma} 

\noindent\begin{proof}  Left-to-right.  Suppose $f: a \to b$ is 
a similarity.  Let $g$ be the singleton image of $f$, namely
$$ g := \{ \langle \{ u\}, \{v\} \rangle : \langle u,v \rangle \in f.$$
The definition is legal since the formula is stratified, giving $u$ and $v$
the same index.  Then $g: \USC(a) \to \USC(b)$ is a similarity.  We omit
the straightforward proof. 
\smallskip

Right-to-left.   Let $g: \USC(a) \to \USC(b)$ be a similarity.
Define
$$ f : = \{ \langle u,v \rangle : \langle \{u\},\{v\} \rangle \in g \}.$$
Again the definition is legal since the formula is stratified, giving $u$ 
and $v$ the same index.  Then $f: a \to b$ is a similarity.  We omit 
the proof.\end{proof}
\smallskip

\begin{lemma}[Specker 2.3] \label{lemma:sscsimilar} For any sets $a$ and $b$
 $$ a \sim b \imp \SSC(a) \sim \SSC(b).$$
\end{lemma}

\noindent\begin{proof}
  Let $f: a \to b$ be a similarity.  Define 
$$g := \{  \langle u, f``u \rangle :  u \in \SSC(a) \}$$
where $f``u$ is the image of $u$ under $f$, i.e., the range of the 
restriction of $f$ to $u$.  Then $g:\SSC(a) \to \SSC(b)$. The fact that 
the values of $g$ are separable subsets of $b$ follows from 
Lemma~\ref{lemma:separable_similarity}. We omit the proof that $g$
is one-to-one. To prove $g$ is onto, let $y \in \SSC(b)$. Then 
define
$$x = \{ u \in a: \exists v\,(v \in y \ \land \ \langle u,v \rangle \in f )\}.$$
The   formula is stratified, giving $u$ and $v$ index 0 and $x$ and $y$ index 1.
Hence $x$ can be defined.  We omit the proof that $g(x) = y$.
($x$ can also be defined using the operations of domain and inverse relation,
which in turn can be defined by stratified comprehension.) 
\end{proof}

\begin{lemma} \label{lemma:usc_subset_ssc} If $a$ has decidable equality,
then $\USC(a) \subseteq \SSC(a)$.
\end{lemma}

\noindent\begin{proof} Let $x \in \USC(a)$.  Then $x = \{u\}$ for some $u \in a$.
Then $x \subseteq a$.  We must show $a = x \cup (a-x)$.  By extensionality,
that follows from $$\forall u\, (u \in a \iff u \in x \ \lor \ u \in a - x),$$
which in turn follows from decidable equality on $a$. 
\end{proof}

\begin{lemma} \label{lemma:usc_subset}
 For all $a,b$, 
$$ a \subseteq b \iff \USC(a) \subseteq \USC(b),$$
\end{lemma}

\noindent\begin{proof}  Left to right.  Suppose $a \subseteq b$
and $t \in \USC(a)$.  We must show $t \in \USC(b)$.  
Then $t = \{ x\}$ for some $x\in a$.  Since $a \subseteq b$ we 
have $x \in b$.  Then $t \in \USC(b)$. That completes the left-to-right
direction.
\smallskip

Right to left.  Suppose $\USC(a) \subseteq \USC(b)$ and $t \in a$.  We
must prove $t\in b$.  Since $t\in a$ we have $\{t\} \in \USC(a)$.  
Then $\{t\} \in \USC(b)$.  Then $\{t\} = \{q\}$ for some $q\in b$.
Then $t=q$. Then $t \in b$ as desired.  
\end{proof}

\begin{lemma} \label{lemma:ssc_subset1}For all $a,b$, 
$$ a \in \SSC (b) \iff UCS (a) \in \SSC (\USC(b)).$$
\end{lemma}
\noindent\begin{proof} Left to right.  Suppose $a \in \SSC(b)$. 
Then $a \subset b$ and 
\begin{eqnarray}
b &=& a \cup (b-a). \label{eq:2041}
\end{eqnarray}
  By Lemma~\ref{lemma:usc_subset},
\begin{eqnarray}
\USC(a) \subseteq \USC(b) \label{eq:2083}
\end{eqnarray}
It remains to show that $UCS (a)$ is a stable subset of $\USC(b)$; 
that is, 
\begin{eqnarray}
\USC(b) &=&  \USC(a) \cup (\USC(b)-\USC(a)).  \label{eq:2046}
\end{eqnarray} 
By extensionality and the definitions of union and set difference,
that is equivalent to 
\begin{eqnarray}
 t \in \USC(b) &\iff& t \in \USC(a) \ \lor \ (t \in \USC(b) \ \land \ t \not \in \USC(a)). \label{eq:2051}
 \end{eqnarray}
Then we need only consider unit classes $ t = \{ x\}$, and using the fact 
that $\{x\} \in \USC(b) \iff t \in b$, and $\{x \} \in \USC(a) \iff t \in a$,
(\ref{eq:2051}) follows from (\ref{eq:2083}).   
\end{proof} 

\begin{lemma} \label{lemma:ssc_subset2}
For all $a,b$, we have
$$  a \in \SSC(b) \iff \SSC(a) \subseteq \SSC(b).$$
\end{lemma}

\noindent\begin{proof} {\em Left to right}. Suppose
$a \in \SSC(b)$.  Then $a \subseteq b$ and 
\begin{eqnarray}
b &=& a \cup (b-a) \label{eq:2109}
\end{eqnarray}
Now let $x \in \SSC(a)$.  We must show $x \in \SSC(b)$.
Since $x \in \SSC(a)$, we have $x \subseteq a$. 
Since $a \subseteq b$ we have $x \subseteq b$.
We have
\begin{eqnarray*}
x \in \SSC(a) && \\
a = x \cup (a - x) && \mbox{\qquad by definition of $\SSC(a)$}\\
b = (x \cup (a-x)) \cup (b - (x \cup (a-x))) &&
                     \mbox{ \qquad by (\ref{eq:2109}) }\\
b = x \cup (b-x) &&\\
x \in \SSC(b) && \mbox{ by definition of $\SSC(b)$} 
\end{eqnarray*}
That completes the left-to-right direction.
\smallskip

{\em Right to left}.   Suppose $\SSC(a) \subseteq \SSC(b)$.
We have to show $a \in \SSC(b)$; but that follows 
from $a \in \SSC(a)$ and the definition of subset.
That completes the right to left direction.
\end{proof}

\begin{lemma} \label{lemma:ssc_subset4} 
Let $b$ be a finite set.  Then the subset relation on $\SSC(b)$
is decidable.  That is,
$$ \forall x,y \in \SSC(b) \, ( x \subseteq y \ \lor\ x \not \subseteq y).$$
\end{lemma}

\noindent\begin{proof} Assume $b\in \FINITE$.
By Lemma~\ref{lemma:finitepowerset}, $\SSC(b) \in \FINITE$.
Then by Lemma~\ref{lemma:finitedecidable}, 
\begin{eqnarray}
\SSC(b) \in \DECIDABLE \label{eq:2132}
\end{eqnarray} 

 We will prove by induction on finite sets y that 
\begin{eqnarray}
y \in \SSC(b) \imp \forall x \in \SSC(b)\, (x \subseteq y \ \lor \ x \not \subseteq y). \label{eq:2150}
\end{eqnarray}
It is legal to proceed by induction, since the formula is stratified.
\smallskip

{\em Base case}.  When $y = \emptyset$, we will prove 
$$ \forall x \in \SSC(b)\, (x \subseteq \emptyset \lor x \not \subseteq \emptyset).$$
Assume $x \in \SSC(b)$. 
We have $x \subseteq \emptyset$ if and only if $x = \emptyset$,  so it suffices
to prove $x = \emptyset \ \lor \ x \neq \emptyset$.  But that
follows from (\ref{eq:2132}).   That completes the base case.

{\em Induction step}.  Let $y = z \cup \{c\}$, with $c \not\in z$ and $z \in \SSC(b)$
and $y \subseteq b$.
Then $c \in b$. 
  The induction hypothesis is
\begin{eqnarray}
z \in \SSC(b) \imp \forall x \in \SSC(b)\, (x \subseteq z \ \lor\ x \not \subseteq z). \label{eq:2151}
\end{eqnarray}
We have to prove 
\begin{eqnarray}
y \in \SSC(b) \imp\forall x \in \SSC(b)\, (x \subseteq y \ \lor \ x \not \subseteq y) \label{eq:2158} 
\end{eqnarray} 
Assume $y \in \SSC(b) $ and $x \in \SSC(b)$.  We have to prove 
$ x \subseteq y \ \lor \ x \not \subseteq y $.
That is,
\begin{eqnarray*} \\
&&  x \subseteq z \cup \{c\}  \ \lor x \not \subseteq z \cup \{c\} 
\end{eqnarray*}
We have
\begin{eqnarray*}
y \in \SSC(b) &&  \mbox{\qquad assumed above}\\
z \cup \{c\} \in \SSC(b) && \mbox{\qquad since $y = z \cup \{c\}$}
\end{eqnarray*}
I say that $z \in \SSC(b)$. To prove that, let $u\in z$.  Since
$z \cup \{c\} \in \SSC(b)$, $u \in z \cup \{c\} \ \lor \ u \not\in z \cup\{c\}$.
Since $c \not\in z$, $u \neq c$.  Therefore $u \in z \lor u \not\in z$.
Then $z \in \SSC(b)$ as claimed.  
\smallskip

I say that also  $x-\{c\} \in \SSC(b)$. Since $b$ is finite, it has decidable
equality by Lemma~\ref{lemma:finitedecidable}.  Then for $y \in b$,
we have $y = c \ \lor \ y \neq c$.  Since $x \in \SSC(b)$ we have 
$y \in x \ \lor \ y \not \in x$.  Then a short argument by cases shows
$y \in x - \{c\} \ \lor \ y \not \in x \{c\}$.  Then $x-\{c\} \in \SSC(b)$,
as claimed. 
\smallskip

By (\ref {eq:2151}) and $z \in \SSC(b)$,  we have
\begin{eqnarray}
\forall x \in \SSC(b)\, (x \subseteq z \ \lor\ x \not \subseteq z). \label{eq:2152}
\end{eqnarray}
Since $x \in \SSC(b)$, we have $c \in x \ \lor \ c \not \in x$.  We 
argue by cases accordingly.
\smallskip

Case~1: $c \in x$.  Then $x \subseteq z \cup \{c\}$ if and only if 
$x-\{c\} \subseteq z$. By (\ref{eq:2152}),
instantiated to $x-\{c\}$ in place of $x$ (which is allowed since
$x-\{c\} \in \SSC(b)$),  we have 
$$ x-\{c\} \subseteq z \ \lor \  x-\{c\} \not\subseteq z.$$
That completes Case~1.
\smallskip

Case~2: $c \not\in x$.  Then $x \subseteq z \cup \{c\} \iff x \subseteq z$,
so (\ref{eq:2158}) follows from  the induction hypothesis (\ref{eq:2151}).
That completes Case~2.   That completes the induction step. 
\end{proof}

\begin{lemma} \label{lemma:ssc_subset3}
Suppose $a$ and $b$ are finite sets.   Then
$$  a \in \SSC(b) \imp \SSC(a) \in \SSC (\SSC (b)) .$$
\end{lemma}

\noindent\begin{proof}  Suppose $a \in \SSC(b)$.
By Lemma~\ref{lemma:ssc_subset2}, 
\begin{eqnarray}
\SSC(a) \subseteq \SSC(b) \label{eq:2101}
\end{eqnarray}
It remains to show that $\SSC(a)$ is a separable
subset of $\SSC(b)$; that is, 
$$\SSC(b) = \SSC(a) \cup (\SSC(b)-\SSC(a)).$$
By extensionality and the definitions of subset and union,
it suffices to show 
\begin{eqnarray}
t \in \SSC(b) & \iff & t \in \SSC(a) \ \lor \ 
                     ( t \in \SSC(b) \ \land \ t \not\in \SSC(a))
            \label{eq:2111}
\end{eqnarray}
The right-to-left direction follows logically from 
(\ref{eq:2101}) and the definition of subset.
\smallskip

Ad the left-to-right direction of (\ref{eq:2111}):
suppose $t \in \SSC(b)$.  Then $t \subseteq b$.  
By Lemma~\ref{lemma:ssc_subset4},
\begin{eqnarray}
t \subseteq a \lor t \not \subseteq a.  \label{eq:2143}
\end{eqnarray} 
We argue by cases using (\ref{eq:2143}).
\smallskip

Case 1: $t \subseteq a$.  It suffices to prove $t \in \SSC(a)$.  It remains
to prove $a = t \cup (a-t)$.  We have
\begin{eqnarray*}
\forall u \in b\, (u \in t \ \lor \ u \not\in t) \mbox{\qquad since $t \in \SSC(b)$}\\
\forall u \in a\, (u \in t \ \lor \ u \not\in t) \mbox{\qquad since $a \subseteq b$}
\end{eqnarray*}
Then $a = t \cup (a-t)$ by the definitions of union and set difference.
That completes Case~1.
\smallskip

Case 2: $t \not \subseteq a$.  Then $t \not\in \SSC(a)$.  Since $t \in \SSC(b)$,
the second disjunct on the right of (\ref{eq:2111}) holds.  That completes
Case~2.   \end{proof}

\begin{lemma} \label{lemma:usc_successor} For all $a$ and $c \not\in a$, 
we have $$ \USC(a \cup \{c\}) = \USC(a) \cup \{ \{c\}\}.$$
\end{lemma}

\noindent\begin{proof} By extensionality it suffices to verify the two
sides have the same members.  
\smallskip

{\em Left to right}\.:  Let $x \in \USC(a \cup \{c\})$.  Then $x = \{ u\}$ for 
some $u \in a \cup \{c\}$.  Then $u \in a \ \lor \ a = c$.  If 
$u\in a$ then $x \in \USC(a)$ and hence $x \in \USC(a) \cup \{\{c\}\}$.
That completes the left-to-right direction.
\smallskip

{\em Right to left}\,: Let $x \in \USC(a) \cup \{\{c\}\}$.  Then $x \in \USC(a) \ \lor \ x=\{c\}$.
If $x \in \USC(a)$, then $x \in \USC(a \cup \{c\})$ by Lemma~\ref{lemma:usc_subset}.
If $x = \{c\}$, then $x \in \USC(a \cup \{c\})$ by definition of $\USC$.
\end{proof} 

\begin{lemma} \label{lemma:usc_dif2} For all $a,b$ we have $$\USC(a-b) = \USC(a) -\USC(b).$$
\end{lemma}

\noindent\begin{proof}  By the definitions of $\USC$ and set difference, using about 50
straightforward inferences, which we choose to omit. \end{proof} 

\begin{lemma} \label{lemma:usc_empty}  $\USC(\emptyset) = \emptyset$.
\end{lemma}

\noindent\begin{proof}  Suppose $x \in \USC(\emptyset)$. By definition 
of $\USC$, there exists $a \in \emptyset$ such that $x = \{a\}$. 
But that contradicts the definition of $\emptyset$.  
\end{proof}

\begin{lemma}\label{lemma:usc_up_down} For every $x$ and $a$,
$$ x \in a \iff \{x\} \in \USC(a).$$
\end{lemma}

\noindent\begin{proof} {\em Left to right}\,: by definition of $\USC(a)$.
{\em Right to left}\,: if $\{x\} \in \USC(a)$, then for some $y \in a$,
$\{x\} = \{y\}$.  Then by extensionality $x = y$. 
\end{proof}

 \begin{lemma} \label{lemma:ssc_empty} $\SSC(\emptyset) = \{\emptyset\}$.
\end{lemma}

\noindent\begin{proof}  The only subset of $\emptyset$ is $\emptyset$, and it is 
a separable subset.  \end{proof}   

\begin{lemma} \label{lemma:similarinhabited} Suppose $a \sim b$ 
and $a$ is inhabited.  Then $b$ is inhabited.
\end{lemma}

\noindent\begin{proof}  Let $f:a \to b$ be a similarity.  Since 
$a$ is inhabited, there exists some $c \in a$.  Fix $c$.  
Then $f(c) \in b$.  Hence $b$ is inhabited.  
\end{proof}

 \begin{lemma}[Bounded DNS] \label{lemma:boundedDNS}  Let $P$ be any set,
and let $y \in \F$. 
Then 
$$ \neg\neg\,\forall x\, (x \in \F \imp x < y \imp x \in P) \iff 
\forall x\, (x \in \F \imp x < y \imp \neg\neg\, x \in P).$$
\end{lemma}

\noindent{\em Remarks}.  This lemma is closely related to 
Lemma~\ref{lemma:finiteDNS}, and can be derived from that lemma,
but here we just prove it directly.   
\medskip

\noindent\begin{proof}  The left-to-right direction is logically valid.
We prove the right-to-left implication by induction on $y$.  The formula to be proved is  stratified, giving $x$ and $y$ index 0, so induction is legal. 
\smallskip

{\em Base case}\,: by Lemma~\ref{lemma:xnotlessthanzero}, $x < 0$ can never hold. 
That completes the base case.
\smallskip

{\em Induction step}\,:  The key fact will be Lemma~\ref{lemma:lessthansuccessor3}:
\begin{eqnarray}
 x < y^+ \iff  x < y \ \lor \ x =y. \label{eq:6263} 
 \end{eqnarray}
Assume $y^+$ is inhabited (as for any proof by induction). 
Then
\begin{eqnarray*}
\forall x\,(x \in \F \imp x < y^+ \imp \neg\neg\,x \in P) &&\mbox{\qquad assumption}\\
\forall x\,(x \in \F \imp ( x < y \ \lor \ x =y) \imp \neg\neg\,x \in P) &&\mbox{\qquad by (\ref{eq:6263}) }\\
\forall x\,(x \in \F \imp   (x < y   \imp \neg\neg\,x \in P)  \ \land\ (x = y \imp \neg\neg\,x \in P)) &&\mbox{\qquad by logic}\\
\forall x\,(x \in \F \imp   (x < y   \imp \neg\neg\,x \in P )) \ \land \ \neg\neg\, y \in P
&&\mbox{\qquad by logic}\\
\neg\neg\,\forall x\,(x \in \F \imp   (x < y   \imp  x \in P ))  \ \land \ \neg\neg\, y\in P &&\mbox{\qquad induction hyp.}\\
  \neg\neg\, \forall x\, (x \in \F \imp x \le y \imp x \in P) &&\mbox{\qquad by (\ref{eq:6263})}
\end{eqnarray*}
That completes the induction step.  \end{proof}

\section{Cardinal exponentiation}

Specker~4.1 follows Rosser in defining $2^m$ for cardinals $m$. They
define $2^m$ to be the cardinal of $\SC(a)$ where $\USC(a) \in m$. 
That definition requires some modification to be of use constructively.
It is {\em separable} subsets of $a$ that correspond to functions
from $a$ to 2, so it makes sense to use $\SSC(a)$, the class of separable
subsets of $a$, instead of $\SC(a)$. 

 \begin{definition} \label{definition:exp}
For finite cardinals $m$, we define
$$ 2^m =   \{ u : \exists a\,( \USC(a)\in m \ \land \ u \sim  \SSC(a)) \}  .$$
\end{definition}

The following lemma shows that our definition is classically equivalent
to Specker's definition.  

\begin{lemma} \label{lemma:expuscssc}
Let $m \in \FregeN$ and $\USC(a) \in m$.  Then 
$\SSC(a) \in 2^m$, and $2^m = \Nc{\SSC(a)}$.
\end{lemma}

\noindent{\em Remark}. This is Specker's definition of $2^m$, but 
our definition avoids a case
distinction as to whether $m$
 does or does not contain a set of the form $\USC(a)$.
 \medskip
 
 \noindent\begin{proof}  Suppose $m \in \FregeN$ and $\USC(a) \in m$.
I say that $2^m$ is a cardinal, i.e., it is closed under similarity.
Suppose $u$ and $v$ are members of $2^m$.  Then there exist $a$ and $b$
such that $\USC(a)$ and $\USC(b)$ are both in $m$ and $u \sim \SSC(a)$
and $v \sim \SSC(b)$.  Then by Lemma~\ref{lemma:finitecardinals2},
$\USC(a) \sim \USC(b)$.  By Lemma~\ref{lemma:uscsimilar}, $a \sim b$.
By Lemma~\ref{lemma:sscsimilar}, $\SSC(a) \sim \SSC(b)$.  By 
Lemma~\ref{lemma:sim}, $u \sim v$.  Hence, as claimed, 
$2^m$ is a cardinal.  
\smallskip

Therefore $2^m$ and $\Nc{\SSC(a)}$ are both closed under similarity.  
Since they both contain $\SSC(a)$,  they each consist of all sets 
similar to $\SSC(a)$.  Hence by extensionality, they are equal.  
\end{proof}
\medskip

\noindent{\em Remark}.  We note that $2^m \neq \emptyset$ 
does not {\em prima facie} imply that $2^m$ is inhabited,  so we 
must carefully distinguish these two statements as hypotheses of lemmas.
$2^m$ is inhabited if $m$ contains a set of the form $\USC(a)$.  
$2^m \neq \emptyset$ means not-not $m$ contains such a set.  

{\em Discussion}.  It is possible, of course, to investigate what happens if we 
use intuitionistic logic, but keep the classical definitions of order and exponentiation.
The most obvious difficulty with this approach is that the integers $\F$ are not closed
under exponentiation.   For example, let us calculate what $2^{\one}$ would be.  We have
$\{ \{\emptyset\}\} = \USC(\{\emptyset\}) \in \one$.  So $2^{\zero}$ would be the cardinal of $\SC(\{\emptyset)\}$,
instead of the cardinal of $\SSC(\{\emptyset\})$.  But $\SC(\{\emptyset\})$ contains every set 
of the form $X_P = \{x: x = \emptyset \ \land \ P\}$,  where $P$ is a stratified formula not 
containing the variable $x$.  Unless we can prove or refute $P$, we cannot prove that $X_P$
is one of the two members of $\SSC(\emptyset)$, and in fact there is no hope of proving $2^{\one}$ 
is an integer.  Hence this notion is useless for constructive mathematics in NF.  Still we 
did investigate the matter further, in the hope that this approach might 
help analyze Specker's proof.  In short, it did not help.  Without the axiom of choice,
one can prove nothing useful about large cardinals.  For example, one cannot prove $2^x = 2^y \imp x =y$
for cardinals; there might even be incomparable $x,y$ such that $2^x = 2^y$.  That might even
be the case with $2^x = 2^\kappa = \kappa$, where $\kappa$ is the cardinal of $\V$.  We consider
this subject no further.
\smallskip

\begin{lemma}\label{lemma:expdefinable} The graph of the exponentiation function
$$ \{ \langle m, 2^m \rangle : m \in \F \}$$
is definable in \INF.
\end{lemma}

\noindent\begin{proof}  We have to show that the relation is definable by 
a formula that can be stratified, giving the two members of ordered pairs 
the same index.  The formula in Definition~\ref{definition:exp} is
$$ 2^m =   \{ u : \exists a\,( \USC(a)\in m \ \land \ u \sim  \SSC(a)) \}  .$$
Stratify it, giving $a$ index 0, $\USC(a)$ and $\SSC(a)$ and $u$ index 1, 
$m$ index 2.   Then $2^m$ gets one index higher than $u$, namely 2, which is 
the same index that $m$ gets.  \end{proof}  

\begin{lemma} \label{lemma:exp_inhabited}  
If $2^m$ is inhabited, then there exists $a$ such that $\USC(a) \in m$
and $\SSC(a) \in 2^m$.
\end{lemma}

\noindent\begin{proof}  Suppose $2^m$ is inhabited. Then by  Definition~\ref{definition:exp},
there exists $a$ with $\USC(a) \in m$, and $2^m$ contains any set similar to $\SSC(a)$.
Since $\SSC(a) \sim \SSC(a)$, by Lemma~\ref{lemma:sim}, we have $\SSC(a) \in 2^m$.
\end{proof}

 \begin{lemma} \label{lemma:finiteexp}
Let $m$ be a finite cardinal.  If $2^m$ is inhabited, then 
$2^m$ is a finite cardinal.
\end{lemma}

\noindent\begin{proof}  Suppose $m$ is a finite cardinal and 
$2^m$ is inhabited.  By Definition~\ref{definition:exp}, there 
exists $a$ such that $\USC(a) \in m$ and $\SSC(a) \in 2^m$.
Then $\Nc{\SSC(a)} = 2^m$, by Definition~\ref{definition:exp}.
We have 
\begin{eqnarray*}
\USC(a) \in \FINITE && \mbox{\qquad by Lemma~\ref{lemma:finitecardinals1}}\\
a \in \FINITE && \mbox{\qquad by Lemma~\ref{lemma:uscfinite}} \\ 
\SSC(a) \in \FINITE &&\mbox{\qquad by Lemma~\ref{lemma:finitepowerset}}\\
\Nc{\SSC(a)} \in \FregeN &&\mbox{\qquad by Lemma~\ref{lemma:finitecardinals3}}\\
2^m \in \FregeN && \mbox{\qquad since $\Nc{\SSC(a)} = 2^m$}
\end{eqnarray*}
\end{proof}

\begin{lemma} \label{lemma:exp_zero}
$2^{\zero} = \one$.
\end{lemma}

\noindent\begin{proof}  $\zero = \{\emptyset\}$.
 It therefore contains $ \emptyset  = \USC(\emptyset)$.  Hence $2^\zero$
 is inhabited and contains $\SSC(\emptyset)$.  But $\emptyset$ has only one 
 subset, namely $\emptyset$, which is a separable subset, so $\SSC(\emptyset) = \{\emptyset\} = \zero$.  Thus $2^\zero = \Nc{\zero}  = \one$.
\end{proof}
 
\begin{lemma} \label{lemma:exp_one}
$2^{\one} = \two$.
\end{lemma}

\noindent\begin{proof}
 $\one$ is the set of all singletons.  It therefore
 contains $\{ \zero \}  = \USC(\zero) $.  Then $2^\one$ contains 
 $\SSC(\zero)$.  There are exactly two subsets of $\{ \emptyset\}$,
 namely $\emptyset$ and $\{ \emptyset\}$, and both are separable.  Hence
 $2^\one$ contains the two-element set $ \SSC(\zero) = \{ \emptyset, \{ \emptyset \} \}$.
 That set belongs to $\two = \one^+$ since it is
 equal to $\{ \{ \emptyset \}\} \cup \{\emptyset\}$, and the singleton
 $\{\{\emptyset\}\}$ belongs to $\one$ and $\{\emptyset\} \not\in \{\{\emptyset\}\}$.
 Therefore $2^\one$ and $\one^+$ have a common element.   
 Both are cardinals, by Lemma~\ref{lemma:finiteexp}.  Then by 
 Lemma~\ref{lemma:cardinalsdisjoint}, $2^\one = \two$.  
 \end{proof}
  
 \begin{lemma}\label{lemma:exp_two} $2^{\two} = \four$.
 \end{lemma}
 
 \noindent\begin{proof} By definition, $\four = \three^+ = \two^{++}$.
 One can show (but we omit the details) that
\begin{eqnarray*}
 \USC(\{\one,\two\}) = \{\{\one\},\{\two\}\} \in \two &&
 \end{eqnarray*}
 Therefore, by the definition of exponentiation,
 \begin{eqnarray*}
\SSC(\{\one,\two\}) = \{\emptyset, \{\one\},\{\two\},\{\one,\two\}\} \in 2^\two &&
\end{eqnarray*}
One can explicitly exhibit the ordered pairs of a similarity between
the last-mentioned set and  the element $\{\one,\two,\three,\four\}$ of $\four$.
We omit the details.
Then by Lemma~\ref{lemma:finitecardinals0}, $2^\two = \four$.
\end{proof}

\begin{lemma} \label{lemma:two_members} We have
$$ u \in \two \iff  \exists a,b\,(a \neq b \ \land \ u = \{a, b\}). $$
\end{lemma}

\noindent\begin{proof} We have $\two = \one^+  $.  
If $a \neq b$ then by Lemma~\ref{lemma:one_members}, $ \{a\} \in \one$,
and $\{a\} \cup \{b\} = \{a, b\} \in \two$.  Conversely,  
If $u \in \two$ then $u = v \cup \{b\}$, where $v \in \one$ and $b \not\in v$.
By Lemma~\ref{lemma:one_members}, $v = \{a\}$ for some $a$, so 
$u = \{a,b\}$.
\end{proof}

\begin{lemma} \label{lemma:three_members} We have
$$ u \in \three \iff  \exists a,b,c\,(a \neq b \ \land b \neq c \ \land a \neq c \ \land \ u = \{a, b, c\}). $$
\end{lemma}

\noindent\begin{proof} We have $\three = \two^+  $.  Assume $a,b,c$ are
pairwise distinct.  Then 
 by Lemma~\ref{lemma:two_members}, $ \{a, b\} \in \two$.  Since
 $\three = two^+$,  $\{a,b\} \cup \{c\} = \{a,b,c\} \in \three$.
 Conversely,  
If $u \in \three$ then $u = v \cup \{c\}$, where $v \in \two$ and $c \not\in v$.
By Lemma~\ref{lemma:two_members}, $v = \{a,b\}$ for some $a,b$ with $a \neq b$.
Since $c \not\in v$,  $a \neq c$ and $b \neq c$.  Therefore  
$u = \{a,b,c\}$ with $a,b,c$ pairwise distinct. 
\end{proof}

\begin{lemma} \label{lemma:smallarith} We have $\zero < \one < \two < \three < \four$.
\end{lemma}

\noindent\begin{proof}  Since each of these numbers is defined as the 
successor of the one listed just before it,  the lemma is a 
consequence of Lemma~\ref{lemma:xlessthansuccessorx}.  
\end{proof}

\begin{lemma}\label{lemma:lessthanone} For $m \in \F$, we have
$m < \one \iff m = \zero$.
\end{lemma}

\noindent\begin{proof}  Let $m \in \F$ and $ m < \one$.  By 
Theorem~\ref{theorem:finitetrichotomy}, $ m < \zero \ \lor \  m = \zero \ \lor \ \zero < m$.  By Lemma~\ref{lemma:xnotlessthanzero}, $m < \zero$ is
ruled out.  It remains to rule out $\zero < m$.  Assume $\zero < m$.
Then 
\begin{eqnarray*}
m < \one  &&\mbox{\qquad by hypothesis}\\
m^+ \le \one  &&\mbox{\qquad by Lemma~\ref{lemma:noinsertions}}\\
m^+ \le \zero^+ &&\mbox{\qquad since $\zero^+ = \one$}\\
m \le \zero  &&\mbox{\qquad by Lemma~\ref{lemma:ordersuccessor}}\\
\zero < \zero && \mbox{\qquad by Lemma~\ref{lemma:le_transitive2}, since $\zero \le m <  \zero$}\\
\zero \not< \zero &&\mbox{\qquad by Lemma~\ref{lemma:xnotlessthanzero}}
\end{eqnarray*}
\end{proof}

\begin{lemma} \label{lemma:lessthantwo} For $m \in \F$, we have
$$ m < \two \iff m = \zero \ \lor\  m = \one .$$
\end{lemma}

\noindent\begin{proof} Left to right.  Assume $m \le \two$.
 We have 
\begin{eqnarray*}
\zero \neq \one &&\mbox{\qquad by Lemma~\ref{lemma:xnotequalsuccessorx}}\\
 \{\zero,\one\} \in \two &&\mbox{\qquad by Lemma~\ref{lemma:two_members}}\\
m \le \two && \mbox{\qquad by Definition~\ref{definition:order}}\\
a \in m \land a \in \SSC( \{\zero,\one\})
           && \mbox{\qquad for some $a$, by Lemma~\ref{lemma:le2}}\\
a \neq \{\zero,\one\} && \mbox{\qquad since $m \neq \two$ }\\
\zero \in a \lor \zero\not\in a && \mbox{\qquad since $a \in \SSC(\{\zero,\one\})$}\\
\one \in a \lor \one \not\in a && \mbox{\qquad since $a \in \SSC(\{\zero,\one\})$}
\end{eqnarray*}
An argument by cases (about 170 steps, which we omit) 
shows that $a = \emptyset$, or $a = \{\zero\}$,
or $a = \{\one\}$.  Then $a = \zero$ or $\one$, by Lemma~\ref{lemma:cardinalsdisjoint}.  That completes the left to 
right direction.
\smallskip

{\em Right to left}.   we have $\zero < \two$ and $\one < \two$ 
by Lemma~\ref{lemma:smallarith}.  \end{proof}

\begin{lemma} \label{lemma:usc_unitclass}
For all $a$,  $a$ is a unit class if and only if $\USC(a)$ is a unit class.
\end{lemma}
 
\noindent\begin{proof} Left to right.  Suppose $a = \{x\}$. 
Then the only unit subset of $a$ is $\{a\}$, so $\USC(a)$ is a unit class.
\smallskip
 
Right to left.  Suppose $\USC(a) = \{ u\}$.  Then $u \in a$.
Let $t \in a$. Then $\{t\}\in \USC(a)$,
so $\{t\} = u$.  Hence every element of $a$ is equal to $u$.  Hence
$a = \{u\}$.  \end{proof}  

\begin{lemma}\label{lemma:le_zero} For all $x \in \F$, $x \le \zero \imp x = \zero$.
\end{lemma}

\noindent\begin{proof}  Suppose $x \in\F$ and $x \le \zero$. 
By the definition of $\le$, there exists $a,b$ such that 
$a \in x$, $b \in \zero$, $a \subseteq b$, and $b = (a \cup b) - a$.
Then 
\begin{eqnarray*}
b = \emptyset &&\mbox{\qquad by definition of $\zero$}\\
a = \emptyset &&\mbox{\qquad since $a \subseteq b$}\\
\emptyset \in x \cap \zero && \mbox{\qquad by definition of $\cap$}\\
x = \zero &&\mbox{\qquad by Lemma~\ref{lemma:cardinalsdisjoint}}
\end{eqnarray*}
\end{proof}
 
\begin{lemma}[Specker~4.6] \label{lemma:mlessthanexpm}\ \\
If $m$ is a finite cardinal and $2^m $ is inhabited, then $m < 2^m$.
\end{lemma}

\noindent{\em Remark}. This version of Specker~4.6 phrases the matter
positively, so it is constructively stronger.
\smallskip

\noindent\begin{proof}
Since $\FregeN$ has decidable equality, by Corollary~\ref{lemma:FregeNdecidable} 
$$m = \zero \ \lor \ m = \one \lor (m \neq \zero \ \land \ m \neq \one.)$$
We argue by cases.
\smallskip

Case 1, $m = \zero$.  Then by Lemma~\ref{lemma:exp_zero}, $2^m = \one$,
and we have to show $\zero < \one$, which follows from the definition of $<$
by exhibiting the separable subset $\emptyset$ of the set $\{\emptyset\}$, and
noting that $\emptyset \in \zero$ while $\{\emptyset\} \in \one$.
\smallskip

Case 2, $m = \one$.  Then by Lemma~\ref{lemma:exp_one}, $2^m = \two$,
and we have to show $\one < \two$, which follows from 
$\zero < \one $ by Corollary~\ref{lemma:cardinalsinhabited} 
and Lemma~\ref{lemma:ordersuccessor},
 or more directly, from the definition of $<$
by exhibiting the separable subset $\{\emptyset\}$ of $\{\{ \emptyset\}, \emptyset\}$,
the former of which belongs to $\one$ while the latter belongs to $\two$.
\smallskip

Case 3, $m \neq \zero$ and $m \neq \one$.
By hypothesis, $2^m$ is inhabited.  Then there exists $a$ such 
that $\USC(a) \in m$.
 Since $m \in \FregeN$, we have
\begin{eqnarray*}
\USC(a) \in \FINITE &&\mbox{\qquad by Lemma~\ref{lemma:finitecardinals1}}\\
a \in \FINITE &&\mbox{\qquad by Lemma~\ref{lemma:uscfinite}}\\
\SSC(a) \in \FINITE &&\mbox{\qquad by Lemma~\ref{lemma:finitepowerset}} \\
a \in \DECIDABLE  && \mbox{\qquad by Lemma~\ref{lemma:finitedecidable}} 
\end{eqnarray*}
Then by Lemma~\ref{lemma:finiteseparable}, $\USC(a)$ is a separable subset
of $\SSC(a)$.  Now with $u = \USC(a)$ and $v = \SSC(a)$ 
we have proved that $u$ is a separable subset of $v$ and 
$u \in m$ and $v \in 2^m$.  Then by Definition~\ref{definition:order}
we have $m \le 2^m$. 
\smallskip

By definition $m < 2^m$ means $m \le 2^m$ and $m \neq 2^m$.
It remains to prove that $m \neq 2^m$.  Suppose $m = 2^m$.
As just proved, we have $\USC(a) \subseteq \SSC(a)$.
I say that it is a proper subset, $\USC(a)\subset \SSC(a)$.  It
suffices to prove $\USC (a)\neq \SSC(a)$. We have to produce an 
element of $\SSC(a)$ that does not belong to $\USC(a)$.  We 
propose $a$ as this element.   We have $a \in \SSC(a)$ since 
$a$ is a separable subset of itself.  It
remains to show that $a \not\in \USC(a)$.  Assume $a \in \USC(a)$.
Then $a$ is a unit class.  By Lemma~\ref{lemma:usc_unitclass}, 
$\USC(a)$ is also a unit class.  Any two unit classes are similar, 
so $\USC (a) \sim \zero$.  Since $\zero \in \one$,
$\USC(a) \in \one$, by Lemma~\ref{lemma:finitecardinals0}. 
Then $m \cap \one$ is inhabited, since it contains $\USC(a)$.
Then by Lemma~\ref{lemma:cardinalsdisjoint}, $m = \one$, contradiction.
That completes the proof that $\USC(a)$ is a proper subset of $\SSC(a)$. 
\smallskip

We have 
\begin{eqnarray*}
\USC(a) \subset \SSC(a) && \mbox{\qquad as proved above}\\
\USC(a) \sim \SSC(a)  && \mbox{\qquad by Lemma~\ref{lemma:finitecardinals2},  
since $\USC(a) \in m$ and $\SSC(a) \in 2^m$}\\
\SSC(a) \mbox{\ is infinite} && \mbox{\qquad since $\SSC(a) \sim \USC(a) \subset \SSC(a)$}\\
\neg\,(\SSC(a) \in \FINITE) &&\mbox{\qquad by Theorem~\ref{theorem:infiniteimpliesnotfinite}}\\
\SSC(a)  \in \FINITE && \mbox{\qquad by Lemma~\ref{lemma:finitepowerset}, since $a\in \FINITE$}
\end{eqnarray*}
 That is a contradiction.
 \end{proof}
 
 \begin{lemma}\label{lemma:mplusone_le_expm} For all $m \in \FregeN$,
 $$ \exists u\,( u \in 2^m) \imp m^+ \le 2^m.$$
 \end{lemma}
 
 \noindent\begin{proof} Suppose $m \in \FregeN$ and $\exists u\,( u \in 2^m)$.  Then
 \begin{eqnarray*}
 m < 2^m  &&\mbox{\qquad by Lemma~\ref{lemma:mlessthanexpm}}\\
 2^m \in \FregeN &&\mbox{\qquad by Lemma~\ref{lemma:finiteexp}}\\
 m^+ \le 2^m &&\mbox{\qquad by Lemma~\ref{lemma:noinsertions}}
 \end{eqnarray*}
 \end{proof}

 \begin{lemma}[Specker~4.8] \label{lemma:exporder}
 Let $m,n \in \FregeN$. 
 If $m \le n$ and $2^n$ is inhabited, then $2^m$ is inhabited and 
 $2^m \le 2^n$.
 \end{lemma}

 \noindent\begin{proof} 
 Suppose $m \le n$ and $2^n$ is inhabited.  Then
 \begin{eqnarray*}
 \exists u\,(u \in n) &&\mbox{\qquad by Corollary~\ref{lemma:cardinalsinhabited}}\\
 \exists b\, (\USC(b) \in n) &&\mbox{\qquad by Lemma~\ref{lemma:exp_inhabited}}\\
 \exists b\, (\USC(b) \in m) &&\mbox{\qquad by Lemma~\ref{lemma:exp_inhabited}}
 \end{eqnarray*}
    Since $m\le n$,
  by Lemma~\ref{lemma:le2} there
 is a separable subset $x$ of $\USC(b)$ such that $x \in m$.
   Let $a = \bigcup x$.
 Then using the definitions of $\bigcup$ and $\USC$, we have
 $x = \USC(a)$.  Therefore $2^m$ is inhabited.  Now 
 $2^m = \Nc{\SSC(a)}$ and $2^n = \Nc{\SSC(b)}$. 
 
 I say that $b$ is finite. We have
 \begin{eqnarray*}
 \USC(b) \in n && \\
 \USC(b) \in \FINITE && \mbox{\qquad by Lemma~\ref{lemma:finitecardinals1} }\\
 b \in \FINITE && \mbox{\qquad by Lemma~\ref{lemma:uscfinite} }
 \end{eqnarray*}
 
 I say that $a$ is also finite.  We have
 \begin{eqnarray*}
x \in \FINITE && \mbox{\qquad by Lemma~\ref{lemma:finitecardinals1}, since $x \in m$  }
 \end{eqnarray*}
 Every member of $x$ is a unit class, since $x = \USC(a)$).   
 Every unit class is finite.  Therefore every member of $x$
 is finite.  Moreover, since the members of $x$ are unit classes,
 distinct members of $x$ are disjoint.   Since $x$ is also finite, 
  $a = \bigcup x$ is a finite union of disjoint  
 finite sets.  Hence $a$ is finite, by Lemma~\ref{lemma:finiteunion}.
 \smallskip

  Since $x = \USC(a)$ is 
 a separable subset of $\USC(b)$,  we have
 \begin{eqnarray*}
  \USC(a) \in \SSC(\USC(b)) &&  \\ 
 a \in \SSC(b)   && \mbox{\qquad by Lemma~\ref{lemma:ssc_subset1} (right to left)}\\
 \SSC(a) \in \SSC(\SSC(b)) && \mbox{\qquad by Lemma~\ref{lemma:ssc_subset3}, since $a$ and $b$ are finite}
 \end{eqnarray*}
 Then
 $\SSC(a)$ belongs to $2^m$, and is a separable subset of $\SSC(b)$,
 which belongs to $2^n$.  
 Therefore, by Definition~\ref{definition:order},  $2^m \le 2^n$. 
 \end{proof} 

\section{Addition}

 Specker uses addition in \sect5 of his paper, and relies on Rosser
for its associativity and commutativity.  Those properties can be proved
(as is very well-known) by induction from the two fundamental 
``defining equations'': 
\begin{eqnarray*}
x + y^+ &=& (x+y)^+  \\
x + \zero &=& x 
\end{eqnarray*}
In the present context, where the main point of the paper is to prove that 
there are infinitely many finite cardinals, we need to bear in mind the possibility
that successor or addition may ``overflow''.  We have arranged that successor
is always defined (for any argument whatever); and if there is a largest
natural number then when we take its successor we get the empty set, which 
can be thought of as the computer scientist's ``not a number.''  We need to 
define addition with similar behavior;  if $x+y$ should ``overflow'', it 
should produce ``not a number'', but still be defined.  Then the equations 
above should be valid without further qualification, i.e.,  without insisting
that $x$ and $y$ should be members of $\FregeN$.  If we assume only that those 
equations are valid for $x,y \in \FregeN$,  then the inductive proofs of 
associativity and commutativity do not go through.  

The proofs of associativity and commutativity proceed via another important
property, ``successor shift'':
\begin{eqnarray*}
x^+ + y = x + y^+
\end{eqnarray*}
Normally this property is proved by induction from the ``defining equations.''
In the present context, that does not work, because if $x$ and $y$ are restricted
to $\FregeN$,  then when we try to use successor shift to prove the associative law,
we need $x+y \in \FregeN$, which we do not want to assume, as the statement of the 
associative law should cover the case when $x+y$ overflows.   Therefore, 
we prove below that successor shift is generally valid, i.e., without restricting
$x$ and $y$ to $\FregeN$.  Once we have these {\em three} equations generally 
valid, then the usual proofs of associativity and commutativity by induction
go through 
without difficulty. But in fact, it is simpler and more general to verify them 
directly from the definition of addition, and then we have associativity and
commutativity of addition for all sets, not just finite cardinals. 

\begin{definition}[Specker 3.1, Rosser\protect\footnotemark]\label{definition:addition}%
For any sets $x$ and $y$ we define
$$ x + y := \{ z : \exists u,v\,(u \in x \ \land \ v \in y \ \land \ 
                                  u \cap v = \emptyset \ \land \ z = u \cup v) \}
$$
\end{definition}
\footnotetext{Theorem~XI.2.9 of \cite{rosser1978}, p.~323}
The formula in the definition is stratified, giving $u$, $v$, and $z$ index 1
and $x$ and $y$ index 2.   Then $x,y$, and $z$ all get the same index,
so addition is definable as a function in \INF.  (See Definition~\ref{definition:triples} for ordered triples.)  

\begin{lemma} \label{lemma:addition2} Addition satisfies the ``defining equations''
and successor shift:
\begin{eqnarray*}
x + \zero &=& x \\
x + y^+ &=& (x+y)^+ \\
x + y^+ &=& x^+ + y 
\end{eqnarray*}
\end{lemma}

\noindent{\em Remark}.  Addition is defined on any arguments, not just on 
$\F$.
\medskip

\noindent\begin{proof}  Ad $x + \zero = x$.  By extensionality, it suffices 
to show $z \in x + \zero \iff z \in x$. 
\smallskip 

{\em Left to right}\,: suppose $z \in x + \zero$.
Then $z = u \cup v$, where $u$ and $v$ are disjoint and $u \in x$ and $v \in \zero$.
Since $\zero = \{\emptyset\}$, we have $v = \lambda$, so $z = u \cup \emptyset = u \in x$.
That completes the left to right implication.  
\smallskip

{\em Right to left}\,: Let $z \in x$. Then $z \cup \emptyset \in x + \zero$, by the 
definition of addition.  Since $z \cup \emptyset = z$, we have $z \in x + \zero$
as desired.  That completes the proof of $x + \zero = x$.
\smallskip

Ad $x + y^+ = (x+y)^+$.  By extensionality, it suffices to show the two sides
have the same members.  
\smallskip

{\em Left to right}\,:  We have
\begin{eqnarray*}
 z \in x+y^+  && \mbox{\qquad assumption} \\
 z = u \cup v &&\mbox{\qquad where $u \in z$ and $v \in y^+$ and $u\cap v = \emptyset$} \\
v = w \cup \{c\} &&\mbox{\qquad where $w \in y$ and $c \not\in w$, by definition of $y^+$} \\
z = (u \cup w) \cup \{c\} &&\mbox{\qquad by associativity of union} \\
u \cup w \in x+y  && \mbox{\qquad by definition of addition} \\
c \not \in u \cup w  && \mbox{\qquad since $c \not\in w$ and $u \cap v = \emptyset$} \\
z \in (x+y)^+  && \mbox{\qquad by definition of successor} 
\end{eqnarray*}
That completes the left to right implication.
\smallskip

{\em Right to left}\,: 
\begin{eqnarray*}
z \in (x+y)^+    &&\mbox{\qquad assumption} \\
z = w \cup \{c\}  && \mbox{\qquad where $c \not \in w$ and $w \in x+y$}\\
w = u \cup v      && \mbox{\qquad where $u \in x$ and $v \in y$ and $u\cap v = \emptyset$}\\
z = u \cup (v \cup \{c\}) &&\mbox{\qquad by the associativity of union} \\
c \not \in v      && \mbox{\qquad since $c \not\in w = u \cup v$}\\
v \cup \{c\} \in y^+  && \mbox{\qquad by definition of successor} \\
u \cap (v \cup \{c\}) = \emptyset && \mbox{ \qquad since  $u \cap v = \emptyset$
                           and $c \not\in u$}\\
u \cup (v \cup \{c\}) \in x + y^+ && \mbox{\qquad by definition of addition} \\
(u \cup v) \cup \{c\} \in x + y^+ && \mbox{\qquad by the associativity of union} \\
z \in x + y^+  && \mbox{\qquad since $z = w \cup \{c\} = (u \cup v) \cup \{c\}$ }
\end{eqnarray*}
That completes the proof of the right to left implication.
That completes the proof of $x + y^+ = (x+y)^+$.
\smallskip

{\em Ad successor shift}.  We must prove 
$$ z \in x + y^+ \iff z \in x^+ + y.$$  
\smallskip

{\em Left to right}\,: 
\begin{eqnarray*}
z \in x + y^+  &&\mbox{\qquad assumption} \\
z = u \cup (v \cup \{c\}) && \mbox{\qquad where $u \in x$, $v \in y$, and $c \not\in v$, and $u \cap (v \cup \{c\}) = \emptyset$} \\
z = (u \cup \{c\}) \cup v  &&\mbox{\qquad by the associativity and commutativity of union }\\
c \not \in u  &&\mbox{\qquad since $u \cap (v \cup \{c\}) = \emptyset$} \\
u \cup \{c\} \in x^+ &&\mbox{\qquad by the definition of successor} \\
(u \cup \{c\}) \cap v = \emptyset &&\mbox{\qquad by the associativity and commutativity of union }\\
z \in x^+ + y  &&\mbox{\qquad by the definition of addition} 
\end{eqnarray*}
That completes the left to right direction.
\smallskip

{\em Right to left}\,: 
\begin{eqnarray*}
z \in x^+ + y &&\mbox{\qquad assumption} \\
z = (u \cup \{c\}) \cup v  &&\mbox{\qquad where $u \in x$, $v \in y$, $c \not \in u$,
and $(u \cup \{c\})\cap v = \emptyset$ } \\
z = (u \cup v) \cup \{c\}   &&\mbox{\qquad by the associativity and commutativity of union }\\
c \not \in u \cup v && \mbox{\qquad since $c \not\in u$ and $(u \cup \{c\}) \cap v \ \emptyset$ }\\
u \cup v \in x + y  && \mbox{\qquad by the definition of addition} \\
z \in (x+y)^+      &&\mbox{\qquad by the definition of successor}
\end{eqnarray*}
That completes the right to left direction.  
\end{proof}

\begin{lemma} \label{lemma:addition3} Addition obeys the 
associative and commutative laws  
and left identity (without restriction to $\FregeN$) 
\begin{eqnarray*}
\zero + x &=& x \\     
(x+y)+z &=& x + (y+z)   \\
x + y &=& y + x 
\end{eqnarray*}
\end{lemma}

\noindent{\em Remark}.  We call attention to the fact that, 
even when $x,y,z$ are assumed to be in $\FregeN$, the expressions
in the equations might ``overflow'',  and the equations contain
implicitly the assertion that the overflows ``match'', i.e., one side
overflows if and only if the other does.  Here ``overflow'' means
to have the value $\emptyset$. 
\smallskip 

\noindent\begin{proof}  These laws are immediate consequences of the 
definition of addition, via the associative and commutative laws of 
set union.  We omit the proofs.
\end{proof}

\begin{lemma} \label{lemma:successorisplusone}
For all $m \in \F$, we have $m^+ = m + \one$.
\end{lemma}

\noindent\begin{proof} We have
\begin{eqnarray*}
m + \one = m + \zero^+  &&\mbox{\qquad by definition of $\one$}\\
m + \one = m^+ + \zero  &&\mbox{\qquad by Lemma~\ref{lemma:addition2}} \\
m + \one = m^+   && \mbox{\qquad by Lemma~\ref{lemma:addition2}}
\end{eqnarray*}
\end{proof} 

\begin{lemma}\label{lemma:oneplusone} $\one + \one = \two$.
\end{lemma}

\noindent\begin{proof} 
\begin{eqnarray*}
\two = \one^+ && \mbox{\qquad by definition of $\two$}\\
\two = \one + \one && \mbox{\qquad by Lemma~\ref{lemma:successorisplusone}}
\end{eqnarray*}
\end{proof}

\begin{lemma} \label{lemma:inhabited_sum} Suppose $\kappa, \mu \in \FregeN$,
and $\kappa + \mu$ is inhabited.  Then $\kappa+ \mu \in \FregeN$.
\end{lemma}

\noindent{\em Remark}.  This lemma addresses the problem of possible ``overflow''
of addition.  If there are enough elements to find disjoint members of 
$\kappa$ and $\mu$ then adding $\kappa$ and $\mu$ will not overflow.
\medskip

\noindent\begin{proof}  By induction on $\mu$, which is legal since the
formula is stratified.
\smallskip

{\em Base case}\,: $\kappa + \zero = \kappa$ is in $\FregeN$ because $\kappa \in \FregeN$.
\smallskip

{\em Induction step}: Suppose $\kappa + \mu^+$ is inhabited and $\mu^+$ is inhabited.
Then $\kappa + \mu^+ = (\kappa + \mu)^+$ is inhabited.  By the induction
hypothesis, $\kappa + \mu \in \FregeN$.  Then by 
Lemma~\ref{lemma:successorF}, $(\kappa+\mu)^+ \in \FregeN$.
Since $(\kappa + \mu)^+ = \kappa + \mu^+$, we have $\kappa + \mu^+ \in \FregeN$.
That completes the induction step. 
\end{proof}
\smallskip

\begin{lemma} \label{lemma:subterms}
 Suppose $p,q,r \in \FregeN$ and $p+q+r \in \FregeN$.  Then $p+q$ and $q+r$ 
are also in $\FregeN$.  Similarly, if $p,q,r,s \in \FregeN$ and $p+q+r+s \in \FregeN$,
then $p + q + r \in \FregeN$.  
\end{lemma}

\noindent\begin{proof}   By Corollary~\ref{lemma:cardinalsinhabited},
$p+q+r$ is inhabited.  Let $u \in p+q+r$.  Then by the definition 
of addition, $u = a \cup b \cup c$ with $a \in p$, $b \in q$, $c \in r$,
and $a,b,c$ pairwise disjoint. Then $a \cup b \in p +q$ and $b \cup c \in q+r$.
Then by Lemma~\ref{lemma:inhabited_sum}, $p+q\in \FregeN$ and $q+r \in \FregeN$.
That completes the proof of the three summand case.  The case of four
summands is treated similarly.  We omit the details.  
\end{proof}
 
 \begin{lemma} \label{lemma:subterms2}
If $p \in \FregeN$  and  $p+q^+ \in \FregeN$,
then $p^+  \in \FregeN$.
\end{lemma}

\noindent{\em Remark}. It is not assumed that $q \in \FregeN$.
\medskip 

\noindent\begin{proof} 
 Suppose $p \in \FregeN$ and $p+q^+ \in \FregeN$.  By 
Corollary~\ref{lemma:cardinalsinhabited}, there exists $u \in p + q^+$.
 Then by Definition~\ref{definition:addition}, 
there exist $a$ and $b$ with $a\in p$ and $b \in q^+$ and $a \cap b = \emptyset$.
By definition of successor, $b = x \cup \{c\}$ for some $x$ and $c$, 
so $c \in b$.  Since $a \cap b = \emptyset$, we have $c \not \in a$.
 Then $a \cup \{ c\} \in p^+$.  Then $p^+ \in \FregeN$. 
\end{proof}

 \begin{lemma} \label{lemma:subterms3}
If $p,q \in \FregeN$  and  $p+q^+ \in \FregeN$,
then $p + q  \in \FregeN$.
\end{lemma}

\noindent\begin{proof}  We have 
\begin{eqnarray*}
p + q^+ \in \F &&\mbox{\qquad by hypothesis}\\
p + q + \one \in \F && \mbox{\qquad by definition of $\one$ and Lemma~\ref{lemma:addition2}}\\
p+ q \in \F &&\mbox{\qquad by Lemma~\ref{lemma:subterms}}
\end{eqnarray*}
\end{proof} 

\begin{lemma}\label{lemma:addorder}
For $a,b,p,q \in \F$,  if $b + q \in \F$ we have
$$a \le b \ \land\ p \le q \imp a + p \le b + q$$
\end{lemma} 

\noindent\begin{proof}  Suppose $a,b,p,q \in \F$ and $b + q \in \F$.
Suppose also  $a \le b$, $p \le q$.   Then
\begin{eqnarray*}
w \in b+q && \mbox{\qquad for some $w$, by 
Corollary~\ref{lemma:cardinalsinhabited}, since $b+q \in \F$}
\end{eqnarray*}
By the definition of addition, there exist $u,v$ with $w = u \cup v$,
$u\in b$, $v\in q$, and
$u \cap v = \emptyset$. 
By Lemma~\ref{lemma:le2}, since $a \le b$ there exists $r \in a$ with 
$r \in \SSC(u)$. By Lemma~\ref{lemma:le2}, since $p \le q$, there exists $s \in p$ with $s \in \SSC(v)$.
Then one can verify that $r \cup s \in \SSC(u \cup v)$. (We omit the details of 
that verification.)  Since $u \cup v = w$ we have $r \cup s \in \SSC(w)$.
We have $r \cap s = \emptyset$, since $r \subseteq u$, $s \subseteq v$,
and $u \cap v = \emptyset$.  Then $r \cup s \in a + p$, by the definition of 
addition.  Then $a + p \le b + q$, as witnessed by $r \cup s \in a + p$,
$r \cup s \in \SSC(w)$, and $w \in b+q$. 
\end{proof}

\begin{lemma}\label{lemma:addorder2}
For $a,p,b,q\in \F$,  if $b + q \in \F$ we have
$$a < b \ \land\ p \le q \imp a + p < b + q$$
\end{lemma} 

\noindent{\em Remark.} It is not assumed that $a+p\in \F$,
which would make the proof easier.
\smallskip

\noindent\begin{proof} 
Suppose $a < b$ and $p \le q$.  Then 
\begin{eqnarray*}
a^+ \le b &&\mbox{\qquad by Lemma~\ref{lemma:noinsertions}}\\
\exists u\, (u \in a^+) &&\mbox{\qquad by the definition of addition}\\
a^+ \in \F  &&\mbox{\qquad by Lemma~\ref{lemma:successorF}}\\
a^+ + p \le b + q  &&\mbox{\qquad by Lemma~\ref{lemma:addorder}}\\
(a+p)^+ \le b + q  &&\mbox{\qquad by Lemma~\ref{lemma:addition3}}\\
\exists u\,(u \in (a+p)^+)  &&\mbox{\qquad by the definition of $\le$}\\
\exists u\,(u \in a+p) &&\mbox{\qquad by definition of successor}\\
a+p \in \F  &&\mbox{\qquad by Lemma~\ref{lemma:inhabited_sum} }\\
(a+p)^+ \in \F   &&\mbox{\qquad by Lemma~\ref{lemma:successorF}}\\  
a+ p < (a+p)^+  &&\mbox{\qquad by Lemma~\ref{lemma:xlessthansuccessorx}}\\
a+p < b+q   &&\mbox{\qquad by Lemma~\ref{lemma:le_transitive2}}
\end{eqnarray*}
\end{proof}

\begin{lemma} \label{lemma:exp_members2} For $m \in \FregeN$ we have 
$$ \USC (x) \in m \imp  \SSC (x) \in 2^m.$$
\end{lemma}

\noindent\begin{proof}  Suppose $\USC(x) \in m$.  By Definition~\ref{definition:exp},
$2^m$ contains all sets similar to $\SSC(x)$.   By Lemma~\ref{lemma:sim},
$\SSC(x)$ is one of those sets, so $\SSC(x) \in 2^m$. 
\end{proof}

\begin{lemma} \label{lemma:expnotzero}
 For all $z$ we have $2^z \neq \zero$.
\end{lemma}

\noindent\begin{proof} Suppose $2^z = \zero$.  Then 
\begin{eqnarray*}
\emptyset \in \zero  &&\mbox{\qquad by Definition~\ref{definition:Fregeonetwo}}\\
\emptyset \in 2^z    && \mbox{\qquad since $2^z = \zero$}\\
\emptyset~\sim \SSC(a) \ \land \ \USC(a) \in x &&\mbox{\qquad by Definition~\ref{definition:exp}}\\
\SSC(a) = \emptyset   && \mbox{\qquad  since only $\emptyset$ is similar to $\emptyset$}
\end{eqnarray*}
But $a \in \SSC(a)$, contradiction.  \end{proof}

\begin{lemma}\label{lemma:extend_similar} Suppose $x \sim y$, 
and $a \not\in x$ and $b \not \in y$.  Then 
\begin{eqnarray*}
 x \cup \{ a\} \sim y \cup \{ b\}.
 \end{eqnarray*}
\end{lemma}

\noindent\begin{proof}  Extend a similarity $f:x \to y$ by defining $f(a) = b$.
We omit the details. 
\end{proof}

\begin{lemma} \label{lemma:cardinality_additive}
Let $p$ and $q$ be disjoint finite sets. 
 Then $\Nc{p \cup q} = \Nc{p} + \Nc{q}$.
\end{lemma}

\noindent\begin{proof} We have 
\begin{eqnarray*}
p \cup q \in \FINITE &&\mbox{\qquad by Lemma~\ref{lemma:union}} \\
\Nc{p \cup q} \in \F &&\mbox{\qquad by Lemma~\ref{lemma:finitecardinals3}}\\
\Nc{p} \in \F &&\mbox{\qquad by Lemma~\ref{lemma:finitecardinals3}}\\
\Nc{q} \in \F &&\mbox{\qquad by Lemma~\ref{lemma:finitecardinals3}}\\
p \cup q \in \Nc{p \cup q} &&\mbox{\qquad by Lemma~\ref{lemma:xinNcx}}\\
p \in \Nc{p} &&\mbox{\qquad by Lemma~\ref{lemma:xinNcx}}\\
q \in \Nc{q}   &&\mbox{\qquad by Lemma~\ref{lemma:xinNcx}}\\
p \cup q \in \Nc{p} + \Nc{q} &&\mbox{\qquad by the definition of addition}\\
\Nc{p \cup q} \cap \Nc{p} + \Nc{q} \neq \emptyset &&\mbox{\qquad since both contain $p \cup q$}\\
\Nc{p} + \Nc{q} \in \F &&\mbox{\qquad by Lemma~\ref{lemma:inhabited_sum}} \\
\Nc{p \cup q} = \Nc{p} + \Nc{q} &&\mbox{\qquad by Lemma~\ref{lemma:cardinalsdisjoint} }
\end{eqnarray*} 
\end{proof}

\begin{lemma} \label{lemma:subtraction}
For $p,q,r \in \F$, if $q+p \in \F$  we have
\begin{eqnarray*}
 q + p &=&r + p \imp q = r \\
 p + q &=& p + r \imp q = r.
\end{eqnarray*}
\end{lemma}

\noindent\begin{proof}  The two formulas are equivalent, by Lemma~\ref{lemma:addition3}.
We prove the first one by induction on $p$, which is legal since the formula 
is stratified.  More precisely we prove by induction on $p$ that 
$$ \forall q,r \in \F\,(q+p \in \F \imp q + p = r + p \imp q = r).$$
{\em Base case}, $p=0$.  Suppose $q + 0 = r + 0$.  Then $q=r$ by the right identity
property of addition, Lemma~\ref{lemma:addition2}.  That completes the base case.
\smallskip

{\em Induction step}.  Suppose $q + p^+ = r + p^+$ and $q + p^+ \in \F$.
 Then 
\begin{eqnarray*}
(q+p)^+ = (r+p)^+  &&\mbox{\qquad by Lemma~\ref{lemma:addition2}}\\ 
q+p \in \F         &&\mbox{\qquad by Lemma~\ref{lemma:subterms3} }\\
r+p \in \F         &&\mbox{\qquad  by Lemma~\ref{lemma:subterms3} } \\ 
\exists u\,(u \in q+p) &&\mbox{\qquad by Corollary~\ref{lemma:cardinalsinhabited}}\\
\exists u\,(u \in r+p) &&\mbox{\qquad by Corollary~\ref{lemma:cardinalsinhabited}}\\
(q+p)^+ = q + p^+ &&\mbox{\qquad by Lemma~\ref{lemma:addition2}}\\
(r+p)^+ = r + p^+ &&\mbox{\qquad by Lemma~\ref{lemma:addition2}}\\
(q+p)^+ \in \F     && \mbox{\qquad equality substitution}\\
(r+p)^+ \in \F     && \mbox{\qquad equality substitution} \\
\exists u\,(u \in (q+p)^+) &&\mbox{\qquad by Corollary~\ref{lemma:cardinalsinhabited}}\\
\exists u\,(u \in (r+p)^+) &&\mbox{\qquad by Corollary~\ref{lemma:cardinalsinhabited}}
\end{eqnarray*}
\begin{eqnarray*}   
q+p = r+p          &&\mbox{\qquad by Lemma~\ref{lemma:successoroneone}, since
                       $(q+p)^+ = (r+p)^+$    } \\
q = r              &&\mbox{\qquad by the induction hypothesis}
\end{eqnarray*}
That completes the induction step.  \end{proof}

\begin{lemma} \label{lemma:ssc_adjoin2} Let $b \in \FINITE$ and $c \not\in b$. 
Then 
$$ \Nc{{\SSC (b \cup \{c\})}} = \Nc{{\SSC(b)}} + \Nc{{\SSC(b)}}.$$
\end{lemma}

\noindent\begin{proof}  Define
$$ R := \{ x \cup \{c\} : x \in \SSC(b)\}.$$
The definition can be rewritten in stratified form, so $R$ can be 
defined in \INF.  Define $f: x \mapsto x \cup \{c\}$,  which can also be defined in \INF:
$$ f := \{ \langle x, x \cup \{c\} \rangle: x \in \SSC(b) \}.$$
The formula is stratified, since all the occurrences of $x$ can be given index 0,
and $\{c\}$ and $\SSC(b)$ are just parameters.  
Then $f: \SSC(b) \to R$ is a similarity. (We omit the 150 steps required to prove that.)
\smallskip

We first note that if $x \in \SSC(b \cup \{c\})$  and $c \in x$,  then 
$x = (x-c) \cup \{c\}$, since $x$ is finite and therefore has decidable equality.
Similarly $b \cup \{c\}$ has decidable equality, so 
every $x \in \SSC(b \cup \{c\})$ either contains $c$ or not.  If $c \in x$ then 
$x \in R$.  If $c \not\in x$ then $x \in \SSC(b)$.  Therefore
\begin{eqnarray*}
\SSC(b \cup \{c\}) =  \SSC(b) \cup R && \\
\SSC(b) \sim R &&\mbox{\qquad since $f:\SSC(b) \to R$ is a similarity}\\
\Nc{\SSC(b)} = \Nc{R} &&\mbox{\qquad by Lemma~\ref{lemma:cardinalequality}} \\
\SSC(b) \cap R = \emptyset &&\mbox{\qquad since $c \not \in b$}\\
\SSC(b) \in \FINITE && \mbox{\qquad by Lemma~\ref{lemma:finitepowerset}} \\
R \in \FINITE && \mbox{\qquad by Lemma~\ref{lemma:finitesimilar}} \\
\Nc{\SSC(b \cup \{c\})} = \Nc{\SSC(b)} + \Nc{R} &&\mbox{  by Lemma~\ref{lemma:cardinality_additive}}\\
\Nc{\SSC(b \cup \{c\})} = \Nc{\SSC(b)}+ \Nc{\SSC(b)} &&\mbox{  since $\Nc{\SSC(b)} = \Nc{R}$}
\end{eqnarray*}
\end{proof} 

\begin{lemma}\label{lemma:exprec} For $p \in \F$, if $2^{p^+} \in \F$,
then
 $2^{p^+} = 2^p + 2^p$.
 \end{lemma}

 \noindent\begin{proof}   Suppose $p \in \F$ and $2^{p^+} \in \F$. Then
 \begin{eqnarray}
 \exists u\,(u \in 2^{p^+}) &&\mbox{\qquad by Corollary~\ref{lemma:cardinalsinhabited}}\nonumber\\
 \USC(a) \in p^+ &&\mbox{\qquad for some $a\in p$, by definition of exponentiation}\nonumber\\
 u \in p^+ \land q\in u &&\mbox{\qquad for some $q,u$, by Lemma~\ref{lemma:successorinhabited}}\nonumber\\
 u \sim \USC(a) &&\mbox{\qquad by Lemma~\ref{lemma:finitecardinals2}, since both are in $p^+$}\nonumber\\
w \in \USC(a)   &&\mbox{\qquad for some $w$, by Lemma~\ref{lemma:similarinhabited}}\nonumber\\
c \in a \land w = \{c\}      &&\mbox{\qquad for some $c$,  by definition of $\USC(a)$}\nonumber\\
\USC(a) \in \FINITE   && \mbox{\qquad by Lemma~\ref{lemma:finitecardinals1}}\nonumber\\
a \in \FINITE && \mbox{\qquad by Lemma~\ref{lemma:uscfinite}}\nonumber\\
a \in \DECIDABLE &&\mbox{\qquad by Lemma~\ref{lemma:finitedecidable}}\nonumber
\end{eqnarray}
\begin{eqnarray}   
b := a-\{c\} &&\mbox{\qquad definition of $b$}\nonumber\\
a = b \cup \{c\} && \mbox{\qquad since $a \in \DECIDABLE$}\\
\USC(a) = \USC(b) \cup \{\{c\}\} &&\mbox{\qquad by Lemma~\ref{lemma:usc_successor}}\nonumber\\
\SSC(b \cup \{c\}) \in 2^{p^+} &&\mbox{\qquad by definition of exponentiation}\nonumber\\
\Nc{{\SSC(b \cup \{c\}}} = 2^{p^+} &&\mbox{\qquad by Lemma~\ref{lemma:xinNcx}}\nonumber\\
\Nc{{\SSC(b \cup \{c\})}} = \Nc{{\SSC(b)}} + \Nc{{\SSC(b)}} &&\mbox{\qquad by Lemma~\ref{lemma:ssc_adjoin2}}\label{eq:4076}\\
\USC(a) \in \DECIDABLE &&\mbox{\qquad by Lemma~\ref{lemma:finitedecidable}}\nonumber\\
\USC(b) = \USC(a) - \{\{c\}\} &&\mbox{\qquad since $\USC(a) \in \DECIDABLE$}\nonumber\\
\USC(b) \in p && \mbox{\qquad by Lemma~\ref{lemma:cardinalpredecessor}}\nonumber\\
\SSC(b) \in 2^p &&\mbox{\qquad by the definition of exponentiation}\nonumber\\
\SSC(b) \in \Nc{\SSC(b)} && \mbox{\qquad by Lemma~\ref{lemma:xinNcx}}\nonumber\\
\SSC(b) \in \FINITE && \mbox{\qquad by Lemma~\ref{lemma:finitepowerset}}\nonumber\\
\Nc{{\SSC(b)}} \in \F  && \mbox{\qquad by Lemma~\ref{lemma:finitecardinals3}}\nonumber\\
\Nc{{\SSC(b)}} = 2^p &&\mbox{\qquad by Lemma~\ref{lemma:cardinalsdisjoint}}\nonumber
\end{eqnarray}
Then $ 2^{p^+} = 2^p + 2^p$ as desired, by (\ref{eq:4076}).
\end{proof}

\begin{lemma}\label{lemma:exponeonebase}
For $m \in \F$,  $2^m = \one \iff m = \zero$.
\end{lemma}

\noindent\begin{proof}
Left to right.  We have
\begin{eqnarray*}
 2^m = \one &&\mbox{\qquad assumption}\\
2^\two = \four&&\mbox{\qquad by Lemma~\ref{lemma:exp_two}}\\
\two \le m \imp 2^\two \le 2^m &&\mbox{\qquad by Lemma~\ref{lemma:exporder}}\\
\two \le m \imp \four \le \one &&\mbox{\qquad by transitivity of $\le$}\\
\one < \four &&\mbox{\qquad by Lemma~\ref{lemma:smallarith}}\\
\two \not\le m &&\mbox{\qquad otherwise $\one < \four\ \land\ \four \le \one$}\\
m < \two \lor \two \le m &&\mbox{\qquad by Theorem~\ref{theorem:finitetrichotomy}}\\
m < \two && \mbox{\qquad since $\two \not \le m$}\\
m = \zero \ \lor \ m = \one &&\mbox{\qquad by Lemma~\ref{lemma:lessthantwo}}\\
2^\one = \two  && \mbox{\qquad by Lemma~\ref{lemma:exp_one}}\\
\one \neq \two && \mbox{\qquad by Lemma~\ref{lemma:xnotequalsuccessorx}}\\
m \neq \one  && \mbox{\qquad since $2^m = \zero$} \\
m = \zero && \mbox{\qquad since $m = \zero \lor m = \one$ but $m \neq \one$}
\end{eqnarray*}

Right to left.  Suppose $m = \zero$.  Then $2^m = 2^\zero = \one$,
by Lemma~\ref{lemma:exp_zero}.  \end{proof}

 \begin{lemma}\label{lemma:exponeone}
 For $n,m \in \FregeN$, if $2^n = 2^m$ and $2^n$ is inhabited, then 
 $n = m$.
 \end{lemma}
 
{\em Remark.} The reader is invited to try a direct proof using the 
definition of exponentiation.  It would work if we had the converse
of Lemma~\ref{lemma:sscsimilar}.  The only proof of that converse that
we know requires this lemma.   Therefore,  we give a more complicated
(but correct) proof by induction.
\medskip

  \noindent\begin{proof} We prove by induction on $n$ that for $n \in \F$ 
  with $2^n$ inhabited, we have 
\begin{eqnarray}
\exists u\,(u \in 2^n) \imp \forall m\in \F\,(2^n = 2^m \imp n = m) \label{eq:2934}
\end{eqnarray}
The formula is stratified giving $n$ and $m$ both index 0, so 
it is legal to proceed by induction.
\smallskip

The base case follows from Lemma~\ref{lemma:exponeonebase}.
\smallskip

{\em Induction step}. Suppose $2^{n^+} = 2^m$ and $n^+$ is inhabited. 
We have $m = \zero \ \lor \ m \neq \zero$, by Lemma~\ref{lemma:FregeNdecidable}.
\smallskip

Case~1, $m = \zero$. Then by Lemma~\ref{lemma:exponeonebase}, $n^+ = \zero$, contradiction.
\smallskip

Case~2,  $m \neq \zero$.  Then
\begin{eqnarray*}
\exists r\in \F\,(m = r^+) &&\mbox{\qquad by Lemma~\ref{lemma:nonzeroissuccessor}}\\
 2^{n^+} = 2^{r^+}  &&\mbox{\qquad since $2^{n^+} = 2^m$}\\
 2^n+2^n = 2^r + 2^r&&\mbox{\qquad by Lemma~\ref{lemma:exprec}}\\
 r < n \ \lor \ r = n \ \lor \ n < r &&\mbox{\qquad by Theorem~\ref{theorem:finitetrichotomy}}
 \end{eqnarray*}
We argue by cases.
\smallskip

Case 1, $r < n$.  Then
\begin{eqnarray*}
2^r \le 2^n  &&\mbox{\qquad  by Lemma~\ref{lemma:exporder}}\\
2^r \neq 2^n &&\mbox{\qquad  by the induction hypothesis}\\
2^r < 2^n    &&\mbox{\qquad by the definition of $<$} \\
2^{n^+} \in \F && \mbox{\qquad by Lemma~\ref{lemma:finiteexp}}\\
2^{r^+} = 2^r + 2^r &&\mbox{\qquad by Lemma~\ref{lemma:exprec}}\\
2^r + 2^r  < 2^n + 2^n &&\mbox{\qquad by Lemma~\ref{lemma:addorder2}}\\
2^n+2^n  = 2^{n+} &&\mbox{\qquad by Lemma~\ref{lemma:exprec}}\\
2^{r^+} <  2^{n+} &&\mbox{\qquad by Lemma~\ref{lemma:le_transitive2}}
\end{eqnarray*}
But that contradicts $2^{n^+} = 2^{r^+}$.  That completes Case~1.
\smallskip

Case 2, $n < r$,  similarly leads to a contradiction.  We omit the steps.
\smallskip

Case 3, $n = r$.  Then  $2^n = 2^r$.  Substituting $2^n$ for $2^r$
in the identity $2^r + 2^r = 2^r + 2^r$, we have 
$2^n+2^n = 2^r + 2^r$.  Then $2^{n^+} = 2^{r^+} = 2^m$ as desired.
That completes the induction step.  
\end{proof}
 
\begin{lemma}\label{lemma:exporderstrict} 
 Let $m,n \in \FregeN$. 
 If $m < n$ and $2^n$ is inhabited, then $2^m$ is inhabited and 
 $2^m < 2^n$.
\end{lemma} 

\noindent\begin{proof}  Suppose $m<n$ and $2^n$ is inhabited.  
Then 
\begin{eqnarray*}
m \le n &&\mbox{\qquad by the definition of $<$}\\
2^m \le 2^n && \mbox{\qquad by Lemma~\ref{lemma:exporder}}\\
m \neq n   &&\mbox{\qquad by the definition of $<$} \\
2^m \neq 2^n &&\mbox{\qquad by Lemma~\ref{lemma:exponeone}} \\
2^m < 2^n   && \mbox{\qquad by the definition of $<$}
\end{eqnarray*}
\end{proof}

\begin{lemma} \label{lemma:orderbyaddition} For $p,q \in \F$ we have
$$ p \le q \iff \exists k\in \F\, (p+k = q).$$
\end{lemma}

\noindent\begin{proof}  By  induction on $q$.  The formula is 
stratified, giving all variables index 0.
\smallskip

{\em Base case}, $p \le \zero \iff \exists k\in \F, p + k = \zero$.
{\em Left to right}\,: Suppose $p \le \zero$.  Then $p = \zero \ \lor \ p < \zero$,
by Lemma~\ref{lemma:letolessthan}.  But $p \not< \zero$ by 
Lemma~\ref{lemma:xnotlessthanzero}.  Hence $p = \zero$.  Then 
$p + k = \zero + k = \zero$ by Lemma~\ref{lemma:addition3}.
\smallskip.
{\em Right to left}.  Suppose $p + k = \zero$. Then by the definition
of addition, there exists sets $a\in p$ and $b \in k$ such that 
$a \cup b \in \zero$.  By definition of $\zero$,  $\zero = \{\emptyset\}$,
so $a \cup b = \emptyset$.  Then $a = \emptyset$.  Then $\emptyset \in p$
and $\emptyset \in \zero$.  Then by Lemma~\ref{lemma:cardinalsdisjoint},
$p  = \zero$.  That completes the base case.
\smallskip

{\em Induction step}.  Assume $q^+$ is inhabited.  We have to show
$$ p \le q^+ \iff \exists k\in \F\,(p+k = q^+).$$
{\em Left to right}. suppose $p \le q^+$.  Then $p=q^+ \ \lor \ p \le q$,
by Lemma~\ref{lemma:lessthansuccessor2}.
\smallskip

Case 1, $p \le q$.
Then by the induction hypothesis, there exists $k \in \F$ such that 
$p + k = q$. We have
\begin{eqnarray*}
\exists u\,(u \in q^+) &&\mbox{\qquad by hypothesis}\\
\exists u\, u \in (p+k)^+ &&\mbox{\qquad  since $p+k=q$}\\
p + (k^+) = (p+k)^+ = q^+ &&\mbox{\qquad by Lemma~\ref{lemma:addition2}}\\
\exists u\, (u \in k^+)  && \mbox{\qquad by the definition
of addition}\\
k^+ \in \F &&\mbox{\qquad by Lemma~\ref{lemma:successorF}}
\end{eqnarray*}
That completes Case~1. 
\smallskip

Case 2, $p = q^+$.  Then taking $k=\zero$ we have 
$$p+k = p + \zero = p  = q^+.$$
That completes Case~2.
\smallskip

Right to left.  Suppose $k \in \F$ and $p+k = q^+$. We have to
show $p \le q^+$.  By definition of addition, there exist $a$ and $b$
with $a \in p$ and $b \in k$, and $a \cap b = \emptyset$ and $a \cup b \in q^+$.
Then $a$ is a separable subset of $a \cup b$, so $p \le q^+$ by 
the definition of $\le$.  That completes the induction step.
\end{proof}

\begin{lemma} \label{lemma:xlessthan_xplusy}
Let $p,q \in \F$ and $p + q \in \F$.  Then $p \le p+q$ and $q \le p+q$.
\end{lemma}

\noindent\begin{proof}  Suppose $p, q \in \F$ and $p+q \in \F$.  We have
\begin{eqnarray*}
u \in q && \mbox{\qquad for some $u$, by Corollary~\ref{lemma:cardinalsinhabited}}\\
\emptyset \in \zero &&\mbox{\qquad by the definition of $\zero$}\\
\emptyset \subseteq u \ \land \ u = \emptyset \cup (u - \emptyset) && \mbox{\qquad by the definitions of subset and difference}\\
\zero \le q  && \mbox{\qquad by the definition of $\le$}\\
p \le p  && \mbox{\qquad by Lemma~\ref{lemma:le_reflexive}} \\
p + \zero \le p + q && \mbox{\qquad by Lemma~\ref{lemma:addorder}}\\
p \le p + q &&\mbox{\qquad by Lemma~\ref{lemma:addition2}}
\end{eqnarray*}
That is the first assertion of the lemma.  By Lemma~\ref{lemma:addition3},
we have $p+q = q+p$, so $q + p \in \F$ and as above we have $q \le q+p$.
Therefore also $q \le p+q$.  \end{proof}

 \begin{lemma}\label{lemma:mplusonelessthanexpm}    
Let $p \in \F$.  Then 
$$ p \neq \zero \imp p \neq \one  \imp 2^p \in \F \imp p^+ < 2^p.$$
\end{lemma}

\noindent{\em Remark.}  Specker~4.6 says $p < 2^p$.  Of course the exponent grows faster than 
linearly, so larger things can be put on the  left side, at the price of small exceptions.
\smallskip

\noindent\begin{proof}  By induction on $p$.  For the base case, there is nothing to prove.   
 For the induction step, assume $p^+$ is inhabited and $2^{p^+} \in \F$ and 
 $p^+ \neq \zero$ and $p^+ \neq \one$.
 We have to prove
 \begin{eqnarray}
 p^{++} < 2^{p^+}  \label{eq:4262}
 \end{eqnarray}
 We have
\begin{eqnarray*}
p \neq \zero  && \mbox{\qquad since $p^+ \neq \one$}
\end{eqnarray*}
Since equality on $\F$ is decidable, $p = \one \ \lor \ p \neq \one$,  
\smallskip

{\em Case 1},  $p = \one$.  Then
\begin{eqnarray*}
 p^{++} = \two^+ = \three  && \mbox{\qquad by definitions of $\two$ and $\three$}\\  
2^{p^+} = \four && \mbox{\qquad by Lemma~\ref{lemma:exp_two}}\\  
\three < \four   && \mbox{\qquad by Lemma~\ref{lemma:xlessthansuccessorx}}\\   
p^{++} < 2^{p+} && \mbox{\qquad since $p^{++} = \three$ and $2^{p^+} = \four$}
\end{eqnarray*}
That completes the case $p = \one$.
\smallskip

{\em Case 2}, $p \neq \one$.    
Then
\begin{eqnarray}
p \neq \zero   && \mbox{\qquad since $p^+ \neq \one$ by hypothesis} \label{eq:h19}\\  
2^{p^+} \in \F  && \mbox{\qquad by hypothesis}  \nonumber \\  
2^{p^+} = 2^p + 2^p  && \mbox{\qquad by Lemma~\ref{lemma:exprec}} \label{eq:h20} \\  
p^+ \in \F     && \mbox{\qquad by Lemma~\ref{lemma:successorF}} \nonumber \\  
2^p < 2^{p+}    &&  \mbox{\qquad by Lemma~\ref{lemma:exporderstrict}, since $2^{p^+} \in \F$ and $p < p^+$} \nonumber \\ 
2^p \in \F     && \mbox{\qquad by Lemma~\ref{lemma:finiteexp}, since it is inhabited}\label{eq:h28}\\
p^+ < 2^p && \mbox{\qquad by the induction hypothesis and (\ref{eq:h19}) and (\ref{eq:h28})}\nonumber \\
p^+ + p^+ <  2^p + 2^p  && \mbox{\qquad by Lemma~\ref{lemma:addorder2}}\nonumber \\ 
p^+ + p^+ < 2^{p^+} && \mbox{\qquad by (\ref{eq:h20})} \nonumber \\  
p^{++} + p  < 2^{p^+} && \mbox{\qquad by the law $x + y^+ = x^+ + y$}\nonumber \\
p^{++} \le p^{++} + p   && \mbox{\qquad  by Lemma~\ref{lemma:xlessthan_xplusy}} \nonumber \\
p^{++} \mbox{\ is inhabited}  && \mbox{\qquad  by the definitions of $\le$ and addition} \nonumber\\ 
p^{++} \in \F  && \mbox{\qquad by Lemma~\ref{lemma:successorF}}  \label{eq:h50}\nonumber \\
p^{++} + p  \in \F  && \mbox{\qquad by the definition of $\le$ and Lemma~\ref{lemma:inhabited_sum}}\nonumber\\
p^{++} < 2^{p^+}   && \mbox{\qquad by Lemma~\ref{lemma:le_transitive3}}\nonumber 
\end{eqnarray}
But that is (\ref{eq:4262}), the desired goal.  That completes the 
induction step.  \end{proof}

\begin{lemma} \label{lemma:lessthansum}
Let $q \in \F$.  Then for all $n \in \F$ and $p \in \F$,
$$ n = p+q \imp \zero < q \imp p < n.$$
\end{lemma}

\noindent{\em Remark}.  This lemma links addition and order.
It probably can be proved directly from the definitions of addition and order,
but here we prove it by induction. Nevertheless we do have to use the definition
of addition directly at one of the steps.
\smallskip

\noindent\begin{proof}  By induction on $q$, which is legal since the formula is stratified.
The formula to be proved includes the quantifiers on $n$ and $p$.
\smallskip

{\em Base case}.  There is nothing to prove because of the hypothesis $q \neq \zero$. 
\smallskip

{\em Induction step}. Suppose $n = p+q^+$ and 
$\zero < q^+$.  
As usual in induction proofs,
we also assume $q^+$ is inhabited. 
 Then 
\begin{eqnarray*}
n = p^+ + q  && \mbox{\qquad by Lemma~\ref{lemma:addition2}}\\ 
q < \zero \ \lor\ q = \zero \ \lor \ \zero < q  && \mbox{\qquad by Theorem~\ref{theorem:finitetrichotomy}}
\end{eqnarray*}

{\em Case~1}, $q < \zero$ is impossible, by Lemma~\ref{lemma:xnotlessthanzero}.
\smallskip

{\em Case~2}, $q = \zero$. Then $q^+ = \one$ so $n = p+\one = p^+$.  Then $p < n$ by 
Lemma~\ref{lemma:xlessthansuccessorx}.
\smallskip

{\em Case~3}, $\zero < q$.  Then  
\begin{eqnarray*}
n \mbox{\ is inhabited}  && \mbox{\qquad by Corollary~\ref{lemma:cardinalsinhabited}}\\ 
p^+ \mbox{\ is inhabited} && \mbox{\qquad by the definition of addition, since $n \ p^+ q$}\\
p^+ < n    && \mbox{\qquad by the induction hypothesis, since $n = p^+ + q$}\\
p < p^+    && \mbox{\qquad by Lemma~\ref{lemma:xlessthansuccessorx}}\\
p < n       && \mbox{\qquad by transitivity} 
\end{eqnarray*} 
That completes the induction step.  
\end{proof} 

\begin{lemma} \label{lemma:propersmaller}  
Let $X$ and $Y$ be finite sets with $X \subseteq Y$ and $Y-X \neq \emptyset$.
Then $\Nc{X} < \Nc{Y}$.
\end{lemma}

\noindent{\em Remarks}.  (1) Of course this is not true without the finiteness
hypotheses. (2) The lemma does not mention addition, but the proof uses it;
hence its placement in the section on addition.
\smallskip

\noindent{\em Proof}.  \begin{eqnarray*}
Y = X \cup (Y-X)         && \mbox{\qquad by Lemma~\ref{lemma:finiteseparable}}\\  
Y-X \in \Nc{Y-X}          && \mbox{\qquad by Lemma~\ref{lemma:xinNcx}}\\  
\zero = \{\emptyset\}    && \mbox{\qquad by definition of $\zero$} \\
\zero \neq \Nc{Y-X}          && \mbox{\qquad since $Y-X \neq \emptyset$} \\
Y-X \in \FINITE              && \mbox{\qquad by Lemma~\ref{lemma:finitedif}} \\ 
\Nc{Y-X} \in \F               && \mbox{\qquad by Lemma~\ref{lemma:finitecardinals3}}\\ 
\neg\,(\Nc{Y-X} < \zero)      && \mbox{\qquad by Lemma~\ref{lemma:xnotlessthanzero}}\\  
\zero < \Nc{Y-X}           && \mbox{\qquad by Theorem~\ref{theorem:finitetrichotomy}}\\  
X \cap (Y-X) = \emptyset    && \mbox{\qquad by the definitions of $-$ and $\cap$} \\
\Nc{Y} = \Nc{X} + \Nc{Y-X}   && \mbox{\qquad by Lemma~\ref{lemma:cardinality_additive}, since $Y = X \cup (Y-X)$} \\  
\Nc{X} < \Nc{Y}              && \mbox{\qquad by Lemma~\ref{lemma:lessthansum}}  
\end{eqnarray*}                  
 That completes the proof of the lemma.

\section{Definition of multiplication} 
Specker did not make any use of multiplication.  If one could manage
to prove that $\F$ is infinite,  one would need multiplication to 
interpret HA in \INF.  But without knowing that $\F$ is infinite,
there are technical difficulties with multiplication.
 Some care is required to make sure that
the equations for multiplication work without assuming $\F$ is finite;
the equations must have the property that if one side is in $\F$, so 
is the other side.  That is, if one side ``overflows'',  so does the other
side.  To arrange this, we must first ensure that
addition has the same property.  This ultimately goes back to the theorem
that successor never takes the value zero, not just on an integer argument
but on any argument whatever.  We carried out those details (and they can 
still be found in earlier versions of this paper on ArXiv), but we have
not included them here.

\section{\texorpdfstring{Results about $\T$ }{Results about T }}
Here we constructivize Specker's \sect 5.  

\begin{definition}\label{definition:T} 
\begin{eqnarray*}
 \T(\kappa) &=& \T \kappa \ = \  \{ u:  \exists x\, (x \in \kappa \ \land \ u \sim \USC(x))\} 
\end{eqnarray*}
\end{definition} 
The formula is stratified, giving $x$ index 0, $u$ and 
$\kappa$ index 1.  
We will use $\T(\kappa)$ only when $\kappa$ is a finite cardinal, although 
that is not required by the definition.   Note that $\T(\kappa)$ has one
type higher than $\kappa$.
Thus we cannot define the graph of $\T$ or the graph of $\T$ restricted to
$\FregeN$. 

\begin{lemma}\label{lemma:Tmembers} If $\kappa \in \F$, then 
$$ x \in \kappa \iff \USC(x) \in \T \kappa.$$
\end{lemma}  

\noindent\begin{proof}  {\em Left to right}.
\begin{eqnarray*}
 x \in \kappa &&\mbox{\qquad by hypothesis} \\
 \USC(x) \sim \USC(x) &&\mbox{\qquad by Lemma~\ref{lemma:sim}}\\
 \USC(x)\in \T(\kappa) &&\mbox{\qquad by Definition~\ref{definition:T}}
\end{eqnarray*}
That completes the left-to-right direction.
\smallskip

{\em Right to left}.   
\begin{eqnarray*}
\USC(x) \in \T \kappa  &&\mbox{\qquad by hypothesis} \\
\exists z\,(z \in \kappa \ \land \ \USC(z) \sim \USC(x)) && \mbox{\qquad by definition of $\T$}\\
z \sim x  &&\mbox{\qquad by Lemma~\ref{lemma:uscsimilar}}\\
x \in \kappa && \mbox{\qquad by   Lemma~\ref{lemma:finitecardinals0}}
\end{eqnarray*}
 That completes the right-to-left direction.  \end{proof}

\begin{lemma} \label{lemma:T} 
If $\kappa \in \FregeN$  then for every  $x \in \kappa$, 
$\T(\kappa) = \Nc{\USC(x)}$.  
\end{lemma} 

\noindent\begin{proof}  
Suppose $\kappa \in \FregeN$.  Then $\kappa$ is inhabited, by 
Corollary~\ref{lemma:cardinalsinhabited}.  Let $x \in \kappa$.  Then
\begin{eqnarray*}
 x \in \FINITE && \mbox{\qquad by Lemma~\ref{lemma:finitecardinals1}}\\
 \USC(x) \in \FINITE  && \mbox{\qquad by Lemma~\ref{lemma:uscfinite}} \\
 \Nc{\USC(x)} \in \FregeN  && \mbox{\qquad by Lemma~\ref{lemma:finitecardinals3}} \\
 \USC(x) \in \T(\kappa)   && \mbox{\qquad by Lemma~\ref{lemma:Tmembers}}\\
 \USC(x) \sim \USC(x)     && \mbox{\qquad by Lemma~\ref{lemma:sim}}\\
 \USC(x) \in \Nc{\USC(x)}  && \mbox{\qquad by Definition~\ref{definition:Nc}}
\end{eqnarray*}
We remark that we cannot finish the proof at this point by 
Lemma~\ref{lemma:cardinalsdisjoint}, because we do not yet know $\T (\kappa) \in \FregeN$.  Instead:  by extensionality it suffices to prove
\begin{eqnarray}
\forall u\, (u \in \T (\kappa) \iff u \in \Nc{\USC(x)})  \label{eq:2620}
\end{eqnarray}
{\em Left to right}. Suppose $u \in \T(\kappa)$.  By definition of $\T$,
there exists $w \in \kappa$ with $u \sim \USC(w)$.  Then 
\begin{eqnarray*}
w \sim x  && \mbox{\qquad by Lemma~\ref{lemma:finitecardinals2},
since $w \in \kappa$ and $x \in \kappa$}\\
\USC(w) \sim \USC(x) && \mbox{\qquad  by Lemma~\ref{lemma:uscsimilar}}\\
u \sim \USC(x) && \mbox{\qquad by Lemma~\ref{lemma:sim} (transitivity of $\sim$), since $u \sim \USC(w)$}
\end{eqnarray*}
That completes the proof of the right-to-left direction
of (\ref{eq:2620}).
\smallskip

{\em Right to left}.   Suppose $u \in \Nc{\USC(x)}$.  Then $u \sim \USC(x)$.
Since $x \in \kappa$, we have $u \in \T(\kappa)$ by the definition of $\T$.
\end{proof}

 \begin{lemma}\label{lemma:Ncdef} If  $\kappa \in \FregeN$ and 
 $x \in \kappa$ then
 $\kappa = \Nc{x}$.
 \end{lemma}
 
 \noindent\begin{proof}  Let $\kappa \in \FregeN$ and $x \in \kappa$.
 By extensionality, it suffices to prove that for all $u$,
 $$u \in \kappa \iff u \in \Nc{x}.$$
 {\em Left to right}. Suppose $u \in \kappa$.  Then 
 \begin{eqnarray*}
 u \sim x && \mbox{\qquad by Lemma~\ref{lemma:finitecardinals2}}\\
 x \sim u && \mbox{\qquad by Lemma~\ref{lemma:sim}} \\
 u \in \Nc{x} &&\mbox{\qquad by Definition~\ref{definition:Nc}} 
 \end{eqnarray*}
\smallskip

{\em Right to left}.  Suppose $u \in \Nc{x}$.  Then
\begin{eqnarray*}
u \sim x &&\mbox{\qquad by Definition~\ref{definition:Nc}}\\
u \in \kappa && \mbox{\qquad by Lemma~\ref{lemma:finitecardinals0}}
\end{eqnarray*}
\end{proof}

\begin{lemma} \label{lemma:SpeckerT} If $\Nc{x} \in \FregeN$, then  $\T(\Nc{x}) = \Nc{\USC(x)}$.
\end{lemma}

\noindent\begin{proof}  By Lemma~\ref{lemma:T}, with $\kappa= \Nc{x}$.
\end{proof}

\begin{lemma} \label{lemma:Tfinite}
If $m \in \FregeN$ then $\T m  \in \FregeN$.
\end{lemma}

\noindent{\em Remark}.  Since the graph of $\T$ is not definable, we 
cannot express the lemma as $\T : \FregeN \to \FregeN$.
\medskip

\noindent\begin{proof} Let $m \in \FregeN$.  By Corollary~\ref{lemma:cardinalsinhabited},
$m$ is inhabited.  Let $a \in m$. Then 
\begin{eqnarray*}
 \USC(a) \in \T m   &&\mbox{\qquad by Lemma~\ref{lemma:Tmembers}} \\
a \in \FINITE      &&\mbox{\qquad by Lemma~\ref{lemma:finitecardinals1}}\\
\USC(a)\in \FINITE  &&\mbox{\qquad by Lemma~\ref{lemma:uscfinite}}\\
\Nc{\USC(a)} \in \FregeN && \mbox{\qquad by  Lemma~\ref{lemma:finitecardinals3}}\\
 \T m  \in \FregeN && \mbox{\qquad by Lemma~\ref{lemma:T}}
\end{eqnarray*}
\end{proof}

\begin{lemma} \label{lemma:Nc_unitclass}
 Every singleton has cardinal $\one$.  That is, 
$\forall x\, (\Nc{\{x\}} = \one)$.
\end{lemma}

\noindent\begin{proof}  By definition, $\one = \zero^+$ and $\zero = \{\emptyset\}$.
Then the members of $\one$ are sets of the form $\emptyset \cup \{ r\}$, by 
the definition of successor.  But $\emptyset \cup \{r\} = \{r\}$. 
Hence the members of $\one$ are exactly the unit classes.
Let $x$ be given; then by definition of $\Nc{\{x\}}$,  $\Nc{\{x\}}$ 
contains exactly the sets similar to $\{x\}$.
By Lemma~\ref{lemma:singletons_similar}, that is exactly the unit classes.
Hence $\Nc{\{x\}}$ and $\one$ have the same members, namely all unit classes.
By extensionality, $\Nc{\{x\}} = \one$.  \end{proof}

\begin{lemma}\label{lemma:Tsuccessor} For all $m \in \FregeN$ with an inhabited
successor, we have 
$$ \T (m^+) = (\T m )^+.$$
\end{lemma}

\noindent\begin{proof} Since $m^+$ is inhabited,  there is an $x \in m$ 
and $a \not\in x$ (so $x \cup \{a\} \in m^+$). 
Then 
\begin{eqnarray*}
m^+ \in \FregeN && \mbox{\hskip1.58in by Lemma~\ref{lemma:successorF}}\\ 
\T(m^+) &=& \Nc{\USC(x \cup \{a\})}\mbox{\qquad\quad by Lemma~\ref{lemma:T}}\\
 &=&  \Nc{\USC(x) \cup \{\{a\}\}} \mbox{\quad\ \  by Lemma~\ref{lemma:usc_successor}} \\
  &=&  (\Nc{\USC(x)})^+   \mbox{ \qquad\qquad by Lemma~\ref{lemma:Ncsuccessor}} \\
  &=& (\T m)^+    \mbox{\hskip1.17in by Lemma~\ref{lemma:T}} 
\end{eqnarray*}
\end{proof}

\begin{lemma}[Specker 5.2] \label{lemma:Tzero}
 $\T(\zero) = \zero$.
\end{lemma}  

\noindent\begin{proof}    We have
$\USC(\emptyset) = \emptyset$ as there are no singleton subsets of $\emptyset$.
Since $\zero = \Nc{\emptyset}$, by Lemma~\ref{lemma:SpeckerT} we have
$\T(\zero) = \Nc{\USC(\emptyset)} = \Nc{\emptyset} = \zero$.  
\end{proof}

\begin{lemma}[Specker 5.2] \label{lemma:Tone}
$\T(\one) = \one$.
\end{lemma}  

\noindent\begin{proof}
\begin{eqnarray*}
\{\emptyset\} \in \one && \mbox{\qquad by definition of $\one$}\\
\T(\one) = \Nc{\USC(\{\emptyset\}} && \mbox{\qquad  by Lemma~\ref{lemma:T}}\\
\T(\one) = \Nc{\{\{\emptyset\}\}} && \mbox{\qquad since $\USC(\{\emptyset\})=\{\{\emptyset\}\}$}\\
\Nc{\{\{\emptyset\}\}} = \one   && \mbox{\qquad by Lemma~\ref{lemma:Nc_unitclass}}\\
\T(\one) = \one  &&   \mbox{\qquad by the two previous lines} 
\end{eqnarray*}
\end{proof}

\begin{lemma}[Specker 5.2] \label{lemma:Ttwo}
$\T(\two) = \two$.
\end{lemma} 

\noindent\begin{proof}  We have
\begin{eqnarray*}
 \T(\two) &=& \T (\one^+) \mbox{\qquad\ \ since $\two = \one^+$} \\
 &=&  ( \T(\one))^+ \mbox{\qquad by Lemma~\ref{lemma:Tsuccessor}}\\
 &=& \one^+ \mbox{\qquad\qquad\ by Lemma~\ref{lemma:Tone}}  \\
 &=& \two.
\end{eqnarray*}
\end{proof}

\begin{lemma}[Specker~5.5]\label{lemma:Torder} Let $m,n \in \F$.  Then 
$$ n < m \imp \T n < \T m.$$
\end{lemma}

\noindent{\em Remarks}.  Specker~5.5 asserts that for cardinal 
numbers $p$ and $q$ we have $p \le q \iff \T p \le \T q$.   
Specker does not prove a version of that lemma with strict inequality.
\medskip

\noindent\begin{proof}  The formula in the lemma is stratified, 
with the relation $<$ occurring as a parameter.  Therefore we can 
prove by induction that for $n \in \F$,
$$ \forall m\in \F\,(n < m \imp \T n < \T m).$$
{\em Base case}, $n = \zero$. Suppose $\zero < m$; we must show $\T \zero < \T m$.  Since $\T \zero = \zero$, we have to show $\zero < \T m$. By 
Theorem~\ref{theorem:finitetrichotomy}, we have 
$$\T m < \zero \ \lor\ \T m = \zero \ \lor \ \zero < \T m $$
and only one of the three disjuncts holds.  Therefore it suffices to 
rule out the first two disjuncts, as the third is the desired conclusion.
By Lemma~\ref{lemma:nothinglessthanzero}, the first one is impossible.
We turn to the second. 
 Suppose $\T m = \zero$.  Since $m \in \F$, by Corollary~\ref{lemma:cardinalsinhabited}
 we have $a \in m$ for some $a$.  Then $\USC(a) \in \T m$, by definition of $\T$.
Since $\T m = \zero$, we have $\USC(a) \in \zero$.  Since $\zero = \{\emptyset\}$,
we have $\USC(a) = \emptyset$.  Then $a = \emptyset$.  Since $\emptyset \in \zero$,
by Lemma~\ref{lemma:cardinalsdisjoint} and the fact that $a \in m$,
 we have $m = \zero$.  But that contradicts the assumption $\zero < m$,
 by Lemma~\ref{lemma:nothinglessthanzero}.
   That completes the base case.
 \smallskip
 
{\em Induction step}.  Suppose $n^+ < m$ and $n^+$ is inhabited.
  We must show $\T (n^+) < \T m$. We have  
\begin{eqnarray*}
m \neq \zero && \mbox{\qquad since $n^+ < m$ and nothing is less than zero} \\
m = r^+  &&\mbox{\qquad for some $r\in \F$, by Lemma~\ref{lemma:nonzeroissuccessor}}\\
n^+ < r^+ && \mbox{\qquad since $n^+ < m$ and $m = r^+$}\\
 n < r  &&\mbox{\qquad by Lemma~\ref{lemma:successorstrict}}\\
\T n < \T r &&\mbox{\qquad by the induction hypothesis} \\
\exists a\, (a \in m) &&\mbox{\qquad by Corollary~\ref{lemma:cardinalsinhabited}}\\
\exists a\, (a \in r^+) &&\mbox{\qquad since $m = r^+$}\\
(\T r)^+ = \T(r^+)  &&\mbox{\qquad by Lemma~\ref{lemma:Tsuccessor}}\\
\exists u\,(u \in n^+) &&\mbox{\qquad by the definition of $\le$, since $n^+ < r^+$}\\
(\T n)^+ = \T(n^+)  &&\mbox{\qquad by Lemma~\ref{lemma:Tsuccessor}, since $n^+$ is inhabited}\\
\T(r^+) \in \F   && \mbox{\qquad by Lemma~\ref{lemma:Tfinite}}\\
\T(n^+) \in \F   && \mbox{\qquad by Lemma~\ref{lemma:Tfinite}}\\
\exists u\,(u \in \T(r^+)) &&\mbox{\qquad by Corollary~\ref{lemma:cardinalsinhabited}}\\
\exists u\,(u \in \T(n^+)) &&\mbox{\qquad by Corollary~\ref{lemma:cardinalsinhabited}}\\
\exists u\, u \in (\T n)^+ &&\mbox{\qquad since $(\T n)^+ = \T(n^+)$ }\\
\exists u\, u \in (\T r)^+ &&\mbox{\qquad since $(\T r)^+ = \T(r^+)$ }\\
 (\T n)^+ < (\T r)^+  &&\mbox{\qquad by Lemma~\ref{lemma:successorstrict}}\\
 \T(n^+) < \T(r^+) &&\mbox {\qquad since $(\T n)^+ = \T(n^+)$ and $(\T r)^+ = \T(r^+)$}\\
 \T(n^+) < \T m  &&\mbox{\qquad since $r^+ = m$}
\end{eqnarray*}
 That completes the
induction step.  \end{proof}

\begin{lemma} [Specker 5.3] \label{lemma:Tsum} Let $m,n\in \FregeN$ and 
suppose $n+m \in \FregeN$.  Then
$$\T(n+m) = \T n  + \T m .$$
\end{lemma}

\noindent{\em Remark}. This theorem can be proved directly from the 
definitions involved, but  we need it only for finite cardinals, and
it is simpler to prove it by induction.
\medskip

\noindent\begin{proof} 
By induction on $m$ we prove 
\begin{eqnarray}
 \forall n \in \FregeN\,( n+m \in \FregeN \imp \T(n+m) = \T n  + \T m ). \label{eq:3060}
\end{eqnarray}
The formula is stratified, since $\T$ raises indices by one.
\smallskip

{\em Base case}, $m = \zero$.  We have to prove $\T(n+\zero) = \T n + \T(\zero)$. Since
$\T(\zero) = \zero$ by Lemma~\ref{lemma:Tzero}, and $n + \zero = n$ by 
Lemma~\ref{lemma:addition2}, that reduces to $\T n  = \T n $.  That completes
the base case.
\smallskip

{\em Induction step}.   The induction hypothesis is (\ref{eq:3060}).  We suppose
that $m^+$ is inhabited and that $n+m^+ \in \FregeN$.  
We must prove $\T(n+m^+) = \T n + \T(m^+)$.
In order to apply the induction hypothesis, we need $n+m \in \FregeN$.
Since $n + m^+ \in \FregeN$, it is inhabited, by Corollary~\ref{lemma:cardinalsinhabited}. By Lemma~\ref{lemma:addition2},
$(n+m)^+$ is inhabited. Hence it has a member, which must be of the form
$x \cup \{a\}$ where $x \in n+m$.  Thus $n+m$ is inhabited.  Then 
by Lemma~\ref{lemma:inhabited_sum}, $n+m \in \FregeN$.  Therefore, by 
the induction hypothesis (\ref{eq:3060}), we have
\begin{eqnarray*}
   \T(n+m) = \T n  + \T m .  
\end{eqnarray*}
Taking the successor of both sides, we have
\begin{eqnarray*}
(\T(n+m))^+ &=& (\T n  + \T m )^+ \\
\T((n+m)^+) &=& (\T n  + \T m )^+  \mbox{\qquad by Lemma~\ref{lemma:Tsuccessor}}\\
\T(n+ m^+) &=& (\T n  + \T m )^+  \mbox{\qquad by Lemma~\ref{lemma:addition2}}\\
&=& \T n  + (\T m)^+ \mbox{\qquad\ \ by Lemma~\ref{lemma:addition2}}\\
&=& \T n  + \T (m^+)  \mbox{\qquad\ \ by Lemma~\ref{lemma:Tsuccessor}}
\end{eqnarray*}
That is the desired goal.  That completes the induction step.
\end{proof}

\begin{lemma}[Specker 5.8] \label{lemma:expTinhabited}
For   $m \in \F$, $2^{\T m }$ is inhabited.
\end{lemma}

\noindent\begin{proof} Let $m \in \F$.  Then
\begin{eqnarray*}
u \in m &&\mbox{\qquad for some $u$, by Corollary~\ref{lemma:cardinalsinhabited}}\\
\USC(u) \in \T m &&\mbox{\qquad by Definition~\ref{definition:T}}\\
\SSC(u) \in 2^{\T m} &&\mbox{\qquad by the definition of exponentiation}
\end{eqnarray*}
\end{proof}

\begin{lemma} \label{lemma:expTinF} For $m \in \F$, $2^{\T m} \in \F$.
\end{lemma}

\noindent\begin{proof}  Suppose $m \in \F$.  Then
$\exists x\, (x \in 2^{\T m})$, by Lemma~\ref{lemma:expTinhabited}.
Then by the definition of exponentiation, for some $u$
we have 
$$\SSC(u) \in 2^{\T m} \ \land \ \USC(u) \in \T m.$$
Then
\begin{eqnarray*}
 \T m \in \F  &&\mbox{\qquad  by Lemma~\ref{lemma:Tfinite}}\\
 \USC(u) \in \FINITE && \mbox{\qquad by Lemma~\ref{lemma:finitecardinals1}}\\
 u \in \FINITE && \mbox{\qquad by Lemma~\ref{lemma:uscfinite}}\\
 \USC(u) \in \T m && \mbox{\qquad by definition of $\T$} \\
 \SSC(u) \in \FINITE && \mbox{\qquad by Lemma~\ref{lemma:finitepowerset}}\\
 \SSC(u) \in 2^{\T m} &&\mbox{\qquad by Lemma~\ref{lemma:exp_members2}}\\
 2^{ \T m} \in \F  && \mbox{\qquad by Lemma~\ref{lemma:finiteexp}}
\end{eqnarray*}
\end{proof}

\begin{lemma} \label{lemma:successorT}
Suppose $m \in \F$.  Then $(\T m)^+ \in \F$.
\end{lemma} 

\noindent\begin{proof}  Suppose $m \in \F$.  Then
\begin{eqnarray*}
2^{\T m} \in \F && \mbox{\qquad by Lemma~\ref{lemma:expTinF}}\\
\T m \in \F   && \mbox{\qquad by Lemma~\ref{lemma:Tfinite}}\\
\exists u\,( u \in 2^{\T m}) && \mbox{\qquad by Corollary~\ref{lemma:cardinalsinhabited}}\\
\T m < 2^{\T m}  && \mbox{\qquad by Lemma~\ref{lemma:mlessthanexpm}}\\
(\T m)^+ \in \F && \mbox{\qquad by Lemma~\ref{lemma:successorbounded}}
\end{eqnarray*}
\end{proof} 

\begin{lemma}[Specker 5.9] \label{lemma:expT}
For  $m \in \F$, if $2^m$ is inhabited, then  $2^{\T m } = \T(2^m)$.
\end{lemma}

\noindent\begin{proof}  Suppose $2^m$ is inhabited.  Then 
there exists $a$ with $\USC(a) \in m$.  Then 

\begin{eqnarray*}
2^m =\Nc{\SSC (a)}  && \mbox{\qquad by Lemma~\ref{lemma:expuscssc}}\\
\SSC(a) \in 2^m       && \mbox{\qquad by Lemma~\ref{lemma:exp_members2}}\\
\USC (\SSC (a)) \in \T (2^m)  && \mbox{\qquad by Lemma~\ref{lemma:Tmembers}}\\
\USC (\USC (a)) \in \T m   && \mbox{\qquad by Lemma~\ref{lemma:Tmembers}}\\
\T m \in \FregeN  && \mbox{\qquad by Lemma~\ref{lemma:Tfinite}}\\
\SSC (\USC(a)) \in 2^{T m }  && \mbox{\qquad by Lemma~\ref{lemma:exp_members2}}\\
2^m \in \FregeN  && \mbox{\qquad by Lemma~\ref{lemma:finiteexp}}\\
2^{\T m } \in \FregeN  && \mbox{\qquad by Lemma~\ref{lemma:finiteexp}}\\
\Nc{\SSC (\USC (a))} = \Nc{\USC (\SSC (a))}   && \mbox{\qquad by Lemma~\ref{lemma:sscusc}}\\
2^{\T m } = \Nc{\SSC (\USC (a))}  && \mbox{\qquad by Lemma~\ref{lemma:Ncdef}}\\
\T (2^m) = \Nc{\USC (\SSC a)}  && \mbox{\qquad by Lemma~\ref{lemma:T}}\\
2^{\T m } = \T (2^m)  && \mbox{\qquad from the last three equations}
\end{eqnarray*}
\end{proof} 

\begin{lemma}\label{lemma:Toneone} 
For $n,m \in \FregeN$, we have 
$$\T n  = \T m \imp n = m $$   
\end{lemma}

\noindent\begin{proof} Suppose $\T n  = \T m $.  By Corollary~\ref{lemma:cardinalsinhabited}, we can find $a \in n$
and $b \in m$.  Then 
\begin{eqnarray*}
\USC(a) \in \T n  &&\mbox{\qquad by definition of $\T$}\\
 \USC(b) \in \T m   &&\mbox{\qquad by definition of $\T$}\\  
 \T n  = \T m  &&\mbox{\qquad by hypothesis}\\
 \USC(a) \in \T n &&\mbox{\qquad by the previous two lines}\\
 \T m \in \F && \mbox{\qquad by Lemma~\ref{lemma:Tfinite}}\\
 \USC(a) \sim \USC(b) &&\mbox{\qquad by Lemma~\ref{lemma:finitecardinals2}}\\
 a \sim b  &&\mbox{\qquad by Lemma~\ref{lemma:uscsimilar}}\\
 b \in n && \mbox{\qquad by Lemma~\ref{lemma:finitecardinals0}}\\
 n = m   &&\mbox{\qquad by Lemma~\ref{lemma:cardinalsdisjoint}}
 \end{eqnarray*}
\end{proof}

\begin{lemma}[Converse to Specker~5.3] \label{lemma:fivepointthree_converse}
Let $a,b,c \in \F$.  Then 
$$  \T a + \T b \in \F \imp \T a + \T b = \T c \imp a+b = c.$$
\end{lemma}

\noindent{\em Remark}. It is not assumed that $a+b \in \F$. 
Indeed, that follows from the stated conclusion.   
\medskip

\noindent\begin{proof}  The formula is stratified, giving $a$, $b$, and $c$
all index zero.  Therefore we may proceed by induction on $b$.
\smallskip

{\em Base case}\,: 
We have
\begin{eqnarray*}
\T a + \T \zero = \T c  && \mbox{\qquad by assumption}\\
\T a + \zero = \T c && \mbox{\qquad by Lemma~\ref{lemma:Tzero}}\\
\T a = \T c && \mbox{\qquad by Lemma~\ref{lemma:addition2}}\\
a = c && \mbox{\qquad by Lemma~\ref{lemma:Toneone}}\\
a + \zero = c && \mbox{\qquad by Lemma~\ref{lemma:addition2}}
\end{eqnarray*}
That completes the base case.
\smallskip

{\em Induction step}\,: We have
\begin{eqnarray*}
\T a + \T (b^+) = \T c  && \mbox{\qquad by assumption}\\
\exists u\, (u \in b^+)  && \mbox{\qquad by assumption}\\
b^+ \in \F   && \mbox{\qquad by Lemma~\ref{lemma:successorF}}\\
\T (b^+) = (\T b)^+  && \mbox{\qquad by Lemma~\ref{lemma:Tsuccessor}}\\
\T a + (\T b)^+ = \T c && \mbox{\qquad by the preceding lines}\\
(\T a + \T b)^+ = \T c && \mbox{\qquad by Lemma~\ref{lemma:addition2}}\\
c \neq \zero  && \mbox{\qquad by Lemmas~\ref{lemma:Tzero} and~\ref{lemma:Fregesuccessoromits0}}\\
c = r^+     && \mbox{\qquad for some $r$, by Lemma~\ref{lemma:nonzeroissuccessor}}\\
(\T a + \T b)^+ = \T (r^+)  && \mbox{\qquad by the preceding two lines}\\
 (\T a + \T b)^+ = (\T  r)^+ && \mbox{\qquad by Lemma~\ref{lemma:Tsuccessor}}\\
 \T a + \T (b^+) \in \F  && \mbox{\qquad by assumption} \\
 (\T a + \T b)^+ \in \F  && \mbox{\qquad by Lemmas~\ref{lemma:Tsuccessor} and~\ref{lemma:addition2}}\\
 \exists u\, (u \in (\T a + \T b)^+) && \mbox{\qquad by Corollary~\ref{lemma:cardinalsinhabited}}\\
  \exists u\, (u \in (\T a + \T b)) && \mbox{\qquad by definition of successor}\\
 \exists u\, (u \in (\T r)^+) && \mbox{\qquad by Corollary~\ref{lemma:cardinalsinhabited}}\\
\T r \in \F && \mbox{\qquad by Lemma~\ref{lemma:Tfinite}}\\
 \T a \in \F  &&\mbox{\qquad by Lemma~\ref{lemma:Tfinite}}\\
 \T b \in \F  &&\mbox{\qquad by Lemma~\ref{lemma:Tfinite}}\\
\T a + \T b \in \F && \mbox{\qquad by Lemma~\ref{lemma:inhabited_sum}}\\
 \T a + \T b = \T r && \mbox{\qquad by Lemma~\ref{lemma:successoroneone}}\\
 a + b = r  && \mbox{\qquad by the induction hypothesis} \\
 (a + b)^+ = r^+  && \mbox{\qquad by the preceding line}\\
  a + b^+ = r^+  && \mbox{\qquad by Lemma~\ref{lemma:addition2}}\\
 a+b^+ = c  && \mbox{\qquad since $r^+ = c$} 
 \end{eqnarray*}
That completes the induction step.
\end{proof}

\begin{lemma}\label{lemma:Tlessthan} 
For $n,m \in \FregeN$, we have 
$$n < m \iff \T n  < \T m.$$
\end{lemma}

\noindent\begin{proof}  Left to right is Lemma~\ref{lemma:Torder}.
\smallskip

Right to left. Suppose $\T n < \T m$.  By Theorem~\ref{theorem:finitetrichotomy},
we have $n < m$ or $n = n$ or $m < n$.  We argue 
by cases.  
\smallskip

Case 1, $n < m$.  Then  we are done, since that is the desired 
conclusion.  
\smallskip

Case 2, $n = m$ then $\T n = \T m$.  By Lemma~\ref{lemma:Tfinite},
$\T n \in \F$ and $\T m \in \F$, so by 
Theorem~\ref{theorem:finitetrichotomy},
$\T n = \T m$ contradicts $\T n < \T m$.  That completes Case~2.
\smallskip

Case~3, $m < n$.  Then $\T m < \T n$ by Lemma~\ref{lemma:Torder}.
\end{proof}

\begin{lemma}[Specker~5.6] \label{lemma:Tonto} Suppose $p,q \in \F$ and $p < \T q$.
Then there exists $r \in \F$ such that $p = \T r$.
\end{lemma}

\noindent\begin{proof}  By induction on $p$ we will prove 
\begin{eqnarray}
\forall q \in \F\, (p < \T q \imp \exists r \in \F\, (p = \T r)) \label{eq:4877}
\end{eqnarray}
The formula is stratified, giving $q$ and $r$ index 0 and $p$ index 1,
so induction is legal.
\smallskip

{\em Base case}, $p = 0$. Then $r = \zero$ satisfies $p = \T r$, by Lemma~\ref{lemma:Tzero}.  That completes the base case.
\smallskip

{\em Induction step}.  The induction hypothesis is (\ref{eq:4877}).
Suppose $p^+ < \T q$ and $p^+$ is inhabited.  Then  
\begin{eqnarray*}
p < p^+ &&\mbox{\qquad  by Lemma~\ref{lemma:lessthansuccessor}}\\
p < \T q &&\mbox{\qquad by Lemma~\ref{lemma:lessthan_transitive}}\\
p = \T r && \mbox{\qquad for some $r$, by (\ref{eq:4877})}
\end{eqnarray*}
Now I say that $r^+$ is inhabited.  To prove that:
\begin{eqnarray*}
\T r = p < p^+ < \T q &&\mbox{\qquad as already proved} \\
\T r < \T q && \mbox{\qquad from the previous line} \\
r < q  && \mbox{\qquad by Lemma~\ref{lemma:Tlessthan}}\\
r^+ \le q && \mbox{\qquad by Lemma~\ref{lemma:noinsertions}} \\
\exists u\,(u \in r^+) &&\mbox{\qquad by the definition of $\le$} 
\end{eqnarray*}
That completes the proof that $r^+$ is inhabited.
Then since $p = \T r$, we have
\begin{eqnarray*}
p^+ = (\T r)^+ = \T(r^+) && \mbox{\qquad by Lemma~\ref{lemma:Tsuccessor}}
\end{eqnarray*}
That completes the induction step.  \end{proof}

\begin{lemma} \label{lemma:Tinexp} Suppose $p \in \F$ and $2^p$ is 
inhabited.  Then $p = \T\,q$ for some $q \in \F$.
\end{lemma}

\noindent\begin{proof}  Suppose $p \in \F$ and $2^p$ is inhabited. 
Then by definition of exponentiation, for some $a$ we have $\USC(a) \in p$
and $\SSC(a) \in 2^p$.  By definition of $\T$ we have $p = \T (\Nc{a})$.
By Lemma~\ref{lemma:finitecardinals3}, we have $\Nc{a} \in \F$.  
\end{proof}

\begin{lemma}[Specker~5.5]\label{lemma:Tlessthanorequal} 
For $n,m \in \FregeN$, we have 
$$n \le m \iff \T n  \le \T m.$$
\end{lemma}

\noindent{\em Remark}. It is also possible to prove this lemma
directly from the definitions of $\le$ and $\T$, instead of 
from Lemma~\ref{lemma:Tlessthan} as we do here,  and then 
prove Lemma~\ref{lemma:Tlessthan} from this lemma.  Or 
we could prove this lemma by induction as we did Lemma~\ref{lemma:Tlessthan}.
\smallskip

\noindent\begin{proof} We have
\begin{eqnarray*}
n \le m \iff n < m \ \lor \ n = m && \mbox{\qquad by Lemma~\ref{lemma:letolessthan}}\\
\T n  \le \T m \iff \T n < \T m \ \lor \ \T n = \T m && \mbox{\qquad by Lemma~\ref{lemma:letolessthan}}
\end{eqnarray*}
Now to prove the desired conclusion:
\smallskip

{\em Left to right}\,: if $n < m$ then $\T n < \T m$ by Lemma~\ref{lemma:Tlessthan}, so $\T n \le \T m$.
And if $n = m$, then $\T n = \T m \le \T m$, by Lemma~\ref{lemma:le_reflexive}.
\smallskip

{\em Right to left}\,: if $\T n < \T m$ then $n < m$ by Lemma~\ref{lemma:Tlessthan}, so $n \le m$.
And if $\T n = \T m$, then $n = m$ by Lemma~\ref{lemma:Toneone}.
\end{proof}

\begin{lemma}\label{lemma:epluse} Let $e \in \F$ and $e + e \in \F$.
Then $e^+ \in \F$.
\end{lemma}

\noindent\begin{proof}  By Theorem~\ref{lemma:FregeNdecidable}, $e = \zero \ \lor \ e \neq \zero$.  If $e=\zero$ then $e^+ = \one$, so we are done by 
Lemma~\ref{lemma:oneF}.  Therefore we may assume $e \neq \zero$.
By Corollary~\ref{lemma:cardinalsinhabited}, $e+e$ is inhabited. 
By the definition of $<$,  there exist $x$ and $y$ with $x \in e$ and $y \in e$ and $x \cap y = \emptyset$.
Then 
\begin{eqnarray*}
y \neq \emptyset && \mbox{\qquad since if $y = \emptyset$ then $e = \zero$, by Lemma~\ref{lemma:cardinalsdisjoint}}\\
y \in \FINITE  && \mbox{\qquad by Lemma~\ref{lemma:finitecardinals1}}\\   
a \in y  && \mbox{\qquad for some $a$, by Lemma~\ref{lemma:empty_or_inhabited}}\\ 
a \not \in x  && \mbox{\qquad since $x \cap y = \emptyset$}\\  
x \cup \{a\} \in e^+  && \mbox{\qquad by definition of successor} \\   
e^+ \in \F   && \mbox{\qquad by Lemma~\ref{lemma:successorF}}   
\end{eqnarray*}
\end{proof}

\begin{lemma} \label{lemma:Teven}  If $\T c$ is even, then     
$c$ is even.  More precisely, if $c, a \in \F$ and $\T c = a + a$ and $a+a\in \F$,
then there exists $b \in \F$ with $c = b+b$.
\end{lemma}

\noindent\begin{proof}  The formula is stratified, giving $b$ and $c$
index 0 and $a$ index 1.  $\F$ is just a parameter, so it does not 
need an index.  Therefore we can proceed by induction on $a$.
\smallskip

{\em Base case}\,: Suppose $\T c   =\zero + \zero$ and $c \in \F$.  We have    
\begin{eqnarray*}
\zero+\zero = \zero  && \mbox{\qquad since $x+\zero = x$}\\
\T (\zero) = \zero    && \mbox{\qquad by Lemma~\ref{lemma:Tzero}}\\
\T (c) = \zero        && \mbox{\qquad since $\T(c) = \zero + \zero = \zero$}\\
c = \zero            && \mbox{\qquad by Lemma~\ref{lemma:Toneone}} \\
\exists b\, (c = b + b) &&\mbox{\qquad namely, $b = \zero$}
\end{eqnarray*}

{\em Induction step}\,:  Suppose $a^+$ is inhabited and $a^+\in \F$ and $\T c = a^+ + a^+$,
and $a \in \F$ and $a^+ + a^+ \in \F$.  (The assumption $a^+ \in \F$ is part of the induction hypothesis,
while the assumptions $a \in \F$ and $a^+$ is inhabited come with every proof by induction on $\F$.)
Then
\begin{eqnarray*}
a^+ + a^+ = (a + a)^{++} && \mbox{\qquad by Lemma~\ref{lemma:addition2}}\\  
(a + a)^{++} \in \F   && \mbox{\qquad since $a^+ + a^+  \in \F$}  
\end{eqnarray*}
I say that 
\begin{eqnarray}
a + a \in \F  && \label{eq:5019}   
\end{eqnarray}
It is surprisingly difficult to prove that.  I had to go back to the definition of addition.
Since $a^+ + a^+ \in \F$, there exists $x \in a^+ + a^+$, by
Corollary~\ref{lemma:cardinalsinhabited}.  By the definition of addition,  $x$ has the form
\begin{eqnarray*}
x = u \cup v     &&  \mbox{\qquad with $u\cap v = \emptyset$ and $u \in a^+ $ and $v \in a^+$}\\
u = z \cup \{p\} \land v = w \cup\{q\}  && \mbox{\qquad with $z \in a$ and $w \in a$, by definition of successor}\\
z \cup w \in a + a  &&  \mbox{\qquad by the definition of addition}\\
a + a \in \F        && \mbox{\qquad by Lemma~\ref{lemma:inhabited_sum}}
\end{eqnarray*}
That completes the proof of (\ref{eq:5019}).  
Similarly, $x \cup u \in a^+ + a$, so 
\begin{eqnarray*}
a^+ \in \F      && \mbox{\qquad by Lemma~\ref{lemma:successorF}, since $a\in F$ and $a^+$ is inhabited}\\   
a^+ + a \in \F  && \mbox{\qquad by Lemma~\ref{lemma:inhabited_sum}}\\
(a+a)^+ = a^+ + a  && \mbox{\qquad by Lemma~\ref{lemma:addition2}}\\
(a+a)^+ \in \F     && \mbox{\qquad by the preceding lines} 
\end{eqnarray*}
 Continuing, we have
\begin{eqnarray*}
\T\, c = (a + a)^{++} && \mbox{\qquad by Lemma~\ref{lemma:addition2}}\\    
\T\, c \neq \zero && \mbox{\qquad by Lemma~\ref{lemma:Fregesuccessoromits0}}\\   
\one \neq (a+a)^{++} && \mbox{\qquad by Lemma~\ref{lemma:successoroneone}}\\
\T\, c \neq \one  && \mbox{\qquad since $T\, c = (a+a)^{++}$} \\  
c \neq \zero  &&\mbox {\qquad by Lemma~\ref{lemma:Tzero}}\\  
c = r^+ && \mbox{\qquad for some $r\in \F$, by Lemma~\ref{lemma:nonzeroissuccessor}} \\ 
r  \neq \zero&& \mbox{\qquad since if $r = \zero$ then $r^+ = c = \one$, so $\T\,c = \one$}\\ 
r = t^+ && \mbox{\qquad for some $t\in \F$, by Lemma~\ref{lemma:nonzeroissuccessor}}\\
c = t^{++} && \mbox{\qquad by the preceding lines}\\
\T\, c = (\T\, t)^{++} && \mbox{\qquad by Lemma~\ref{lemma:Tsuccessor}}\\
(a+a)^{++} = (\T\, t)^{++} &&\mbox{\qquad since $\T c = (a+a)^{++}$}\\
\T\, t \in \F  && \mbox{\qquad by Lemma~\ref{lemma:Tfinite}}\\
\T\, t= a+a &&\mbox{\qquad by Lemma~\ref{lemma:successoroneone}}\\   
t = e+e && \mbox{\qquad for some $e \in \F$, by the induction hypothesis}\\
t^{++} = (e^+ + e^+) && \mbox{\qquad by Lemma~\ref{lemma:addition2}}\\
c = b + b && \mbox{\qquad with $b = e^+$, by the preceding lines} \\
e^+ \in \F && \mbox{\qquad by Lemma~\ref{lemma:epluse}, since $e+e = t \in \F$}
\end{eqnarray*}
That completes the induction step.  \end{proof}

\begin{lemma}[Specker 5.4] \label{lemma:fivepointfour} Let $m \in \NC$.  Then
$$ m \neq \T (m) + \one $$   
\end{lemma}

\noindent{\em Remark}.  And so on, with $\one$ replaced by $\two$ or  23, 457, and any number you could name.
If $\T m \neq m$,  $\T m$ must be a non-standard distance away from $m$.
\smallskip

\noindent\begin{proof} We give the proof for $\one$.
Recall that $m$ is even if $m = p+p$ for some $p \in \F$, and odd
if $m = p + p + \one$ for some $p \in \F$.  Then if $m$ is even, $m+ \one$ is odd, and vice
versa.  One can verify by induction that every integer is either even or odd, and not both. 
Suppose $\m = \T \m + \one$.  If $m$ is even, then $\T m$ is even, by Lemma~\ref{lemma:Teven},
so $\T (m) + \one$ is odd, contradiction.  If $m$ is odd, then $m = k^+ = k + \one$ for some $k \in \F$, 
since $\zero$ is even. Then $k$ is even.  Then $\T(m) = \T (k^+) = \T(k + \one) = \T(k) + \one$, which
is odd since $\T(k)$ is even.  Then $\T(m) + \one$ is even, contradiction, since $m$ is odd and equal to $\T(m) + \one$.
\end{proof}

\begin{lemma}\label{lemma:adds_to_zero} For all $p,q$, if $ p+q = \zero$ then $ p = \zero$.
\end{lemma}

\noindent{\em Remark}. No additional hypothesis is needed.
\medskip

\noindent\begin{proof}  By definition, $\zero = \{ \emptyset\}$.  
By the definition of addition, there exist $a$ and $b$ with $a \in p$ and 
$b \in q$ and $a \cap b = \emptyset$, such that $a \cup b \in \zero$.
Then $a \cup b = \emptyset$.  It follows that $a = \emptyset$ and $b = \emptyset$.
On the other hand, if $a$ or $b$ had a non-empty member, then by the definition
of addition, $a+b$ would have a non-empty member, so $\zero$ would have a 
non-empty member.  Therefore $p = q = \{\emptyset\} = \zero$.  
\end{proof}

\begin{lemma} \label{lemma:dividebytwo} For $x,y \in \F$, 
$ x + x = y + y \imp x = y$.
\end{lemma}

\noindent\begin{proof}  The formula is stratified;  we prove it by 
induction on $x$, in the form 
$$ \forall y\in F\, (x+x = y+y \imp x = y).$$

{\em Base case}\,: Suppose $\zero+\zero = y+y$.  Then $\zero = y+y$.  By 
Lemma~\ref{lemma:adds_to_zero}, $y = \zero$.  That completes the base case.
\smallskip

{\em Induction step}\,: Suppose $x^+ + x^+ = y + y$, and suppose (as always in 
induction proofs) that $x^+$ is inhabited.  Then 
\begin{eqnarray*}
x^+ + x^+ = (x + x)^{++} &&\mbox{\qquad by Lemma~\ref{lemma:addition2}}\\
y \neq \zero && \mbox{\qquad by Lemma~\ref{lemma:Fregesuccessoromits0}}\\
y = r^+ && \mbox{\qquad for some $r$, by Lemma~\ref{lemma:nonzeroissuccessor}}\\
(x+x)^{++} = (r+r)^{++}  &&\mbox{\qquad by Lemma~\ref{lemma:addition2}}\\
x+x = r+r   &&\mbox{\qquad by Lemma~\ref{lemma:successoroneone}}\\
x = r &&\mbox{\qquad by the induction hypothesis} \\
x^+ = r^+ && \mbox{\qquad by the preceding line} \\
x^+ = y && \mbox{\qquad since $y = r^+$} 
\end{eqnarray*}
That completes the induction step.
\end{proof} 

\begin{lemma} \label{lemma:expandT}  Let $p \in \F$.  Then 
$$   2^p \in \F  \iff \exists q \in \F \,(p = \T q).$$
\end{lemma}

\noindent\begin{proof}  Suppose $p \in \F$.  Left to right:
\begin{eqnarray*}
2^p \in \F  &&\mbox{\qquad assumption}\\
\exists u\, (u \in 2^p) && \mbox{\qquad by Corollary~\ref{lemma:cardinalsinhabited}}\\
\USC(a)\in p && \mbox{\qquad for some $a$, by the definition of exponentiation}\\
\USC(a) \in \FINITE && \mbox{\qquad by Lemma~\ref{lemma:finitecardinals1}} \\
a \in \FINITE && \mbox{\qquad by Lemma~\ref{lemma:uscfinite}}
\end{eqnarray*}   
\begin{eqnarray*}
\Nc{a} \in \F && \mbox{\qquad by Lemma~\ref{lemma:finitecardinals3}}\\
a \in \Nc{a} && \mbox{\qquad by Lemma~\ref{lemma:xinNcx}}\\
\USC(a) \in \T (\Nc{a}) && \mbox{\qquad by definition of $\T$} \\
\T (\Nc{a}) \in \F && \mbox{\qquad by Lemma~\ref{lemma:Tfinite}}\\
\USC(a) \in p \cap \T(\Nc{a}) && \mbox{\qquad by definition of intersection}\\
p = \T (\Nc{a}) && \mbox{\qquad by Lemma~\ref{lemma:cardinalsdisjoint}}\\
\exists q \in \F\, (p = \T a) && \mbox{\qquad namely $q = \Nc{a}$}
\end{eqnarray*}
That completes the proof of the left-to-right direction.
\smallskip

{\em Right to left}.  Suppose $p = \T q$ and $q \in \F$.  Then
\begin{eqnarray*}
u \in q && \mbox{\qquad for some $u$, by Corollary~\ref{lemma:cardinalsinhabited}}\\
\USC(u) \in \T q && \mbox{\qquad by definition of $\T$ } \\
\T q \in \F  && \mbox{\qquad by Lemma~\ref{lemma:Tfinite}}\\
\SSC(u) \in 2^{\T q}  && \mbox{\qquad by definition of exponentiation}\\
\SSC(u) \in 2^p && \mbox{\qquad since $p = \T q$}\\
2^p \in \F && \mbox{\qquad by Lemma~\ref{lemma:finiteexp}}
\end{eqnarray*}
That completes the proof of the right-to-left direction.  
\end{proof} 
\smallskip

\begin{definition} \label{definition:imageT}  Let $X$ be any set of cardinals.   
Then we define
$$\T``(X) = \{ \T(u): u \in X\}.$$
\end{definition}
or more explicitly
$$ \T``(X) = \{ \T(u): u \in X\} = \{ y: \exists u\in X\, (y = \T u)\}.$$
The formula in the definition is stratified, giving $u$ index 0 and $y$ and $X$ index 1.
Actually, $X$ is just a parameter and does not even need an index.  Therefore the 
definition is legal in \INF.  We note that it is not a function definable in \INF.  It
is just an abbreviation for a comprehension term.   Note also that the set $X$ can
be finite or not, and the cardinals in $X$ can be finite or not.

In general images commute with union.  For images under $\T$ we have
\begin{lemma} \label{lemma:Timageunion}  $\T``(X \cup Y) = \T``(X) \cup \T``(Y)$.
\end{lemma}

\noindent\begin{proof} This is proved in a few short steps from the definitions
of $T``(X)$ and $\cup$.
\end{proof}

\begin{lemma} \label{lemma:NCsum} Let $a$ and $b$ be finite disjoint sets.  Then 
$$ \Nc{a \cup b} = \Nc{a} + \Nc{b}.$$
\end{lemma}

\noindent\begin{proof}  {\em Left to right}. Suppose $t \in \Nc{a \cup b}$.  Then
\begin{eqnarray*}
a \cup b \in \Nc{a \cup b} && \mbox{\qquad by Lemma~\ref{lemma:xinNcx}}\\  
a \in \Nc{a} && \mbox{\qquad by Lemma~\ref{lemma:xinNcx}}\\  
b \in \Nc{b} && \mbox{\qquad by Lemma~\ref{lemma:xinNcx}}\\  
a \cap b = \emptyset && \mbox{\qquad by hypothesis}\\  
a \cup b \in \FINITE  && \mbox{\qquad by Lemma~\ref{lemma:union}}\\   
\Nc{a} \in \F          && \mbox{\qquad by Lemma~\ref{lemma:finitecardinals3}} \\  
\Nc{b} \in \F           && \mbox{\qquad by Lemma~\ref{lemma:finitecardinals3}} \\  
\Nc{ a \cup b} \in \F   && \mbox{\qquad by Lemma~\ref{lemma:finitecardinals3}} \\  
\Nc{a} + \Nc{b} \in \F   && \mbox{\qquad by Lemma~\ref{lemma:inhabited_sum}} \\ 
a \cup b \in \Nc{a} + \Nc{b} && \mbox{\qquad by the definition of addition}\\  
\Nc{a \cup b} =  \Nc{a} + \Nc{b} && \mbox{\qquad by Lemma~\ref{lemma:cardinalsdisjoint}}
\end{eqnarray*}
\end{proof}

\begin{lemma} \label{lemma:Timage}      
Let $X$ be a finite set of cardinal numbers.   Then 
$$ \Nc{\T``(X)} = \T (\Nc{X}).$$
\end{lemma}

\noindent\begin{proof}  
The displayed formula in the lemma
is stratified, giving $X$ index 1; then $\Nc{X}$ gets index 2 and $\T \Nc{X}$ gets index 3.
On the left, the members of $\T``(X)$ are $\T u$ for $u \in X$, so $u$ gets index 0, and 
$\T u$ gets index 1, so $\T``X$ gets index 2, so $\Nc{\T``(X)}$ gets index 3, the same 
as the right side of the equation.  So it is stratified, as claimed.

The  part of the lemma involving $X$ is 
$$ \forall X\, (X \in \FINITE \imp X \subset \NC \imp \Nc{\T``(X)} = \T (\Nc{X})),$$
and this is also stratified, since $\FINITE$  and $\NC$ are just parameters.
Therefore we can prove it by induction on finite sets.
\smallskip

{\em Base case}, $X = \emptyset$.  On the right, $\Nc{\emptyset} = \zero$, so $\T(\Nc{\emptyset}) = \T (\zero) = \zero$.  On 
the left, $\T``(\emptyset) = \emptyset$, so $\Nc{\T``(\emptyset)} = \zero$.  That completes the base case.   
\smallskip

{\em Induction step}. Suppose $X$ is finite and $c \not\in X$. We have to show
\begin{eqnarray*}
 \Nc{{\T``(X \cup \{c\}}} = \T (\Nc{X \cup \{c\}}).   
\end{eqnarray*}
We have
\begin{eqnarray}
X \cup \{c\} \in \FINITE  && \mbox{by Lemma~\ref{lemma:finite_adjoin}, since $c \not\in X$} \nonumber \\ 
\T``(X \cup \{c\}) = \T``(X) \cup \{\T(c)\}  &&  \mbox{by Lemma~\ref{lemma:Timageunion}}\label{eq:h30} \\    
\Nc{\T``(X \cup \{c\})} = \Nc{\T``(X) \cup \{\T(c)\}} &&  \mbox{ by the preceding line} \nonumber \\    
\T c \not\in \T``(X)   && \mbox{by Lemma~\ref{lemma:Toneone}, since $c \not\in X$} \label{eq:h31}\\    
\Nc{\T``(X)} = \T (\Nc{X})  && \mbox{ by the induction hypothesis} \label{eq:h21}
\end{eqnarray}
\begin{eqnarray*}
\Nc{X} \in \F           && \mbox{\qquad by Lemma~\ref{lemma:finitecardinals3}, since $X \in \FINITE$} \\ 
\T (\Nc{X}) \in \F         && \mbox{\qquad by Lemma~\ref{lemma:Tfinite}}\\  
\T``(X) \in \FINITE    && \mbox{\qquad by Lemma~\ref{lemma:finitecardinals1}} \\ 
\{ \T(c)\} \in \FINITE   && \mbox{\qquad by Lemma~\ref{lemma:singletons_finite}} \\ 
\T``(X) \cap \{T(c)\} = \emptyset &&\mbox{\qquad by (\ref{eq:h31}) }   
\end{eqnarray*}
\begin{eqnarray*}
\Nc{\T``(X) \cup \{\T(c)\})} = \Nc{\T``(X)} + \Nc{\{\T(c)\}} &&  \mbox{\qquad by Lemma~\ref{lemma:NCsum}}\\ 
\Nc{\T``(X) \cup \{\T(c)\})} = \Nc{\T``(X)} + \one   && \mbox{\qquad by Lemma~\ref{lemma:Nc_unitclass}}  
\end{eqnarray*}
\begin{eqnarray*}
\Nc{\T``(X) \cup \{\T(c)\})} = \T (\Nc{X}) + \one && \mbox{\qquad by the induction hypothesis (\ref{eq:h21})}\\  
\Nc{\T``(X) \cup \{\T(c)\})} = \T (\Nc{X}) + \T(\one) && \mbox{\qquad since $\T(\one) = \one$}\\ 
\Nc{X} + \one = \Nc{X \cup \{c\}}    && \mbox{\qquad by Lemma~\ref{lemma:NCsum} since $c \not \in X$} \\
\Nc{X} + \one \in \F                                  && \mbox{\qquad by Lemma~\ref{lemma:inhabited_sum}}\\ 
\Nc{{\T``(X) \cup \{\T(c)\})}} = \T (\Nc{X} + \one) && \mbox{\qquad  by Lemma~\ref{lemma:Tsum}} \\
\Nc{{\T``(X) \cup \{\T(c)\}}}= \T (\Nc{X \cup \{c\}}) && \mbox{\qquad by the preceding lines} \\
\Nc{{\T``(X \cup \{c\})}} = \T (\Nc{X \cup \{c\}}) && \mbox{\qquad by (\ref{eq:h30})} 
\end{eqnarray*}
That completes the induction step.  \end{proof}

\begin{lemma} \label{lemma:Timagefinite}  Let $X$ be a finite set of cardinals.  Then $\T``(X)$ is finite.
\end{lemma}  

\noindent\begin{proof}  Let $X$ be a finite set of cardinals.  Then 
\begin{eqnarray}
X \subseteq \FINITE  && \mbox{\qquad by hypothesis} \nonumber\\ 
\T``(X) \in \FINITE && \mbox{\qquad by hypothesis} \nonumber\\
\Nc{\T``(X)} = \T (\Nc{X})  && \mbox{\qquad by Lemma~\ref{lemma:Timage}} \label{eq:h300}\\ 
\Nc{X} \in \F       && \mbox{\qquad by Lemma~\ref{lemma:finitecardinals3}}\nonumber\\ 
\T (\Nc{X}) \in \F   && \mbox{\qquad by Lemma~\ref{lemma:Tfinite}} \label{eq:h5} \\   
\Nc{\T``(X)} \in \F   && \mbox{\qquad by (\ref{eq:h300}) and (\ref{eq:h5})}\nonumber\\  
\T``(X) \in \Nc{\T``(X)} && \mbox{\qquad by Lemma~\ref{lemma:xinNcx}}\nonumber \\   
\T``(X) \in \FINITE  && \mbox{\qquad by Lemma~\ref{lemma:finitecardinals1}}\nonumber 
\end{eqnarray}
\end{proof}

\section{Cartesian products}

The Cartesian product of two sets is defined as usual; the definition
is stratified, so it can be given in \INF.  But because ordered pairs
raise the types by two,  the cardinality of $A \times B$  is not 
the product of the cardinalities of $A$ and $B$, but instead  
it is the product of $\T^2$ of those cardinalities.  In this section
we provide a proof of this fact, in the interest of setting down the 
fundamental facts about the theory of finite sets. 

\begin{lemma} \label{lemma:productofunion}  Let $X$, $Y$, and $Z$ be finite sets.
Then $$(X \cup Y) \times Z = (X \times Z) \cup (Y \times Z).$$    
\end{lemma}

\noindent\begin{proof}  This follows in a few steps from extensionality,
the definition of $\times$, and the logical fact that 
$$(P \ \lor \ Q) \ \land \ R \iff (P \ \land \ R) \ \lor \ (Q \ \land \ R).$$
\end{proof}

\begin{lemma} \label{lemma:productfinite_helper2} Let $Y$ be   
a finite set and let $a$ be any set.  Then 
$\{a \} \times Y$ is finite.   If 
 $\kappa = \Nc{Y}$ then  $\T^2 \kappa = \Nc{\{a\} \times Y}$.
\end{lemma}

\noindent\begin{proof}  Consider the map $f: \USC^2(Y) \to \{a\}\times Y$
defined by 
\begin{eqnarray*}
f = \{ \langle \{\{ y\}\}, \langle a,y \rangle \rangle: y \in Y\}. 
\end{eqnarray*}
The formula is stratified, giving $y$ and $a$ index 0, so $\langle a,y\rangle$
gets index 2, as does $\{\{y\}\}$.  $Y$ gets index 1.  Since the 
formula is stratified,  $f$ can be defined in \INF.  

One then proves without any surprises that $f$ 
is a similarity from  $\USC^2(Y)$  to $\{a\} \times Y$.  We omit the straightforward
196-line verification of that fact.   
\smallskip

Then we have
\begin{eqnarray*}
\USC^2(Y) \sim \{a\} \times Y  && \mbox{\qquad since $f$ is a similarity}\\
\USC(Y) \in \FINITE && \mbox{\qquad by Lemma~\ref{lemma:uscfinite}}\\
\USC^2(Y) \in \FINITE && \mbox{\qquad by Lemma~\ref{lemma:uscfinite}}\\
\{a\} \times Y \in \FINITE && \mbox{\qquad by Lemma~\ref{lemma:similar_to_finite}}\\
\Nc{{\USC^2(Y)}} = \T^2 \kappa && \mbox{\qquad by definition of $\T$} \\
\Nc{{\{a\} \times Y}} = \T^2 \kappa && \mbox{\qquad by Lemma~\ref{lemma:finitecardinals0}}
\end{eqnarray*}
\end{proof}

\begin{lemma} \label{lemma:productfinite2}   Let $X$ and $Y$ be   
finite sets. Then $X \times Y$ is finite.
\end{lemma}

\noindent\begin{proof}    
The formula to be proved is
\begin{eqnarray*}
 X \in \FINITE   \imp  \forall Y \in \FINITE\, 
 X \times Y \in \FINITE 
\end{eqnarray*}
That formula (and the hypotheses listed before it) are
 stratified, giving $X$ and $Y$ index 1; then $X \times Y$
 gets index 3, $\Nc{X \times Y}$ gets index 4, $\kappa = \Nc{X}$ gets index 2, 
 and $\T^2(\kappa)$ gets index 4;  since multiplication is a function, the 
 whole left-hand side gets index $4$.  
 $\FINITE$ is just 
 parameter.  Therefore we may proceed by induction on finite sets $X$.
 \smallskip
 
 {\em Base case}.  We have to show $\emptyset \times Y \in \FINITE$.    
 One shows $\emptyset \times Y = \emptyset$ using the definition of $\times$,
 and then $\emptyset \in \FINITE$ by Lemma~\ref{lemma:lambda_finite}.
 \smallskip
 
 {\em Induction step}.  Assume $X$ is finite and $a \not\in X$.  The 
 induction hypothesis is 
 \begin{eqnarray}
 \forall Y \in \FINITE\,( X \times Y \in \FINITE) \label{eq:2870}
 \end{eqnarray}
 Assume $X \cup \{a\} \in \FINITE$.  We have to prove $(X \cup \{a\}) \times Y \in \FINITE$.
 We have
 \begin{eqnarray*}
 X \in \FINITE && \mbox{\qquad by hypothesis}\\ 
 X \times Y \in \FINITE && \mbox{\qquad by the induction hypothesis (\ref{eq:2870})} \\ 
 \{a\} \times Y \in \FINITE && \mbox{\qquad by Lemma~\ref{lemma:productfinite_helper2}}\\ 
 (X \cup \{a\}) \times Y = (X \times Y) \cup (\{a\} \times Y) && \mbox{\qquad by Lemma~\ref{lemma:productofunion}} \\ 
  (X \times Y) \cap (\{a\} \times Y) = \emptyset && \mbox{\qquad since $a \not\in X$}\\  
  (X \cup Y) \cup (\{a\} \times Y) \in \FINITE && \mbox{\qquad by Lemma~\ref{lemma:union}} \\ 
  (X \cup \{a\}) \times Y \in \FINITE  && \mbox{\qquad by the preceding lines}
 \end{eqnarray*}
 That completes the induction step.  \end{proof}
 
 \begin{lemma} \label{lemma:productfinite3}   Let $X$ and $Y$ be
finite sets. If $\kappa = \Nc{X}$ and $\mu = \Nc{Y}$,  then 
$$ \T^2( \kappa) \cdot  \T^2( \mu) = \Nc{X \times Y}.$$
\end{lemma}

\noindent{\em Remarks.}  Without  $\T^2$, the formula is 
not stratified.  It is not necessary to {\em assume} that $(\T^2 \kappa) 
 \cdot \T^2( \mu) \in \F$.  That will, of course, be a consequence, by Lemma~\ref{lemma:finitecardinals3}.   
\smallskip

\noindent\begin{proof}  By induction on finite sets, like Lemma~\ref{lemma:productfinite2}.  We omit
the proof, since we never use this lemma.  It is included only because it illustrates the general
situation that arises from using Kuratowski pairing, which increases the type.
\end{proof}

\section{Onto and one-to-one for maps between finite sets}
In this section, we prove the well-known theorems that for maps $f$ from 
a finite set $X$ to itself, $f$ is one-to-one if it is onto, and vice-versa.
These theorems are somewhat more difficult to prove constructively than classically,
but they are provable.  

 In treating this subject rigorously
one has to distinguish the relevant concepts precisely.  Namely, we have
\begin{eqnarray*}
f:X \to Y \\
Rel(f) \\
f \in \FUNC \\
oneone(f,X,Y) 
\end{eqnarray*}
$Rel(f)$ means that all the members of $f$ are ordered pairs.  $f \in \FUNC$ means that 
two ordered pairs in $f$ with the same first member have the same second member.  (Nothing
is said about possible members of $f$ that are not ordered pairs.)  $f:X \to Y$ means that 
if $x \in X$, there is a unique $y$ such that $\langle x,y \rangle \in f$ and that $y$ is in $Y$.
(But nothing is said about $\langle x, y\rangle \in f$ with $x \not \in X$.)  ``$f$ is one-to-one
from $X$ to $Y$'',  or $oneone(f,X,Y)$,  means $f:X \to Y$ and in addition, if $\langle x,y \rangle \in f$
and $\langle u,y\rangle \in f$ then $x = u$,  and if $y \in Y$ then $x \in X$.  (So $x=u$ does
not require $y \in Y$ or $x \in X$.)  In particular, $f:X \to Y$ does not require $dom(f) \subseteq X$,
so the identity function maps $X$ to $X$ for every $X$; but the identity function (on the universe)
has to be restricted to $X$ before it is one-to-one.  

\begin{definition}\label{definition:permutation} 
$f$ is a {\bf permutation} of a finite set $X$ if and only if $f:X \to X$,
and $Rel(f)$ and $f \in \FUNC$, and $dom(f) \subseteq X$,
and $f$ is both one-to-one and onto from $X$ to $X$.
\end{definition} 

In this section we will prove that for finite $X$, either one of the conditions ``one-to-one''  and ``onto''
implies the other, if all the other conditions are assumed.
\medskip

\noindent{\em Remark}. We do not need to specify $range(f) \subseteq X$, because that follows
from $dom(f) \subseteq X$ and $f:X\to X$.  The reader can check that none of the  conditions
in the definition are superfluous.
\medskip

\begin{lemma} \label{lemma:finitefunction} Let $A$ and $B$ be finite sets,
and let $f$ be a function with domain $A$, and  $f:A \to B$.  Then $f$ is finite.  
\end{lemma}

\noindent\begin{proof}  By induction on finite sets $A$ we prove that for all finite
sets $B$,  if the domain of $f$ is $A$ and $f: A \to B$, then $f$ is finite.   
\smallskip

{\em Base case}.  A function with domain $\emptyset$ is the empty function, which is finite.
\smallskip

{\em Induction step}.  Let $A$ and $B$ be finite sets, and let $c \not\in A$, and suppose
$f: A \cup \{c\} \to B$.   Then
\begin{eqnarray*}
\langle c, y \rangle \in f && \mbox{\qquad for some $y \in B$ } 
\end{eqnarray*}
Let $g:= f - \{\langle c,y \rangle \}$.  One can verify that $g:A \to B$
and the domain of $g$ is $A$.%
\footnote{Formalizing this sort of lemma makes one appreciate the informal 
functional notation; this lemma took 330 lines of Lean and several hours.
I changed ``One can easily verify''  to the present ``One can verify.''
} 
Then by the induction hypothesis, $g$ is finite.  Since $A$ and $B$ are finite,
equality on $A$ and $B$ is decidable, so any member of $f$ is either equal to 
$\langle c,y \rangle$ or not.  Therefore
$$ f = g \cup \{ \langle c,y \rangle\}.$$ 
Since $g$ is finite and $\{ \langle c,y \rangle\} \not\in f$, 
$f$ is also finite.  \end{proof}

\begin{lemma} [Decidable image] \label{lemma:decidable_image} Let $X$ and $Y$ be finite sets. Let
$f: X \to Y$  and suppose the domain of $f$ is  $X$.   Then the set $P$  defined by
$$ f(X)   = \{ y \in Y :  \exists x \in X\, \langle x,y \rangle \in f\} $$
is a decidable subset of $Y$.   
\end{lemma}

\noindent\begin{proof} Let $y \in X$.  Define
$$ Z:= \{ x\in X : \exists y \in Y\, (\langle x,y\rangle \in f) \}.$$
The formula is stratified, giving $x$ and $y$ index 0, $f$ index 3, and $X$ index 1.
Therefore the definition is legal.  
Then
\begin{eqnarray*}
f \subseteq X \times Y  && \mbox{\qquad since $dom (f) = X$}\\
f \in \FINITE  &&\mbox{\qquad by  Lemma~\ref{lemma:finitefunction}}\\
X \in \DECIDABLE && \mbox{\qquad by Lemma~\ref{lemma:finitedecidable}}\\
X \times Y \in \FINITE && \mbox{\qquad by Lemma~\ref{lemma:productfinite2}}\\
f \mbox{ \ is a separable relation on $X$} && \mbox{\qquad by Lemma~\ref{lemma:finiteseparable}}\\
Z \in \FINITE && \mbox{\qquad by Lemma~\ref{lemma:boundedquantification2}}\\
Z = \emptyset \ \lor \ \exists x\,(x \in Z) && \mbox{\qquad by  Lemma~\ref{lemma:empty_or_inhabited}}
\end{eqnarray*}
Putting in the definition of $Z$, we have the formula
in  the conclusion of the lemma.  \end{proof}

 \begin{theorem} \label{theorem:dedekind1}
Let $X$ be a finite set,  
and let $f: X \to X$ be a one-to-one 
function.  Then $f$ is onto.
\end{theorem}

\noindent
\begin{proof}  By induction on finite sets, we prove
that if $f:X \to X$ is one-to-one, then
 $f$ is onto.   By Lemma~\ref{lemma:finitedecidable}, 
$X$ has decidable equality.
\smallskip

{\em Base case}\,:
The only function defined on the empty set is the empty 
function, which is both one-to-one and onto. 
\smallskip

{\em Induction step}\,:  Let $X = B \cup \{a\}$, where 
$a \not \in B$, and $B$ is finite.  Suppose $f:X \to X$ is 
one-to-one.  We have to prove
\begin{eqnarray}
\forall y \in X\,\exists x \in X\, (\langle x,y \rangle \in f) \label{eq:1284}
\end{eqnarray}
By Lemma~\ref{lemma:decidable_image}, $a \in range(f) \ \lor \ a \not\in range(f)$.  Explicitly,
$$ \exists x \in X\, (\langle x,a \rangle \in f) \ \lor \
   \neg\, \exists x \in X\,  (\langle x,a \rangle \in f).$$
We argue by cases accordingly.
\smallskip

Case 1, $\exists x \in X\, (\langle x, a \rangle \in f)$.   Fix $c$ such that $c \in X$ and 
$ \langle c, a \rangle \in f$.   Since $X$ has decidable equality, we have $c = a \ \lor \ c \neq a$.
We argue by cases.
\smallskip

Case 1a, $c = a$.  Then $f:B \to B$.   Let $g$ be $f$ restricted to $B$.  Then $g$ is one-to-one,
since $f$ is one-to-one.  By the induction hypothesis, $g: B \to B$ is onto.  Now let $y \in X$.
Then $y = a \ \lor \ y \in B$.  If $y = a$, then $\langle a,a \rangle \in f$.  If $y \in B$,
then since $g$ is onto, there exists $x \in B$ with $\langle x,y\rangle \in B$.  Then
$\langle x, y \rangle \in f$.  That completes Case~1a.
\smallskip

Case 1b, $c \neq a$.  
 Since $f:X \to X$, there exists $b \in X $
such that $\langle a, b \rangle \in f$.  Then $a \neq b$, since $\langle c,a\rangle \in f$ 
and $\langle a,b \rangle \in f$, so if $a=b$ then $\langle a,a\rangle \in f$; then since
$f$ is one-to-one we have $a=c$, contradiction. 
Define
$$ g:=  (f - \{ \langle c, a\rangle\}- \{\langle a,b\rangle \} )\cup \{ \langle c, b\rangle\}.$$

{\em Remark}.  In case the formal use of ordered pairs is difficult for the reader accustomed
to functional notation, we  put the matter informally:  We have $a = f(c)$ and $b = f(a)$, and we make $g$ agree with $f$
except at $a$ and $c$, where we make $g(c) = b$. Thus having eliminated $a$ from both
domain and range, we will be able to show $g:B \to B$. The formal details follow.
\smallskip
 
We have $Rel(g)$, since by hypothesis $Rel(f)$.
  I say $dom(g) = B$.   By extensionality, it 
suffices to show 
\begin{eqnarray}
 \exists y\, (\langle t, y \rangle \in g)  \iff  t \in B \label{eq:6162}
\end{eqnarray}
{\em Left to right}.  Assume $\langle t, y \rangle \in g$.  
Then 
$$(\langle t, y\rangle \in f \ \land \ \langle t,y \rangle \neq \langle c,a \rangle \ 
\land\ \langle t,y \rangle \neq \langle a,b \rangle) \ \lor \ (t = c \ \land \ y = b).$$
  If the second 
disjunct holds, then $t = c$, and $c \in X$ but $c \neq a$, so $c \in B$; so $t \in B$.
Therefore we may assume the first disjunct holds:
$$(\langle t, y\rangle \in f \ \land \ \langle t,y \rangle \neq \langle c,a \rangle \ 
\land\ \langle t,y \rangle \neq \langle a,b \rangle).$$
Then $t \in X$ since $dom(f) = X$.   
  Since $\langle t,y \rangle \neq \langle a,b \rangle$, we have $y \neq b$.
Since $\langle a,b \rangle \in f$ and $\langle t,y \rangle \in f$ it follows that $t \neq a$.
Since $X = B \cup \{a\}$, we have $t \in B$.  That completes the left-to-right
direction of (\ref{eq:6162}).
\smallskip

{\em Right to left}. Suppose $t \in B$.  Since $dom f = X$ and $B \subseteq X$,  there exists $z$ such 
that $\langle t,z \rangle \in f$.  Unless $t = c$ or $t = a$,  we have $\langle t,z \rangle \in g$.
If $t = c$ we can take $y = b$.   Since $t \in B$ we do not have $t=a$.  That completes
the proof of (\ref{eq:6162}).  That completes the proof that $dom(g) = B$. 
\smallskip

Now I say that $g:B \to B$.  Suppose $x \in B$.  We must show there exists $y$
with $\langle x,y \rangle \in g$.  Since $f:X \to X$, there exists $y \in X$
such that $\langle x,y\rangle \in f$.  Then $x= c \ \lor \ x \neq c$.
 If $x \neq c$ then $\langle x, y \rangle \in g$.
 If $x = c$ then $\langle x, b\rangle \in g$.  That completes the proof that 
 $\exists y\, (\langle x,y \rangle \in g)$.  We must also show that 
 if $\langle x,y \rangle \in g$ and $\langle x,z \rangle \in g$ then $y = z$.
If $x \neq c$ then $\langle x,y \rangle \in f$ and $\langle x,z \rangle \in f$,
so $y = z$.  If $x = c$ then $y = b$ and $z = b$, so $y = z$.  That completes
the proof that $g:B \to B$. 
\smallskip

Now I say that $g$  is one-to-one.  Suppose $g(u) = g(v)$.  If $u \neq c$ and $v \neq c$,
then $g(u) = f(u)$ and $g(v) = f(v)$, so $u = v$ since $f$ is one-to-one.  If $u = c$ and $v \neq c$
then $g(u) = b$.  Since $v \neq c$, $g(v) = f(v) = b$.  Since $f$ is one-to-one, $v = a$. 
But $v \not\in B$,  so $\langle v, b \rangle \not\in g$, since $dom(g) = B$.  Similarly if 
$v = c$ and $u \neq c$.  That completes the proof that $g$ is one-to-one.
\smallskip

 By the induction hypothesis, $g$ is onto.
Now I say that $f$ is onto.  Let $y \in X$.  Then if $y = a$, we have $\langle c, y\rangle \in f$.
If $y=b$ we have $\langle a,y \rangle \in f$.  If $y \neq a$ and $y \neq b$, then $y = g(x) = f(x)$
for some $x$.  Since $X$ has decidable equality, these cases are exhaustive.  That completes
Case~1b.
\smallskip

Case 2, $\neg\,\exists x \in X\, (\langle x, a \rangle \in f)$.   
Let $g$ be $f$ restricted
to $B$.  Then $Rel(g)$, and $dom(g) = B$, and $g$ is one-to-one, and $g:B \to B$.  
Then by the induction hypothesis, $g$ is onto. Since $f:X \to X$, there exists
some $b \in X$ such that $\langle a, b \rangle \in f$.  By hypothesis $b \neq a$.
Then $b \in B$.  Since $g$ is onto, there exists $x \in B$ such that $\langle x,b\rangle \in g$.
Then $\langle x,b \rangle \in f$.  Since $f$ is one-to-one, we have $x = a$.  But $x \in B$,
while $a \not \in B$.  That contradiction completes Case~2. 
\end{proof}

\begin{lemma} \label{lemma:adjoin_cardinality}   
Let $B \in \FINITE$ and $a \not\in B$.
Then $\Nc{B \cup \{a\}} = (\Nc{B})^+$.
\end{lemma}

\noindent\begin{proof}  We have
\begin{eqnarray*}
B \in \Nc{B} && \mbox{\qquad by Lemma~\ref{lemma:xinNcx}}\\
B \cup \{a\} \in \Nc{B \cup \{a\}} && \mbox{\qquad by Lemma~\ref{lemma:xinNcx}}\\
B \cup \{a\} \in (\Nc{B})^+ && \mbox{\qquad by definition of successor}\\
B \cup \{a\} \in \FINITE && \mbox{\qquad by Lemma~\ref{lemma:finite_adjoin}}\\
\Nc{B \cup \{a\}} \in \F && \mbox{\qquad by Lemma~\ref{lemma:finitecardinals3}}\\
\Nc{B} \in \F && \mbox{\qquad by Lemma~\ref{lemma:finitecardinals3}}\\
(\Nc{B})^+ \in \F && \mbox{\qquad by Lemma~\ref{lemma:successorF}}\\
B \cup \{a\} \in \Nc{{B \cup \{a\}}} \cap (\Nc{B})^+ && \mbox{\qquad by the definition of intersection}\\
\Nc{{B \cup \{a\}}} = (\Nc{B})^+ && \mbox{\qquad by Lemma~\ref{lemma:cardinalsdisjoint}}
\end{eqnarray*}
\end{proof}

\begin{lemma} \label{lemma:nothingbetween}
Let $m,n \in \F$ and $m + n \le m^+$  and $m+n \in \F$ and $n \neq \zero$.  Then $n = \one$.
\end{lemma}

\noindent\begin{proof} 
\begin{eqnarray*}
 n = r^+ &&\mbox{\qquad for some $r \in \F$, by Lemma~\ref{lemma:nonzeroissuccessor}} \\
 m + r^+ \le m^+ && \mbox{\qquad since $m+n \le m^+$ and $n = r^+$}\\
a \in m + n \ \land \ b \in m^+ &&\mbox{\qquad for some $a$ and $b$, by definition of addition}\\
m^+ \in \F &&\mbox{\qquad by Lemma~\ref{lemma:successorF}}\\
m + r^+ +  k = m^+ &&\mbox{\qquad for some $k \in \F$, by Lemma~\ref{lemma:orderbyaddition}}\\
(m+r+k)^+ = m^+ &&\mbox{\qquad by Lemma~\ref{lemma:addition2}}\\
m+r\in \F && \mbox{\qquad by Lemma~\ref{lemma:subterms2}}\\
m+r+ k^+ = m^+ &&\mbox{\qquad by Lemma~\ref{lemma:addition2}}\\
m+r+k^+ \in \F && \mbox{\qquad since $m+r+k^+ = m^+ \in \F$}\\
m+r+k \in \F && \mbox{\qquad by Lemma~\ref{lemma:subterms3}}\\
m+r+k = m && \mbox{\qquad by Lemma~\ref{lemma:successoroneone}}\\
r+k+m = \zero + m && \mbox{\qquad by Lemma~\ref{lemma:addition2}}\\
m+r \in \F \ \land \ r+k \in \F &&\mbox{\qquad by Lemma~\ref{lemma:subterms}}\\
r+k+m \in \F &&\mbox{\qquad by commutativity and associativity, since $m+r+k \in \F$}\\
r+k = \zero && \mbox{\qquad by Lemma~\ref{lemma:subtraction}}\\
(m+r)^+ \le m^+ &&\mbox{\qquad by Lemma~\ref{lemma:addition2}}\\
m+r^+ \in \F  && \mbox{\qquad since $m+n \in \F$}\\
m+r \in \F && \mbox{\qquad by Lemma~\ref{lemma:subterms3}}\\
m+r = m   &&\mbox{\qquad by Lemma~\ref{lemma:successoroneone}}\\
m+r = m+\zero && \mbox{\qquad by Lemma~\ref{lemma:addition2}}\\
r = \zero && \mbox{\qquad by Lemma~\ref{lemma:subtraction}} \\
n = r^+ = \zero^+ = \one && \mbox{\qquad since $\one =\zero^+$} \\
r = \zero &&\mbox{\qquad by Lemma~\ref{lemma:adds_to_zero}}\\
r^+ = \one && \mbox{\qquad by the definition of $\one$}\\
n = \one && \mbox{\qquad since $n = r^+$}
\end{eqnarray*}
\end{proof}

\begin{lemma} \label{lemma:separableNc} Let $X \in \FINITE$
and let $Z$ be a separable subset of $X$.  Then 
$$ \Nc{Z} \le \Nc{X}.$$
\end{lemma}

\noindent\begin{proof} We have
\begin{eqnarray*}
\Nc{X} \in \F && \mbox{\qquad by Lemma~\ref{lemma:finitecardinals3}}\\
Z \in \FINITE && \mbox{\qquad by Lemma~\ref{lemma:separablefinite}}\\
\Nc{Z} \in \F && \mbox{\qquad by Lemma~\ref{lemma:finitecardinals3}}\\
X \in \Nc{X} && \mbox{\qquad by Lemma~\ref{lemma:xinNcx}}\\
Z \in \Nc{Z} && \mbox{\qquad by Lemma~\ref{lemma:xinNcx}}\\
\Nc{Z} \le \Nc{X} && \mbox{\qquad by the definition of $\le$}
\end{eqnarray*}
\end{proof}

\begin{theorem}\label{theorem:dedekind2}
Let $X$ be a finite set,  
and let $f: X \to X$ be onto, with $dom (f) \subseteq X$.
 Then $f$ is one-to-one.
\end{theorem}

\noindent\begin{proof}  We prove the more general fact that if $X$ and $Y$
are finite sets with $\Nc{X} \le \Nc{Y}$,  and $f:X \to Y$ is onto, then $f$ 
is one-to-one.  (The theorem follows by taking $Y=X$).  More explicitly,
we will prove by induction on finite sets $Y$ that 
\begin{eqnarray*}
&&\forall Y \in \FINITE\, \forall X \in \FINITE\, (\Nc{X} \le \Nc{Y} \imp 
\forall f\,(f \in \FUNC  \\
&&\imp Rel(f) \imp dom (f) \subseteq X   \\
&&\imp \forall x \in X\,\exists y \in Y\,(\langle x,y\rangle \in f) \\
&&\imp \forall y \in Y\, \exists x \in X\, (\langle x,y\rangle\in f) \\
&&\imp \forall y \in Y\,\forall x,z \in X\, (\langle x,y\rangle \in f \imp \langle z,y\rangle \in f 
\imp x = z)))
\end{eqnarray*}
The formula is stratified, giving $x,y,z$ index 0, $f$ index 3, $X$ and $Y$ index 1, and
$\Nc{X}$ and $\Nc{Y}$ index 2. $\FUNC$ and $\FINITE$ are parameters;
$Rel(f)$ is stratified giving $f$ index 3; $dom(f) \subseteq X$ can be expressed as 
$\forall x,y\,(\langle x,y \rangle \in f \imp x \in X)$, which is stratified.
 Therefore we may proceed by induction on finite sets $Y$.
\smallskip

{\em Base case}, $Y = \emptyset$.  Then (in the last line) $y \in Y$ is impossible, so 
the last line holds if the previous lines are assumed.  That completes the base case.
\smallskip

{\em Induction step}, $Y = B \cup \{a\}$ with $a \not\in B$ and $B \in \FINITE$.  Suppose $X \in \FINITE$,
and $f:X \to Y$
is onto, and $f \in \FUNC$ and $Rel(f)$ and $dom(f) \subseteq X$.  We must prove $f: X \to Y$ is one-to-one.
  Define
\begin{eqnarray}
 Z := \{x \in X : \langle x,a\rangle \in f\}.  \label{eq:6192}
 \end{eqnarray}
The formula is stratified, giving $x$ and $a$ index 0 and $f$ index 3,  so the definition is legal.
Since $f$ is onto,  $Z$ is inhabited.  I say that $Z$ is a separable subset of $X$.  That is,
\begin{eqnarray}
 \forall x \in X\, (\langle x,a \rangle \in f \ \lor \ \langle x,a \rangle \not \in f).\label{eq:6193}
 \end{eqnarray}
To prove that, let $x \in X$.  Since $f:X \to Y$, there exists $y \in Y$ with $\langle x,y \rangle \in f$.
Since $f \in \FUNC$, we have $\langle x, a \rangle \in f \iff y = a$.  Since $Y$ is finite,
we have $y = a \ \lor y \neq a$ by Lemma~\ref{lemma:finitedecidable}.  That completes
the proof of (\ref{eq:6193}).   Then by Lemma~\ref{lemma:separablefinite}, $Z \in \FINITE$ 
and $X - Z \in \FINITE$.  
\smallskip

Let $g$ be $f$ restricted to $X - Z$.  Then $g:X-Z \to B$ and $g$ is onto $B$.   I say that
\begin{eqnarray}
\Nc{X-Z} \neq \Nc{X} \label{eq:6392}
\end{eqnarray}
To prove that, assume $\Nc{X-Z} = \Nc{X}$.  Then
\begin{eqnarray*}
\Nc{X-Z} \in \F && \mbox{\qquad by Lemma~\ref{lemma:finitecardinals3}}\\
\Nc{X} \in \F && \mbox{\qquad by Lemma~\ref{lemma:finitecardinals3}}\\
X \sim X-Z  && \mbox{\qquad by Lemma~\ref{lemma:finitecardinals2}}\\
u \in Z  && \mbox{\qquad for some $u\in X$, since $f$ is onto $Y$}\\
X - Z \subseteq X  && \mbox{\qquad by the definition of $Z$}\\
X \neq X-Z   && \mbox{\qquad since $u \not\in X- Z$ but $u \in X$}
\end{eqnarray*}
Therefore $X$ is similar to a proper subset of $X$.  Then by 
Definition~\ref{definition:infinite}, $X$ is infinite.  Then 
by Theorem~\ref{theorem:infiniteimpliesnotfinite}, $X$ is not finite.
But that contradicts the hypothesis.  That completes the proof of 
(\ref{eq:6392}).

Now  I say 
that $\Nc{X-Z} \le \Nc{B}$.    To prove that:
\begin{eqnarray*}
\Nc{X-Z} \le \Nc{X} && \mbox{\qquad by Lemma~\ref{lemma:separableNc}}\\
\Nc{X-Z} < \Nc{X}  && \mbox{\qquad by (\ref{eq:6392}) and the definition of $<$}\\
\Nc{X} \le \Nc{B \cup \{ a\}} && \mbox{\qquad by hypothesis} \\
\Nc{B \cup \{a\}} = (\Nc{B})^+ && \mbox{\qquad since $a \not\in B$} \\
\Nc{X-Z} < \Nc{B}^+ && \mbox{\qquad by the previous two lines}\\
\Nc{X-Z} \le \Nc{B} && \mbox{\qquad by Lemma~\ref{lemma:nothingbetween}}
\end{eqnarray*}
Therefore we can apply the induction hypothesis to $g$.  Hence $g:X-Z \to B$ is one-to-one.
Therefore $g$ is a similarity.  Then
\begin{eqnarray*}
\Nc{X-Z} = \Nc{B}  &&\mbox{\qquad by Lemma~\ref{lemma:finitecardinals2} and ten omitted steps}\\
\Nc{X} = \Nc{X-Z} +  \Nc{Z} && \mbox{\qquad  by Lemma~\ref{lemma:cardinality_additive}}\\
\Nc{X} = \Nc{B} + \Nc{Z} && \mbox{\qquad by the previous two lines}\\
\Nc{X} \le \Nc{Y} && \mbox{\qquad by hypothesis} \\
 \Nc{B} + \Nc{Z} \le \Nc{Y} && \mbox{\qquad by the previous two lines}\\
\Nc{Y} = \Nc{B}^+ &&\mbox{\qquad since $Y = B \cup \{a\}$ and $a \not\in B$}\\
  \Nc{B} + \Nc{Z} \le \Nc{B}^+&& \mbox{\qquad by the previous two lines}\\
\Nc{Z} = \one   && \mbox{\qquad  by Lemma~\ref{lemma:nothingbetween}}
\end{eqnarray*}

By Lemma~\ref{lemma:one_members}, $Z$ is a unit class $\{c\}$ for some $c$.  By
(\ref{eq:6192}),   
$$\forall x\,( \langle x, a \rangle \in f \iff x = c).$$
I say that 
$f$ is one-to-one.  To prove that, let $u,v \in X$ and $\langle u,y\rangle \in f$
and $\langle v,y \rangle \in f$.  We must prove $u = v$.  Since $Y$
has decidable equality, we have $y = a \ \lor \ y \neq a$.  We argue by cases
accordingly.
\smallskip

Case 1, $y = a$.  Then $u \in Z$ and $v \in Z$.  Then $u = c$ and $v = c$, so $u = v$.
That completes Case~1.

Case 2, $y \neq a$. Then 
$u \not \in Z$ and $v \not\in Z$, so $\langle u,y\rangle \in g$ and $\langle v,y\rangle \in g$.
Since $g$ is one-to-one, we have $u = v$ as desired. That completes Case~2.  That completes
the induction step.  \end{proof}

\begin{theorem}\label{theorem:dedekind3}  Let $X$ and $Y$ be finite sets, and 
suppose $f: X \to Y$ is onto, and the domain of $f$ is $X$.   Then $\Nc{Y} \le \Nc{X}$. 
\end{theorem}

\noindent\begin{proof} By induction on finite sets $X$, we prove the theorem for all $Y$.
\smallskip

\noindent{\em Base case}.  If $f: \emptyset \to Y$ has domain $\emptyset$ and is onto $Y$
then $Y = \emptyset$, so $$\Nc{X} = \Nc{Y} = \Nc{\emptyset} = \zero.$$ 
\smallskip

\noindent{\em Induction step}.  Suppose $c \not \in X$ and $f$ has domain $X \cup \{c\}$,
and $f: X \cup \{c\} \to Y$ is onto.  Let $g$ be $f$ restricted to $X$, which is 
conveniently defined as $f \cap X \times Y$.  Then the domain of $g$ is exactly $X$. 

We have $f: X \cup \{c\} \to Y$, from which it follows in a few steps that also $g: X  \to Y$.
Then by Lemma~\ref{lemma:decidable_image}, 
 the image $g(X)$ of $X$ under $g$ is a decidable subset of $Y$. (That lemma requires
 that the domain of $g$ be exactly $X$, not larger,  which is why we had to use $g$ instead of $f$.)  
 That is,
 \begin{eqnarray*}
 (\exists x \in X\, g(x) = f(c)) \ \lor \ \neg\, \exists x \in X\, g(x) = f(c).   
\end{eqnarray*}
We argue by cases, as justified by that disjunction.
\smallskip

{\em Case~1}, $ \exists x \in X\, g(x) = f(c)$.  Then $g: X \to Y$ is onto.   
Then
\begin{eqnarray*}  
\Nc{Y} \le \Nc{X}  && \mbox{\qquad by the induction hypothesis} \\   
\Nc{X} <  \Nc{X}^+  && \mbox{\qquad by Lemma~\ref{lemma:lessthansuccessor}}\\  
\Nc{X}^+ = \Nc{X \cup \{c\}} && \mbox{\qquad by Lemma~\ref{lemma:Ncsuccessor} }\\
\Nc{Y} \le \Nc{X \cup \{c\}}  && \mbox{\qquad by the preceding lines}
\end{eqnarray*}
That completes Case~1.   
\smallskip

{\em Case~2}, $\neg\, \exists x \in X\, g(x) = f(c)$.  Let $t = f(c)$. 
 Then $g: X \to Y-\{t\}$ is onto.  
 We have 
 \begin{eqnarray*}
 Y - \{t\} \in \FINITE  && \mbox{\qquad by Lemma~\ref{lemma:oneout}} \\  
 g: X \to Y - \{t\} \mbox{\ and $g$ is onto } && \mbox{\qquad as one can check}  
 \end{eqnarray*}
 and the domain of $g$ is $X$.  
 Then by the induction
hypothesis, 
$$ \Nc{Y- \{t\}} \le \Nc{X}.$$    
We want to take the successor of both sides, but to do that we have to check that 
those successors are inhabited.
\begin{eqnarray*}
Y \mbox {\ has decidable equality} && \mbox{\qquad by Lemma~\ref{lemma:finitedecidable}}\\
(Y - \{t\}) \cup \{t\} = Y   && \mbox{\qquad by decidable equality on $Y$}\\ 
\exists u\, (u \in \Nc{X}^+)  && \mbox{\qquad namely $u = X \cup \{c\}$} \\ 
\exists u\, (u \in \Nc{Y-\{t\}}^+) && \mbox{\qquad namely $u = (Y - \{t\}) \cup \{t\} = Y$} \\ 
\Nc{X}^+ \in \F  && \mbox{\qquad by Lemma~\ref{lemma:successorF}}\\ 
\Nc{Y-\{t\}}^+ \in \F && \mbox{\qquad by Lemma~\ref{lemma:successorF}} 
\end{eqnarray*}
Now we can take the successors:
\begin{eqnarray}
\hskip-1.3cm \Nc{Y-\{t\}}^+ \le \Nc{X}^+ && \mbox{\qquad by Lemma~\ref{lemma:ordersuccessor}} \label{eq:3539} 
\end{eqnarray}
\begin{eqnarray*}
Y \in \Nc{Y}               && \mbox{\qquad by Lemma~\ref{lemma:xinNcx}}  \\  
(Y - \{t\}) \cup \{t\} \in \Nc{Y-\{t\}}^+   && \mbox{\qquad by definition of successor}  \\  
Y \in  \Nc{Y-\{t\}}^+    && \mbox{\qquad since $(Y - \{t\}) \cup \{t\} = Y$}   \\   
\Nc{Y-\{t\}}^+  = \Nc{Y}  &&  \mbox{\qquad by Lemma~\ref{lemma:cardinalsdisjoint}}  \\
\Nc{Y} \le \Nc{X}^+ && \mbox{\qquad by (\ref{eq:3539}) and the preceding line}\\  
\Nc{X \cup \{c\}} = \Nc{X}^+  && \mbox{\qquad by Lemma~\ref{lemma:Ncsuccessor}}\\  
\Nc{Y} \le \Nc{X \cup \{c\}} && \mbox{\qquad since $\Nc{Y} \le \Nc{X}^+  = \Nc{X \cup \{c\}}$} 
\end{eqnarray*}
That completes the induction step. \end{proof}

\begin{lemma} \label{lemma:ontocardinals}  Let $A$ and $B$ be finite sets, and let 
$f$ be a function mapping $A$ onto $B$.  
Then $\Nc{B} \le \Nc{A}.$
\end{lemma}

\noindent\begin{proof}  We may assume without loss of generality that $A$ is the 
domain of $f$.   Then 
\begin{eqnarray*}
f \in \FINITE  && \mbox{\qquad by Lemma~\ref{lemma:finitefunction}}\\
A \times B \in \FINITE  && \mbox{\qquad by Lemma~\ref{lemma:productfinite2}} \\
f \in \FINITE && \mbox{\qquad by Lemma~\ref{lemma:finitefunction}}\\
f \in \SSC(A \times B) && \mbox{\qquad by Lemma~\ref{lemma:finiteseparable}}
\end{eqnarray*}
That is, $f$ is a decidable relation on $A \times B$.
Define
$$Z:= \{s \in B : \exists b \in A\, (\langle b,s \rangle \in f).$$
By Lemma~\ref{lemma:decidable_image}, since $f$ is a decidable relation
on $A \times B$,  $Z$ is a separable subset of $B$.  That is, 
\begin{eqnarray}
\forall t\in B\,(t \in Z \ \lor \ t \not \in Z) \label{eq:3436}
\end{eqnarray}

Now we will proceed by induction on finite sets $A$ to prove that for all finite 
sets $B$ and all $g: A \to B$ onto,  $\Nc{B} \le \Nc{A}$.   
\smallskip

{\em Base case}: If $g: \emptyset \to B$ is onto, then $B = \emptyset$, so $\Nc{A} = \Nc{B} = \Nc{\emptyset}$.
\smallskip

{\em Induction step}: Let $g: A \cup \{c\} \to B$ be onto, where $c \not\in A$.  Let 
$t = g(c)$ and let $f = g - \{\langle c,t\rangle \}$.  Then $f: A \to B$. 
By (\ref{eq:3436}), $t \in Z \ \lor \ t \not\in Z$.  That is, 
$$ \exists b \in A\,(\langle b,t \rangle \in f) \ \lor \ \neg\exists b \in A\,(\langle b,t \rangle \in f).$$
 We may therefore argue by 
these two cases.
\smallskip

{\em Case~1}.  If there 
exists $b \in A$ with $f(b) = t$,  then $f:A\to B$ is onto, so by the induction
hypothesis 
$$\Nc{B} \le \Nc{A} < (\Nc{A})^+ = \Nc{A \cup {c}}$$
as desired. 
\smallskip

{\em Case~2}.  If there does not exist such a $b$ then $f:A \to (B -\{t\})$ is onto.
Also $B - \{t\}$ is a finite set, by Lemma~\ref{lemma:finitedif}.  Hence by the 
induction hypothesis, $\Nc{B - \{t\}} \le \Nc{A}$.  Then
$$ B = (B - \{t\} ) \cup \{ t\}$$
since equality on the finite set $B$ is decidable, so 
\begin{eqnarray*}
\Nc{B} = \Nc{B-\{t\}}^+    && \mbox{} \\
\Nc{B-\{t\}}^+  \le (\Nc{A})^+  && \mbox{\qquad by Lemma~\ref{lemma:ordersuccessor}}\\
\Nc{B}  \le (\Nc{A})^+  && \mbox{\qquad by the previous two lines}\\
\Nc{A}^+ = \Nc{A \cup \{c\}}   && \mbox{\qquad by Lemma~\ref{lemma:Ncsuccessor}, since $c \not\in A$} \\
\Nc{B} \le \Nc{A \cup \{c\}}   && \mbox{\qquad by the previous two lines}
\end{eqnarray*}
That completes Case~2, and that completes the induction step.
\end{proof}

\begin{lemma} \label{lemma:finiteunion2} Let $X$ be a finite set and let $a$ and $b$
be finite subsets of $X$.  Then $a \cup b$ is finite.   
\end{lemma}

\noindent{\em Remark}.  We cannot prove the union of two finite sets is finite without
some additional hypothesis, for consider $\{p\} \cup \{q\}$, where we do not know whether
$p = q$ or not, e.g., $p = \emptyset$ and $q = \{x: x = \{\emptyset\} \ \land\ P\}$, where $P$ 
is Goldbach's conjecture or the Riemann hypothesis.  Does the union contain one or two elements?
\smallskip

\noindent\begin{proof}  We have
\begin{eqnarray*}
a \in \SSC(X)  && \mbox{\qquad by Lemma~\ref{lemma:finiteseparable}} \\
b \in \SSC(X)  && \mbox{\qquad by Lemma~\ref{lemma:finiteseparable}} \\
\forall x \in X\,(x \in a \ \lor \ x \not\in a) && \mbox{\qquad by the definition of $\SSC(X)$}\\
\forall x \in X\,(x \in b \ \lor \ x \not\in b) && \mbox{\qquad by the definition of $\SSC(X)$}\\
\forall x \in X\, (x \in a \cup b \ \lor x \not \in a \cup b) && \mbox{\qquad by the preceding lines and logic}\\
a \cup b \subset X  &&\mbox{\qquad by the definition of $\subseteq$}\\
a \cup b \in \SSC(X) && \mbox{\qquad by the definition of $\SSC(X)$} \\
a \cup b \in \FINITE && \mbox{\qquad by Lemma~\ref{lemma:separablefinite}}
\end{eqnarray*}
\end{proof}

\begin{lemma} \label{lemma:finiteunion3}  Let $X$ be a finite set and let $y$ be a finite
subset of $\SSC(X)$ (that is, the members of $y$ are separable subsets of $X$).  Then the union
of $y$ is a finite set.  That is,
$$ \bigcup y \in \FINITE.$$     
\end{lemma}

\noindent\begin{proof}  By induction on finite sets $y$ (for fixed $X$). 
\smallskip

{\em Base case}.  When $y = \emptyset$, the union of $y$ is also $\emptyset$, which is finite. 
\smallskip

{\em Induction step}.  Suppose $c \not\in y$ and $y \cup \{c\} \subseteq \SSC(X)$.  Then
we have (in a few steps from the definitions of $\bigcup$ and $\cup$)
\begin{eqnarray}
  \bigcup\, (y \cup \{ c\}) = \left(\bigcup y \right)  \cup c  \label{eq:3507}  
\end{eqnarray}
Then 
\begin{eqnarray*}
\bigcup y \in \FINITE  && \mbox{\qquad by the induction hypothesis} \\ 
c \in \SSC(X)  && \mbox{\qquad since  $y \cup \{c\} \subseteq \SSC(X)$}\\ 
c \in \FINITE && \mbox{\qquad by Lemma~\ref{lemma:separablefinite}}\\ 
y \subseteq \SSC(X) && \mbox{\qquad since $y \cup \{c\} \subseteq \SSC(X)$} \\ 
\bigcup y \subseteq X  && \mbox{\qquad since  $y \subseteq \SSC(X)$} \\ 
\bigcup y  \cup c  \in \FINITE && \mbox{\qquad by Lemma~\ref{lemma:finiteunion2}} \\ 
\bigcup\, (y \cup \{c\})\in \FINITE && \mbox{\qquad by (\ref{eq:3507})}   
\end{eqnarray*}
That completes the induction step.
\end{proof}

\section{\texorpdfstring{The initial segments of $\F$}{The initial segments of F}}
Next we begin to investigate the possible cardinalities of finite sets.
The set of integers less than a given integer is a canonical example 
of a finite set.

 \begin{definition}\label{definition:J}
 For $k \in \F$, we define
$$ \J(k) = \{ x \in \F : x < k\}$$
$$ \bar \J(k) = \{x \in \F : x \le k\}.$$
\end{definition}
The definition is stratified, so $\J(k)$ can be defined, but $\J(k)$ 
gets index 1 if $x$ gets index 0, so $\J$ is not definable as a function on $\F$.

\begin{lemma} \label{lemma:Jsuccessor} 
For each $m \in \F$, if $m^+ \in \F$ then 
$$ \J(m^+) = \J(m) \cup \{m\}$$
$$ \bar \J (m^+) = \bar \J(m) \cup \{m^+\}.$$
\end{lemma}

\noindent\begin{proof}  By the definitions of $\J$ and $\bar \J$, and 
the fact that for $x \in \F $ we have $$x < m^+ \iff x < m \ \lor \ x = m,$$
by Lemma~\ref{lemma:lessthansuccessor3}. 
\end{proof}

\begin{lemma} \label{lemma:Jfinite}
For $m \in \F$,  $\J(m)$ and $\bar \J(m)$ are finite sets.
\end{lemma}

\noindent\begin{proof}  By induction on $m$.  The formulas to be proved, namely
$$ \forall m\,(m \in \F \imp \J(m) \in \FINITE)$$
and similarly for $\bar \J$, are stratified, giving $m$ index 0.  $\F$ and $\FINITE$
are parameters and do not require an index.
\smallskip

{\em Base case}, $m = \zero $. Then $\J(\zero) = \emptyset$, by Lemma~\ref{lemma:nothinglessthanzero}.  
By Lemma~\ref{lemma:lambda_finite}, 
$\emptyset \in \FINITE$.  That completes the base case for $\J$.  For $\bar J$,
we have $x \le \zero \iff x = \zero$, so $\bar J(\zero) = \{\zero\}$, which
is finite by Lemma~\ref{lemma:finite_adjoin}.  That completes the base case.
\smallskip

{\em Induction step}.  Suppose $m\in \F$ and $m^+$ is inhabited.  By induction hypothesis,
$\J(m)$ and $\bar \J(m)$ are finite.   By Lemma~\ref{lemma:Jsuccessor},
$\J(m^+) = \J(m) \cup \{m\}$, so by Lemma~\ref{lemma:finite_adjoin}, $J(m^+) \in \FINITE$.
Similarly for $\bar \J(m)$.  That completes the induction step.  
\end{proof}  

\begin{lemma} \label{lemma:Jcardinality}
 Suppose $m \in \F$.  Then 
 $ \Nc{\J(m)} = \T^2 m$. 
\end{lemma}

\noindent\begin{proof}  The formula of the lemma is stratified, giving 
$m$ index 0, since then $\T^2 m$ gets index 2, while $\J(m)$ gets index 1
and $\Nc{\J(m)}$ gets index 2,  so the two sides of the equation both get
index 2.   Therefore the lemma may be proved by induction.
\smallskip

{\em Base case}\,:  $\J(\zero) = \emptyset$, by Lemma~\ref{lemma:xnotlessthanzero}.
We have $\Nc{\emptyset} = \zero$, by Lemma~\ref{lemma:xinNcx} and the definition 
of $\zero$.  By Lemma~\ref{lemma:Tzero},  we have $\T^2 \zero = \zero$.  That completes
the base case.
\smallskip

{\em Induction step}\,: We have
\begin{eqnarray*}
\J(m^+) = \J(m) \cup \{m^+\} && \mbox{\qquad by Lemma~\ref{lemma:Jsuccessor}}\\
\Nc{\J(m)} = \T^2 m  && \mbox{\qquad by the induction hypothesis}\\
\J(m) \in \T^2 m  && \mbox{\qquad by Lemma~\ref{lemma:xinNcx}}\\
\exists u\, (u \in m^+) &&  \mbox{\qquad assumed for proof by induction}\\
m^+ \in \F  && \mbox{\qquad by Lemma~\ref{lemma:successorF}}\\
m \not \in \J(m) && \mbox{\qquad by definition of $\J(m)$}\\
\J(m) \cup \{ m \} \in (\T^2 m)^+ && \mbox{\qquad by definition of successor}\\
(\T  m)^+ = \T (m^+) && \mbox{\qquad by Lemma~\ref{lemma:Tsuccessor}}\\
\T (m^+) \in \F   && \mbox{\qquad by Lemma~\ref{lemma:Tfinite}}\\
(\T m)^+ \in \F  && \mbox{\qquad by the preceding two lines}\\
\exists u\,(u \in (\T m)^+) && \mbox{\qquad by Corollary~\ref{lemma:cardinalsinhabited}}\\
(\T^2 m)^+ = \T^2(m^+) && \mbox{\qquad by Lemma~\ref{lemma:Tsuccessor}}\\
\J(m^+) \in \T^2(m^+)  && \mbox{\qquad by the preceding lines} \\
\J(m^+) \in \Nc{\J(m^+)} && \mbox{\qquad by Lemma~\ref{lemma:xinNcx}}\\
\J(m^+) \in \T^2(m^+) \cap \Nc{\J(m^+)} && \mbox{\qquad by definition of intersection}\\
\Nc{\J(m^+)} = \T^2(m^+) && \mbox{\qquad by Corollary~\ref{lemma:cardinalsinhabited}}
\end{eqnarray*}
That completes the induction step.  
\end{proof}

 \section{Rosser's Counting Axiom}
Rosser introduced the ``counting axiom'',  which is
$$ m \in \F \imp \J(m) \in m.$$
(See~\cite{rosser1978}, p.~485.)  In view of Lemma~\ref{lemma:Jcardinality},
that is equivalent to 
$$ m \in \F \imp \T m = m.$$
Since $2^{\T m}$ is always defined for $m \in \F$,  the counting axiom implies
that $2^m$ is always defined for $m \in \F$.  In particular                                                      then the   
set of iterated powers of $2$  starting from $\zero$  is an infinite set.
That is the conclusion of Specker's proof (but without assuming the counting
axiom).  The point here is that the counting axiom eliminates the need to 
constructivize Specker's proof: if we assume it, there remain only surmountable
difficulties to interpreting HA in \INF.  But the counting axiom is stronger
than NF~\cite{orey1964}, so this observation does not help with the problem of finiteness in \INF.  

\section{Infinity in intuitionistic NF}
We use Dedekind's definition, that a set is infinite if it is similar to a proper subset.
The ``axiom of infinity''  says there is an infinite set.  Before going further,
we remind the reader that with intuitionistic logic, ``not finite''  does not imply ``infinite''.
There are two obvious candidates for infinite sets:  $\V$ and $\F$.  Specker showed that,
with classical logic, $\V$ is not finite; we will discuss that proof below.
\smallskip

If $\F$ is finite, then by Lemma~\ref{lemma:finitemaximal}, there is a maximal finite
cardinal $\m$.  Then by Corollary~\ref{lemma:cardinalsinhabited}, $\m$ has a member $U$,
and by Lemma~\ref{lemma:finitecardinals3}, $U$ is finite.  If we could find some $c \not\in U$,
then $\m^+$ would be inhabited and hence in $\F$, contradicting the maximality of $\m$. 
Therefore $\forall x\,\neg\neg\,(x \in U)$; that is, $\V$ is the double complement of $U$.
However unlikely this may seem, nobody has yet been able to find anything contradictory about it,
without using classical logic.  The following lemma states this remarkable result.

\begin{lemma}\label{lemma:unaargeable2} Suppose $\m$ is a maximal element of $\F$, and $U \in \m$.
Then $\forall x\,\neg\neg\,(x \in U)$.   
\end{lemma}

\begin{lemma} \label{lemma:Tmax}  Let $\m$ is a maximal element of $\F$ and $n \in \F$.  Then
$\T \m < n$ implies $2^n = \emptyset$.  
\end{lemma}

\noindent\begin{proof}  Suppose $\T \m < n$ and $2^n$ is inhabited; we must derive a contradiction.
\begin{eqnarray*}
\USC(u) \in n && \mbox{\qquad for some $u$, by definition of exponentiation}\\
u \in \Nc{u}  && \mbox{\qquad by Lemma~\ref{lemma:xinNcx}}\\
\USC(u) \in \T (\Nc{u})   && \mbox{\qquad by definition of $\T$}\\
\T n = \T(\Nc{u}) && \mbox{\qquad by Lemma~\ref{lemma:cardinalsdisjoint}}\\
\T \m < \T (\Nc{u}) && \mbox{\qquad since $\T \m < n$}\\
\m < \Nc{u}        && \mbox{\qquad by Lemma~\ref{lemma:Tlessthan}}
\end{eqnarray*}
But that contradicts the maximality of $\m$.  \end{proof}

\begin{lemma} If $\V$ is infinite then $\F$ is not finite.
\end{lemma}

\noindent{\em Remark}.  Note that Specker's proof shows $\V$ is not finite, but not that $\V$ 
is infinite, which is stronger.
\smallskip

\noindent\begin{proof}
Suppose $\V$ is infinite and $\F$ is finite, with maximal integer $\m$ and $U \in \m$ and 
$f: \V \to \V$ with $c$ not in the range of $f$.  Then 
\begin{eqnarray*}
\forall x\, \neg\neg\,x \in U   && \mbox{\qquad since $U \in \m$} \\
\forall x\,(x \in U \imp \neg\neg\,(f(x) \in U)) &&\mbox{\qquad by the previous line}\\
\neg\neg\forall x\,(x \in U \imp f(x) \in U) && \mbox{\qquad by Lemma~\ref{lemma:finiteDNS}}\\
\neg\neg\,(f: U \to U)   && \mbox{\qquad by definition of $f:U \to U$}\\
\neg\neg\,(c \in U)  && \mbox{\qquad since $\forall x\, \neg\neg\,x \in U$}
\end{eqnarray*}
That implies that $U$ is not not infinite.  But since $U$ is finite, it is not infinite,
by Theorem~\ref{theorem:infiniteimpliesnotfinite}.   
\end{proof}
\smallskip

\begin{lemma} With classical logic, if $\V$ is not finite then $\F$ is not finite.
\end{lemma}

\noindent\begin{proof} Suppose $\V$ is not finite and $\F$ is finite. Let $\m$ be the 
maximal integer and $U \in \m$.  Then $U$ is finite and $\forall x\, \neg\neg\,(x \in U)$.
Then by classical logic, $\V = U$,  contradiction, since $U$ is finite and $\V$ is not.
\end{proof}
\smallskip

But constructively, the situation is more complicated:  we can prove $\V$ is not finite,
but it is an open problem whether $\F$ is finite or not.

To prove $\F$ is infinite,  we would hope to prove that successor maps $\F$ into $\F$,
so it is of some interest whether that follows from the apparently weaker proposition that
$\F$ is not finite.  We cannot answer that question:  it is an open problem whether
$$ \F \in \FINITE \imp \forall x \in \F\,(x^+ \in \F).$$
In other words, as far as we know,  it might be that $\forall\,U \in \FINITE\,(\V - U \neq \emptyset)$,
but nevertheless we cannot prove $\forall\,U \in \FINITE\, \exists x\,(x \in \V - U)$.  The 
former is equivalent to successor being nonempty on $\F$, the latter to successor being inhabited on $\F$.
We cannot shift the double negation left through $\neg\neg$.  (We shall see below that $\FINITE$ is 
not finite, so Lemma~\ref{lemma:finiteDNS} is no use here.)
\smallskip

Nevertheless, if we did somehow prove  that $\F$ is not finite,  we {\em could} prove that Heyting's 
arithmetic HA is interpretable in $\INF$.  Here is how we would do that:
\smallskip

 Recall that $\F$ is the least set containing $\zero$ and closed under
inhabited successor.  Now define $\H$ to be the least set containing $\zero$ and closed under
nonempty successor.  Then we can prove things using $\H$-induction, in which at the induction
step one is allowed to assume $x^+ \neq \emptyset$, instead of the usual $\exists u\,(u \in x^+)$.
Assume that $\F$ is not finite. 
We do not give all the details, but here is a sketch: First we prove $\F \subseteq \H$,
by $\F$-induction.   Then by $\H$-induction, we prove $\forall x \in \H\,(\neg\neg\,x \in \F)$,
then $\emptyset \not \in \H$; then that $\H$ is closed under successor and has decidable equality,
and that successor is one-to-one on $\H$.  Then we 
could use $\H$ as the interpretation of the variables of HA.  But that would {\em still} not 
prove that $\F$ is closed under successor! 

In the sea of open problems, there is an island:  the theorem of Specker that $\V$ is 
not finite.  This theorem, proved classically in~\cite{specker1953},
 is widely acknowledged as constructively correct,
for reasons I will now explain.  Let $P$ be any stratified formula and let 
$X_P = \{x \in \{\emptyset\} : P\}$.  Then $X_P$  is $\zero$ or $\emptyset$  according as $P$ or $\neg P$.
If $\V$ is finite then $\V$ has decidable equality,  so by deciding whether $X_P =\emptyset$ or not, we 
decide  $P \ \lor \ \neg P$.  That is, $\V$ finite implies the stratified law of excluded middle.
Then,  folklore has it, Specker's proof of infinity
uses classical logic only for stratified formulas,  so it will go through 
under the assumption that $\V$ is finite, and produce a contradiction.  
\smallskip

While this metamathematical argument is appealing, it still requires checking the details of Specker's proof to
ensure that classical logic is used only for stratified formulas.  
I studied Specker's proof, trying to 
make it constructive, and using Lean to check my proofs. 
Assume there is a maximal integer $\m$.  Then 
$\m$ has a member $U$, which is ``unenlargeable'', as discussed above.
I thought that perhaps $U$ could be made to play the role that $\V$ plays in Specker's proof.
That plan did not succeed, unless we assume $\V$ is finite, in 
which case Specker's proof does provide a Lean-checkable proof that $\V$ is 
not finite.   I chose not to present it here.%
\footnote{
It is not very short; the details are in no doubt; it leads to an even lengthier discussion of 
the problem of infinity, but not to a solution of that problem.
}

Rosser, in an appendix to~\cite{rosser1978} (but not the first edition~\cite{rosser1953}), 
gave another proof that $\V$ is not finite, in which 
Specker's ideas are recognizable.  Rosser proves $\V$ is not finite and then immediately
concludes that $\F$ is not finite, since classically $m \in \F$ and $U \in m$, $U$ is finite
so $\V-U$ is inhabited, so $m^+$ is inhabited.  The proof that $\V$ is not finite might well
be constructive.  I did not check it in Lean, since I already checked Specker's proof in Lean.

Once we know that $\V$ is not finite,  we can try to prove other sets are
not finite.  For example, $\FINITE$ is not finite, as we shall prove soon.

\begin{lemma} \label{lemma:uscseparable} $\forall x\,(x \in \FINITE \imp x \in \USC(\V) \ \lor\ x \not\in \USC(\V))$.
\end{lemma}

\noindent\begin{proof} 
 A set $x$ is a singleton if and only if $\Nc{x} = \one$.  That is,
\begin{eqnarray*}
\forall x\,(x \in \USC(\V) \iff \Nc{x} = \one) && \mbox{\qquad by the definitions of $\USC$ and $\one$}   
\end{eqnarray*}
 Since
equality on $\F$ is decidable,  it is decidable whether a finite set is 
a singleton or not.  Therefore
\begin{eqnarray*}
\forall x\, (x \in \FINITE \imp x \in \USC(\V) \ \lor \ x \not\in \USC(\V))&& 
\end{eqnarray*}
\end{proof}

\begin{lemma}\label{lemma:finitenotfinite}  $\FINITE$ is not finite.
\end{lemma}

\noindent{\em Remark}. This depends on the fact that $\V$ is not finite, which 
we do not list as a hypothesis, since it is a theorem, even if the proof has not 
been presented here.
\smallskip


\noindent\begin{proof}  Assume $\FINITE$ is finite.  We must derive
a contradiction. We have
\begin{eqnarray*}
\USC(\V) \subseteq \FINITE && \mbox{\qquad by Lemma~\ref{lemma:singletons_finite}} \\
\USC(\V) \in \SSC(\FINITE)  && \mbox{\qquad by Lemma~\ref{lemma:uscseparable}} \\
\USC(\V) \in \FINITE   && \mbox{\qquad by Lemma~\ref{lemma:separablefinite}, since $\FINITE \in \FINITE$} \\
\V \in \FINITE  && \mbox{\qquad by Lemma~\ref{lemma:uscfinite}} 
\end{eqnarray*}
 \end{proof}
\vspace{-1cm}

\section{Conclusions}
This paper lays the foundations for future studies of intuitionistic \NF\ set theory \INF, 
by providing coherent definitions for the basic concepts, including order,
exponentiation, addition, finite sets, and $\T$.  The concept of separability plays an 
important role in order and power set, and hence in exponentiation as well.  The 
theory presented here---if supplemented by a proof that the set of integers is not 
finite---would serve well as a basis for formalizing constructive mathematics in the style of Bishop.
These basic theorems will surely be both useful and necessary for deeper investigations of the 
metamathematical properties of \INF.  That subject has yet to begin, as at present we cannot
even show that the law of the excluded middle is not provable in \INF.

\printbibliography

\end{document}